\documentclass{amsart}

\usepackage{amsmath,amsthm,amsfonts,amssymb,bm,graphicx, mathrsfs}

\usepackage
{hyperref}
\hypersetup{colorlinks=true,citecolor=blue,linkcolor=blue,urlcolor=blue,
pdfstartview=FitH }

\theoremstyle{plain}
\newtheorem{theorem}{Theorem}
\newtheorem{defi}{Definition}
\newtheorem{lemma}{Lemma}
\newtheorem{sublemma}{Sublemma}

\numberwithin{equation}{section}

\theoremstyle{definition}

\renewcommand{\geq}{\geqslant}
\renewcommand{\leq}{\leqslant}

\newcommand{\abs}[1]{\lvert#1\rvert}

\newcommand{\RR}{\mathbb{R}}

\newcommand{\changed}[1]{{\color{black} #1}}

\usepackage{calc}
\newsavebox\CBox
\newcommand\hcancel[2][0.5pt]{%
  \changed{\ifmmode\sbox\CBox{$#2$}\else\sbox\CBox{#2}\fi%
  \makebox[0pt][l]{\usebox\CBox}%
  \rule[0.5\ht\CBox-#1/2]{\wd\CBox}{#1}}}
\newcommand{\Oneb}{2/561}
\newcommand{\Halfb}{1/561}
\newcommand{\Halfbfrac}{\frac{1}{561}}

\newcommand{\Quarterbfrac}{\frac{1}{1122}}
\newcommand{\Eighthb}{1/2244}
\newcommand{\Eighthbfrac}{\frac{1}{2244}}
\newcommand{\deltaFinal}{1/30000}
\newcommand{\deltaHalf}{\frac{1}{120000}}

\newcommand{\Deltavar}{\changed{\Delta}}
\newcommand{\Xvar}{\changed{X}}

\makeatletter
\DeclareRobustCommand\widecheck[1]{{\mathpalette\@widecheck{#1}}}
\def\@widecheck#1#2{%
    \setbox\z@\hbox{\m@th$#1#2$}%
    \setbox\tw@\hbox{\m@th$#1%
       \widehat{%
          \vrule\@width\z@\@height\ht\z@
          \vrule\@height\z@\@width\wd\z@}$}%
    \dp\tw@-\ht\z@
    \@tempdima\ht\z@ \advance\@tempdima2\ht\tw@ \divide\@tempdima\thr@@
    \setbox\tw@\hbox{%
       \raise\@tempdima\hbox{\scalebox{1}[-1]{\lower\@tempdima\box
\tw@}}}%
    {\ooalign{\box\tw@ \cr \box\z@}}}
\makeatother

\begin{document}

\author{Valentin Blomer}

\author{Jack Buttcane}
\address{Mathematisches Institut, Bunsenstr. 3-5, 37073 G\"ottingen, Germany} \email{blomer@uni-math.gwdg.de}
\email{buttcane@uni-math.gwdg.de}
 
 \title{On the subconvexity problem for $L$-functions on ${\rm GL}(3)$ }

\thanks{First  author   supported  by the  Volkswagen Foundation and a Starting Grant of the European Research Council. Second author supported  by  a Starting Grant of the European Research Council.}


\begin{abstract} Let $f$ be a cusp form for the group $SL(3, \Bbb{Z})$ with Langlands parameter $\mu $ and associated $L$-function $L(s, f)$.  If $\mu$ is in generic position, i.e.\ away from the Weyl chamber walls and away from self-dual forms, we prove the subconvexity bound $L(1/2, f) \ll  \| \mu \|^{\frac{3}{4} - \deltaHalf}$. 
\end{abstract}


\setcounter{tocdepth}{2}  \maketitle 

\maketitle

\section{Introduction}

\subsection{The main result} Analytic number theory on higher rank groups has recently seen substantial advances. One of the most  challenging touchstones for the strength of available techniques is the subconvexity problem for automorphic $L$-functions. We recall that subconvexity refers to an estimate of an automorphic $L$-function on the critical line that is superior (usually with a power saving) to the generic convexity bound in one or more of the defining parameters of the underlying automorphic form. This has been achieved for ${\rm GL}(2)$ in full generality over arbitrary number fields \cite{MV}. In higher rank, the available results become very sporadic. \\

For a fixed, self-dual Maa{\ss} form for ${\rm SL}_3(\Bbb{Z})$, the first breakthrough was achieved by X.\ Li \cite{Li} who solved the subconvexity problem in the $t$-aspect. This was generalized by Munshi \cite{Mu2} to arbitrary fixed Maa{\ss} forms. Similar results are available for  twists 
 by Dirichlet characters \cite{Bl1, Mu1, Mu3}.   All of these results fall into the category of ${\rm GL(1)}$ twists of a fixed Maa{\ss} form and  use mainly ${\rm GL}(1)$ and ${\rm GL(2)}$ tools (enhanced by the ${\rm GL(3)}$ Voronoi formula).  Subconvexity in terms of  genuine parameters of a ${\rm GL(3)}$ automorphic $L$-function (level or spectral parameter) has resisted all attempts so far and seems to require a completely new set of methods. \\

In this paper we go, for the first time,  beyond ${\rm GL}(1)$ twists and prove a prototype of a genuine ${\rm GL}(3)$ subconvexity result  using the spectral theory of automorphic forms on ${\rm GL(3)}$. For  an automorphic representation    $\pi$ we denote its Langlands parameter by $\mu= \mu_{\pi} = (\mu_1, \mu_2, \mu_3)$. This is a triple of complex numbers satisfying $\mu_1 +\mu_2 + \mu_3 = 0$, normalized such that the Ramanujan predicts $\mu \in (i \Bbb{R})^3$. Let $\pi_0$ be an everywhere unramified  automorphic representation  with   Langlands  parameter $\mu_0 = (\mu_{0, 1}, \mu_{0, 2}, \mu_{0, 3})$.  We assume that $\mu_0$ is in \emph{generic} position, i.e.\   there exist  constants $C > c > 0$ such that
\begin{equation}\label{cC}
c \leq \frac{|\mu_{0, j}|}{\| \mu\|} \leq C \quad (1 \leq j \leq 3), \quad \text{and} \quad c \leq \frac{|\mu_{0, i} - \mu_{0, j}|}{\| \mu\|} \leq C \quad (1 \leq i< j \leq 3).
\end{equation}
 This set describes two cones in each Weyl chamber away from the walls and away from the self-dual forms, and covers $99\%$   of all Maa{\ss} forms (choosing $c$ and $C$ appropriately).  For the rest of the paper we fix $c$ and $C$, and all implied constants may depend on them. The convexity bound for $L$-functions associated with such representations states $L(s, \pi_0) \ll   \| \mu_0 \|^{3/4+\varepsilon}.$ 
 
 \begin{theorem}\label{thm1} Let $\pi_0   \subseteq L^2({\rm SL}_3(\Bbb{Z})\backslash \Bbb{H}_3)$ be an irreducible cuspidal representation  with large Langlands  parameter $\mu_0$   in generic position.  Then 
 $$L(1/2, \pi_0 ) \ll  \| \mu_0 \|^{\frac{3}{4} - \deltaHalf}.$$
 \end{theorem}
 
We remark that the same proof works almost literally for any fixed point on the critical line and produces
 $L(1/2+it, \pi_0  ) \ll_{  t} \| \mu_0 \|^{\frac{3}{4} - \deltaHalf}$ with polynomial dependence in $t$. It also works for Maa{\ss} forms for fixed congruence subgroups $\Gamma_0(N) \subseteq {\rm SL}_3(\Bbb{Z})$, again with polynomial dependence on $N$.\\ 
 
 The main tool is the ${\rm GL}(3)$ Kuznetsov formula that was successively refined, most notably in \cite{Bu}, and has recently been used for a variety of applications. The starting point is an amplified fourth moment, averaged over representations with Langlands parameter in an $O(T^{\varepsilon})$-ball   about $\mu_0$. We insert an approximate functional equation and apply Poisson summation in all four variables. It is instructive to compare this with the ${\rm GL(2)}$ version, which was worked out by Iwaniec \cite[Theorem 4]{Iw} more than 20 years ago.  While for ${\rm GL(2)}$  a \emph{second} moment suffices,   in rank 2 a fourth moment is necessary, and the method requires an extremely delicate analysis of Kloosterman sums and special functions. There are several other new phenomena in higher rank that will be discussed in due course. On the technical side, we need very precise estimates for the four-fold Fourier transform of the kernel function of the Kuznetsov transform associated to the long Weyl element.  Ultimately this amounts to the analysis of a multi-dimensional oscillatory integral with degenerate and non-degenerate stationary points to which we apply, among \changed{other} things, Morse theory in the form of a theorem of Milnor and Thom. Several auxiliary results on special functions and integral transforms associated with the group ${\rm GL}_3(\Bbb{R})$ may be useful in other situations. \\ 
  
  
 
 The excluded situations in Theorem \ref{thm1}, i.e.\  forms close to self-dual forms and close to the walls of the Weyl chambers, are exceptional for two different reasons: for self-dual forms the conductor of the $L$-function drops so that instead of a fourth moment a sixth moment would be necessary (Theorem \ref{thm1} remains true in the self-dual case, too, but is worse than the convexity bound). Close to the Weyl chambers, on the other hand, the spectral measure drops, so that the spectral average becomes less powerful.  Notice  that possible  exceptional spectral parameters (i.e.\ violating the Ramanujan conjecture) lie on the Weyl chamber walls, so that these are in particular excluded; this simplifies some of the forthcoming arguments, but is not essential to the method. 
 
 
\subsection{A heuristic roadmap} It might be useful to give a short informal description of the proof which reflects reality -- if at all -- only in a very vague sense, but may guide the reader through the argument. The mean value
$$\sum_{\mu = \mu_0 + O(1)} |L(1/2, \pi)|^4$$
contains about $T^3$ terms, where $T = \| \mu_0 \|$. If we can show that the off-diagonal term is $\ll T^{3 - \delta}$ for some $\delta > 0$, then the amplification method will prove subconvexity. Our amplifier has length $L = T^{\lambda}$ for some very small $\lambda > 0$, but for simplicity we suppress the amplifier in the present discussion. By an approximate functional equation we have
$$|L(1/2, \pi)|^4 \approx  T^{-3} \sum_{m_1, m_2, n_1, n_2 \asymp T^{3/2}}  A_{\pi}(m_2, n_1) \overline{A_{\pi}(m_1, n_2)},$$ 
where here and \changed{throughout the section} we do not display smooth weight functions. The contribution of the long Weyl element of the Kuznetsov formula is roughly of the shape
\begin{equation}\label{heur}
T^{-3} \sum_{m_1, m_2, n_1, n_2 \asymp T^{3/2}}  \sum_{D_1, D_2} \frac{S(n_1, m_2, m_1, n_2; D_1, D_2)}{D_1D_2} \Phi\left(\frac{n_1m_1D_2}{D_1^2}, \frac{n_2m_2D_1}{D_2^2}\right),
\end{equation}
where $S(n_1, m_2, m_1, n_2, D_1, D_2)$ is a certain Kloosterman sum and $\Phi(y) = \Phi(y_1, y_2)$ is given by an  integral of the form
$$\Phi(y) = \int_{\mu = \mu_0 + O(1)} K(y; \mu) \text{spec}(\mu) d\mu,$$
where $\text{spec}(\mu) d\mu \approx \| \mu\|^3 d\mu$ is the spectral measure and $K$ is the kernel function of the ${\rm GL}(3)$ Kuznetsov transform, an analogue of a Bessel $K_{2it}$ or $J_{2it}$ function. The specific shape of this function was clarified recently in \cite{Bu}, and it is given by a double Mellin transform in Definition \ref{MB}. A useful alternative  representation of independent interest as an integral over a product of two Bessel functions is derived in Lemma \ref{stade}. This formula suggests that the typical size of $K(y; \mu)$ is $T^{-3/2}$: each Bessel function saves $T^{1/2}$, and the $u$-integral also saves $T^{1/2}$ by a stationary phase argument. Therefore the typical size of $\Phi(y; \mu)$ is roughly $T^{3/2}$, the square-root of the spectral measure. 

Our first aim is to show that $\Phi(y)$ is much smaller for small arguments, so that we can truncate the $D_1, D_2$-sums. In Lemma \ref{lem4} below we show that the expression becomes negligible for $D_1, D_2 \gg T$. This is not obvious; a direct integration by parts argument would only show $D_1, D_2 \gg T^2$, see \cite{Bl}. 

We need to show now that the $n_1, n_2, m_1, m_2$-sums have (almost) square-root cancellation. 
To this end, we apply Poisson summation in all four variables. It  follows from Lemma \ref{lem4} that the dual variables, say $x_1, x_2, y_1, y_2$, can be truncated at size $\max(D_1, D_2)^{1/2} \ll T^{1/2}$. This is already a step forward since we have shortened the variables, but this alone is not sufficient, and we also need an important Diophantine feature of the Kloosterman sums: 
the Fourier transform of the Kloosterman sum  does not oscillate, but is roughly the characteristic function on $x_1y_1 \equiv D_2$ (mod $D_1$) and $x_2y_2 \equiv D_1$ (mod $D_2$), see Lemma \ref{lem2}. One could now hope that each Fourier integral saves a factor $T^{1/2}$ by a stationary phase argument, so that we get a total saving of $T^2$ and   are left with
$$T^{-3} \cdot  T^{\frac{3}{2} \cdot 4} \sum_{D_1, D_2 \ll T}   \sum_{\substack{x_1, y_1, x_2, y_2 \ll \max(D_1, D_2)^{1/2}\\ x_1y_1 \equiv D_2 \, (\text{mod  }D_1)\\ x_2y_2 \equiv D_1 \, (\text{mod  }D_2)}} \frac{1}{D_1D_2}   \cdot \frac{T^{3/2}}{T^2}.$$
Here we can glue together $z_1 =x_1y_1$, $z_2 = x_2y_2$, and notice that $z_1, z_2 \ll \max(D_1, D_2)$, so that they are essentially fixed by the congruence condition, at least if $D_1, D_2$ are roughly of the same size. This gives a total bound of $T^{5/2}$, and we win. 

The previous discussion is much oversimplified, and real life is   more complex. First, some of the $x$ and $y$ variables can be zero, in which case the divisor argument, implicit in the change of variables $z_j = x_jy_j$, is not possible. Indeed, experience has shown that the central terms in the Poisson summation formula require special care, and in the present situation it turns out that the central Poisson term $x_1 = x_2 = y_1 = y_2 = 0$ is of order $T^3$ and furnishes an additional  off-diagonal main term in the asymptotic formula of the fourth moment. This phenomenon cannot happen with the ${\rm GL}(2)$ Kuznetsov formula, and we refer to the remark after Lemma \ref{lem2} for further discussion.  In particular, our initial hope to prove a bound $O(T^{3-\delta})$ for the off-diagonal term cannot be fulfilled. This off-diagonal main term can be computed explicitly, and it turns out that we can save not in the $T$-aspect, but in the $L$-aspect of the amplifier, which itself is a small power of $T$, see Section \ref{centralsec}.  This is not obvious and follows after non-trivial manipulations from the existence of an accidental zero in the Mellin transform of the Kuznetsov kernel that becomes only apparent after piecing together various terms in the Kuznetsov formula.

Secondly, when $D_1, D_2$ are highly imbalanced,   we need extra savings, since the congruences become less powerful. This is a serious issue and requires a fine-scale analysis of the four-fold Fourier transform of $\Phi(y)$. Finally and most importantly, the desired $T^{1/2}$-savings by stationary phase are very hard to show and do not   happen in general, as there are several degenerate stationary points with smaller savings. This phenomenon can already be seen, for instance, in the one-dimensional case by the function $K_{it}(y)$ in the transitional range $t \approx y$. Much more badly behaved phenomena appear in higher rank, and in addition the stationary points are given by implicit  algebraic expressions that cannot be used for explicit calculations. We must therefore argue  more indirectly. The key result is Lemma \ref{lem6}, where we will   show that the four-fold Fourier transform of $\Phi(y )$ is bounded by $T^{-\Eighthb}$ in typical ranges. This is weaker than  our idealistic (and incorrect) treatment above  with an estimate $T^{3/2}/T^2 = T^{-1/2}$, but  just suffices for a subconvexity estimate.  

\subsection{Notation} Unless noted otherwise, we will use ``$\varepsilon$-convention'': the letter $\varepsilon$ denotes a sufficiently small positive quantity that may change from line to line. There are certain places in the argument, however, where it is important to play off some $\varepsilon$ against another. We will then announce explicitly that  $\varepsilon$-convention will not be in force. For two quantities $A, B$ (positive or negative) we write $A \asymp B$ to mean that there are positive constants $c_1, c_2$ such that $c_1 A \leq B \leq c_2 A$. These constants are absolute and depend only on the constants $c, C$ in \eqref{cC} and the support of the various compactly supported weight functions occurring in the argument. We will sometimes use the phrase ``negligible error'' by which we mean an error term $O_B(T^{-B})$ for an arbitrary   constant $B > 0$. 
 
\section{Preparing the stage} \label{stage}

For $0 \leq c \leq \infty$ let  
$$\Lambda_c := \{\mu \in \Bbb{C}^3 \mid \mu_1 + \mu_2 + \mu_3 = 0,\, |\Re \mu_j| \leq c\}$$
and 
$$\Lambda'_c := \left\{\mu \in \Lambda_c \mid \{-\mu_1, -\mu_2, -\mu_3 \} = \{\overline{\mu_1}, \overline{\mu_2}, \overline{\mu_3}\}\right\}.$$
In the Lie algebra $\mathfrak{a}_{\Bbb{C}}^{\ast} = \Lambda_{\infty}$ we will simultaneously use $\mu$ and $\nu =  (\nu_1, \nu_2, \nu_3)$ coordinates, defined by
 \begin{equation}\label{nu}
 \nu_1 = \frac{1}{3}(\mu_1 - \mu_2), \quad \nu_2 = \frac{1}{3}(\mu_2 - \mu_3), \quad \nu_3 = \frac{1}{3}(\mu_3 - \mu_1).
 \end{equation}
 The latter are already implicit in \eqref{cC}. 
Throughout the paper the letter $\mu =  (\mu_1 , \mu_2 , \mu_3 )$ is reserved for an element in $\mathfrak{a}_{\Bbb{C}}^{\ast}$. 

By unitarity and  the standard Jacquet-Shalika bounds, the Langlands parameter of  an arbitrary irreducible representation $\pi \subseteq L^2({\rm SL}_3(\Bbb{Z})\backslash \Bbb{H}_3)$    is contained in $\Lambda'_{1/2} \subseteq \Lambda_{1/2}$, and the non-exceptional para\-meters are in $\Lambda'_0 = \Lambda_0$. Let $$\mathcal{W} := \left\{I,  w_2 = \left(\begin{smallmatrix}  1 && \\ & & 1\\ & 1 &\end{smallmatrix}\right), w_3 = \left(\begin{smallmatrix}   &1& \\1 & & \\ &  &1\end{smallmatrix}\right),w_4 = \left(\begin{smallmatrix}   &1& \\ & & 1\\ 1&  &\end{smallmatrix}\right), w_5 = \left(\begin{smallmatrix}   && 1\\ 1 & & \\ & 1 &\end{smallmatrix}\right), w_6 = \left(\begin{smallmatrix}   &&1 \\ & 1& \\ 1 &  &\end{smallmatrix}\right)  \right\}$$ be the Weyl group. It acts on $\mu$ by permutation, which defines a corresponding action on $\nu$.
\changed{In particular, the action of the 3-cycles are given by
\begin{equation}\label{3cycle}
	w_4(\mu) = (\mu_3,\mu_1,\mu_2), \qquad w_5(\mu) = (\mu_2,\mu_3,\mu_1).
\end{equation}}


Let $\pi_0 \subseteq L_{\text{cusp}}^2({\rm SL}_3(\Bbb{Z})\backslash \Bbb{H}_3)$ be our preferred irreducible cuspidal automorphic representation with Hecke eigenvalues $A_{\pi_0}(1, n)$ and Langlands parameter $\mu_{0} = (\mu_{0, 1}, \mu_{0, 2}, \mu_{0, 3}) \in \Lambda_0$, and assume   that   
\begin{equation}\label{munu}
  |\mu_{0, j}| \asymp |\nu_{0, j}| \asymp T \quad (j = 1, 2, 3)
\end{equation}
 for some sufficiently large parameter $T$. (As mentioned in the introduction, this implies in particular that $\mu_0 \in \Lambda_0$.)   Recall that in general $A_{\pi}(1, m) = \overline{A_{\pi}(m, 1)}$, see e.g. \cite[p.\ 230]{Go}.   By a standard approximate functional equation (\cite[Section 5.2]{IK}) we have 
 $$L(1/2, \pi_0   ) = \sum_n \frac{A_{\pi_0}(1, n) }{\sqrt{n}} V\left(\frac{n}{T^{3/2}}\right) + \kappa  \overline{\sum_n \frac{A_{\pi_0}(1, n) }{\sqrt{n}}  V\left(\frac{n}{T^{3/2}}\right)},$$
 where $|\kappa| = 1$ and $V$ is a smooth function satisfying the uniform bounds
 $$x^j V^{(j)}(x) \ll_B (1 +x)^{-B}.$$
  Inserting a smooth partition of unity, this shows
$$L(1/2, \pi_0) \ll   \sum_{2^j \leq T^{3/2+\varepsilon}}2^{-j/2}  \Bigl|\sum_n  A_{\pi_0}(1, n)   W\left(\frac{n}{2^j}\right) \Bigr|  $$
(up to a negligible error) for some fixed, smooth, compactly supported function $W$. 
 Using the Hecke relation (\cite[Section 6.4]{Go})
$$A_{\pi_0}(1, n) \overline{A_{\pi_0}(1, m)} = \sum_{d\mid (n, m)} A_{\pi_0}\left(\frac{m}{d}, \frac{n}{d}\right)$$
and the Cauchy-Schwarz inequality, we obtain
$$|L(1/2, \pi_0)|^2 \ll  (\log T) \sum_{2^j \leq T^{3/2+\varepsilon}} 2^{-j} \sum_d \sum_{n, m}  A_{\pi_0}(m, n)  W\left(\frac{nd}{2^j}\right)\overline{W\left(\frac{md}{2^j}\right)}.$$
For $M \gg 1$ and an arbitrary    $\pi \subseteq L^2({\rm SL}_3(\Bbb{Z})\backslash \Bbb{H}_3)$  (potentially generated by an Eisenstein series) let
$$\mathcal{L}_M(\pi) := \frac{1}{M} \Bigl|\sum_{n, m}  A_{\pi}(m, n)   W\left(\frac{n}{M}\right)\overline{W\left(\frac{m}{M}\right)}\Bigr|.$$
Then clearly 
\begin{equation}\label{basic1}
|L(1/2, \pi_0)|^2 \ll T^{\varepsilon} \max_{M \leq T^{3/2+\varepsilon}} \mathcal{L}_M(\pi_0). 
\end{equation}
It follows from \cite[Theorem 2]{Li1} or \cite[Corollary 2]{Br}  that
$$\sum_{m \leq x} |A_{\pi_0}(m, 1)|^2 \ll x(xT)^{\varepsilon},$$
which together with the Hecke relations and the Cauchy-Schwarz inequality easily implies the trivial bound 
\begin{equation}\label{etasmall}
  \mathcal{L}_M(\pi_0) \ll M (MT)^{\varepsilon}.
\end{equation}
We will use this bound if $M$ is small. We   fix some small $0 < \eta < 1$ and assume from now on 
\begin{equation}\label{sizeM}
  T^{3/2 - \eta} \leq M \leq T^{3/2+\varepsilon}.
\end{equation}\\
 
For any     $\pi \subseteq L^2({\rm SL}_3(\Bbb{Z})\backslash \Bbb{H}_3)$   
  we have the Hecke relation (\cite[Section 6.4]{Go})
$$A_{\pi}(1, \ell) A_{\pi}(1, \ell^2) = A_{\pi}(1, \ell^3) + A_{\pi}(1, \ell)A_{\pi}(\ell, 1) - 1,$$
from which we conclude
$$\max(|A_{\pi}(1, \ell)|, |A_{\pi}(1, \ell^2)|, |A_{\pi}(1, \ell^3)|) \geq 1/2.$$
This allows us to construct an amplifier.  
Let $x(n) := \text{sgn}(A_{\pi_0}(1, n)) \in S^1 \cup \{0\}$. Fix some sufficiently small $0 < \lambda < 1/20$, and for
\begin{equation}\label{lambda}
L = T^{\lambda} \leq T^{1/20}
\end{equation}
 define
$$\mathcal{A}(\pi) = \sum_{j=1}^3 \Bigl| \sum_{\substack{L \leq \ell \leq 2L\\ \ell \text{ prime}}} A_{\pi}(1, \ell^j) \overline{x(\ell^j)}\Bigr|^2.$$
Clearly $\mathcal{A}(\pi_0) \geq \frac{1}{2} |\{L \leq \ell \leq 2L \mid \ell \text{ prime}\}|^2 \gg L^{2-\varepsilon}$. 

Let $h$ be a non-negative function on 
$ \Lambda'_{1/2}$ that is rapidly decaying as $|\Im \mu_j| \rightarrow \infty$ for $j=1, 2, 3$ and  satisfies  $h(\mu_0) \gg 1$. Let $\mathcal{N}(\pi)$ be some positive quantities (they will later be some normalizing factors) such that   $\mathcal{N}(\pi) \ll \| \mu_{\pi}\|^{\varepsilon}$ for cuspidal $\pi$.  Then clearly
$$\mathcal{L}_M(\pi_0)^2  \ll T^{\varepsilon} \frac{\mathcal{A}(\pi_0)}{L^2}  \mathcal{L}_M(\pi_0)^2 \ll \frac{T^{\varepsilon}}{L^2}   \int  \mathcal{A}(\pi)  \mathcal{L}_M(\pi)^2  \frac{h(\mu_{\pi})}{\mathcal{N}(\pi)} d\pi,$$
where here and in the following the notation $\int (\cdots) d\pi$ is understood as a combined sum/integral over an orthonormal basis of spectral components of $L^2({\rm SL}_3(\Bbb{Z})\backslash \Bbb{H}_3)$,  which effectively runs over Hecke-Maa{\ss} cusp forms and Eisenstein series. The  precise shape of the spectral decomposition is given explicitly, for instance, in \cite[Theorem 4]{Bu}.   We have
\begin{displaymath}
\begin{split}
\mathcal{A}(\pi)  \mathcal{L}_M(\pi)^2   & = \frac{1}{M^2}\sum_{j=1}^3\sum_{\substack{\ell_1, \ell_2\asymp L\\ \ell_1, \ell_2 \text{ prime}}}\overline{x(\ell_1^j)}  x(\ell_2^j)  \sum_{n_1, n_2, m_1, m_2}  W\Bigl(\frac{ n_1}{M}\Bigr) W\Bigl(\frac{m_1}{M}\Bigr)  \overline{W\Bigl(\frac{ n_2}{M}\Bigr)W \Bigl(\frac{m_2}{M} \Bigr) }   \\
& \quad\quad\quad\quad\quad   \times A_{\pi}(m_2, n_1) \overline{A_{\pi}(m_1, n_2)}A_{\pi}(1, \ell_1^j)\overline{A_{\pi}(1, \ell_2^j)} .
\end{split}
\end{displaymath}
By the Hecke relations (\cite[Section 6.4]{Go}), the second line equals
$$\sum_{\substack{r_0r_1r_2 = \ell_1^j\\r_0\mid m_2, r_2 \mid n_1 }} \sum_{\substack{s_0s_1s_2 = \ell_2^j\\ s_0 \mid m_1, s_1 \mid n_2}}  A_{\pi}(m_2r_2/r_0, n_1r_1/r_2) \overline{A_{\pi}(m_1s_1/s_0, n_2s_2/s_1)}, $$
so that 
\begin{equation}\label{basic}
\begin{split}
\mathcal{L}_M(\pi_0)^2  \ll  \frac{T^{\varepsilon}}{M^2L^2} & \sum_{j=1}^3
 \sum_{\substack{\ell_1, \ell_2 \asymp L\\ \ell_1, \ell_2 \text{ prime}}}  \sum_{\substack{r_0r_1r_2 = \ell_1^j \\ s_0s_1s_2 = \ell_2^j }}   \Bigl|  \sum_{n_1, n_2, m_1, m_2} W\Bigl(\frac{r_2n_1}{M}\Bigr) W\Bigl(\frac{s_0m_1}{M}\Bigr)   \\
 & \times \overline{W\Bigl(\frac{s_1n_2}{M}\Bigr)W \Bigl(\frac{r_0m_2}{M} \Bigr) }  \int  A_{\pi}(m_2r_2, n_1r_1) \overline{A_{\pi}(m_1s_1, n_2s_2)} \frac{h(\mu_{\pi})}{\mathcal{N}(\pi)} d\pi \Bigr|. 
 \end{split}
\end{equation}
This is our basic inequality to which we can apply the Kuznetsov formula.

\section{The Kuznetsov formula}

\subsection{Normalizing factors} In this subsection we choose the normalizing factors $\mathcal{N}(\pi)$ as the (square of the) ratio between Hecke eigenvalues and Fourier coefficients of $L^2$-normalized automorphic forms appearing in the spectral decomposition. An inspection of \cite[Theorem 4]{Bu} shows that 
 for a cuspidal automorphic representation $\pi \subseteq  L^2({\rm SL}_3(\Bbb{Z})\backslash \Bbb{H}_3)$ we need to define 
$$\mathcal{N}(\pi) := \| \phi \|^2  \prod_{j=1}^3  \cos\left(\frac{3}{2} \pi \nu_{\pi, j} \right),  $$ where $\phi$ is the arithmetically normalized Maa{\ss} form $\phi$ generating $\pi$ and $\nu$ is given by \eqref{nu}. That is, $\phi$ is given by the Fourier expansion
$$\phi(z) = \sum_{\gamma    \in U \backslash {\rm SL}_2(\Bbb{Z})}  \sum_{m_1=1}^{\infty} \sum_{ m_2 \not= 0} \frac{A_{\pi}(m_1, m_2)}{|m_1m_2|}  \mathscr{W}^{\text{sgn}(m_2)}_{\nu}\left(\left(\begin{smallmatrix} |m_1m_2| & & \\ & m_1 & \\ & & 1 \end{smallmatrix}\right)\left(\begin{smallmatrix} \gamma  & \\ & 1 \end{smallmatrix}\right)z\right), $$
where $U = \left\{\left(\begin{smallmatrix} 1 & *\\ & 1 \end{smallmatrix}\right)\in {\rm SL}_2(\Bbb{Z})\right\}$ and $\mathscr{W}^{\pm}_{\nu}(z) = e(x_1 \pm x_2) W^{\ast}_{\nu}(y_1, y_2)$, where $W^{\ast}_{\nu}$ is the standard \emph{completed} Whittaker function as in \cite[Def.\ 5.9.2]{Go}, and $A_{\pi}(1, 1) = 1$.  
By Rankin-Selberg theory  in combination with Stade's formula (see e.g.\ \cite[Section 4]{Bl}) and \cite[Theorem 2]{Li1}, it is easy to see that
$$\mathcal{N}(\pi) \asymp  \underset{s=1}{\text{res}} L(s, \pi \times \tilde{\pi}) \ll \| \mu_{\pi} \|^{\varepsilon}$$
with  implied constants depending at most on $\varepsilon$. 
For non-cuspidal $\pi$, one can check that the proper analogue of $\mathcal{N}(\pi)$ is given by $$
\frac{1}{16} \prod_{j=1}^3 |\zeta(1 + 3\nu_{\pi, j})|^2$$
if $\pi$ is generated by a minimal Eisenstein series (see \cite[Chapter 7]{Bum}), and
$$ 8 L(1, \text{Ad}^2 u) |L(1 + 3s, u)|^2$$
if $\pi$ is generated by a maximal Eisenstein series $E(z, 1/2 + s, u)$ associated to an ${\rm SL}_2(\Bbb{Z})$ cusp form $u$, although this plays no role in our situation. 

\subsection{Kloosterman sums}
For $n_1, n_2, m_1, m_2, D_1, D_2 \in \Bbb{N}$  we define the two relevant types of  Kloosterman sums  by
\begin{displaymath}
	\tilde{S}(n_1,n_2,m_1;D_1,D_2) := \sum_{\substack{C_1 (\text{mod }D_1), C_2 (\text{mod }D_2)\\(C_1,D_1)=(C_2,D_2/D_1)=1}} e\left(n_2\frac{\bar{C_1}C_2}{D_1}+m_1\frac{\bar{C_2}}{D_2/D_1}+n_1\frac{C_1}{D_1}\right)
\end{displaymath}
for $D_1\mid D_2$,   and
\begin{displaymath}
\begin{split}
&S(n_1, m_2, m_1, n_2; D_1, D_2)\\
& = \sum_{\substack{B_1, C_1 \, ({\rm mod }\, D_1)\\B_2, C_2 \, ({\rm mod }\,  D_2)\\ D_1C_2 + B_1B_2 + D_2C_1 \equiv 0 \, ({\rm mod }\, D_1D_2)\\ (B_j, C_j, D_j) = 1}} e\left(\frac{n_1B_1 + m_1(Y_1 D_2 - Z_1 B_2)}{D_1} + \frac{m_2B_2 + n_2(Y_2 D_1 - Z_2B_1)}{D_2}\right),
\end{split}
\end{displaymath}
where $Y_jB_j + Z_jC_j \equiv 1 \, (\text{mod }D_j)$ for $j = 1, 2$. We have the standard (Weil-type) bounds
\begin{equation}\label{tildeS}
\tilde{S}(n_1, n_2,m_1; D_1, D_2) \ll \left((m_1, D_2/D_1) D_1^2, (n_1, n_2, D_1)D_2\right)(D_1D_2)^{\varepsilon}
\end{equation}
and
\begin{equation}\label{S}
S(n_1, m_2, m_1, n_2; D_1, D_2) \ll (D_1D_2)^{1/2+\varepsilon} \left((D_1, D_2)(m_1n_1, [D_1, D_2])(m_2n_2, [D_1, D_2])\right)^{1/2}.
\end{equation}
The first bound is due to Larsen \cite[Appendix]{BFG}, the second due to Stevens (see  \changed{\cite[p.\ 383]{Bu1}}). 

\subsection{Integral kernels}\label{intrep} Following \cite[Theorem 2 \& 3]{Bu}, we define the following  integral kernels in terms of Mellin-Barnes representations. For $s \in \Bbb{C}$, $\mu \in \Lambda_{\infty}$ define the meromorphic function
$$\tilde{G}^{\pm}(s, \mu) := \frac{  \pi^{-3s}}{12288 \pi^{7/2}}\Biggl(\prod_{j=1}^3\frac{\Gamma(\frac{1}{2}(s-\mu_j))}{\Gamma(\frac{1}{2}(1-s+\mu_j))} \pm i   \prod_{j=1}^3\frac{\Gamma(\frac{1}{2}(1+s-\mu_j))}{\Gamma(\frac{1}{2}(2-s+\mu_j))} \Biggr), $$
and for $s = (s_1, s_2) \in \Bbb{C}^2$, $\mu \in \Lambda_{\infty}$ define the meromorphic function
$$G(s, \mu) :=  \frac{1}{\Gamma(s_1 + s_2)} \prod_{j=1}^3 \Gamma(s_1 - \mu_j) \Gamma(s_2 + \mu_j).$$
The latter is essentially the double Mellin transform of the ${\rm GL(3)}$ Whittaker function. We also define the following trigonometric functions
\begin{displaymath}
\begin{split}
& S^{++}(s, \mu) := \frac{1}{24 \pi^2} \prod_{j=1}^3 \cos\left(\frac{3}{2} \pi \nu_j\right),\\
&  S^{+-}(s, \mu) :=  -\frac{1}{32 \pi^2} \frac{\cos(\frac{3}{2} \pi \nu_2)\sin(\pi(s_1 - \mu_1))\sin(\pi(s_2 + \mu_2))\sin(\pi(s_2 + \mu_3))}{\sin(\frac{3}{2} \pi \nu_1)\sin(\frac{3}{2} \pi \nu_3) \sin(\pi(s_1+s_2))}, \\
& S^{-+}(s, \mu) :=-\frac{1}{32 \pi^2}  \frac{\cos(\frac{3}{2} \pi \nu_1)\sin(\pi(s_1 - \mu_1))\sin(\pi(s_1 - \mu_2))\sin(\pi(s_2 + \mu_3))}{\sin(\frac{3}{2} \pi \nu_2)\sin(\frac{3}{2} \pi \nu_3)\sin(\pi(s_1+s_2))}, \\
& S^{--}(s, \mu) := \frac{1}{32 \pi^2}  \frac{\cos(\frac{3}{2} \pi \nu_3) \sin(\pi(s_1 - \mu_2))\sin(\pi(s_2 + \mu_2))}{\sin(\frac{3}{2} \pi \nu_2)\sin(\frac{3}{2} \pi \nu_1)}. 
  \end{split}
\end{displaymath}

\begin{defi}\label{MB} For  $y \in \Bbb{R} \setminus \{0\}$ with $\text{{\rm sgn}}(y) = \epsilon$ 
let
$$K_{w_4}(y; \mu) =  \int_{-i\infty}^{i\infty} | y|^{-s} \tilde{G}^{\epsilon}(y)(s, \mu) \frac{ds}{2\pi i}. $$
For $y = (y_1, y_2) \in (\Bbb{R}\setminus \{0\})^2$ with $\text{{\rm sgn}}(y_1) = \epsilon_1$, $\text{{\rm sgn}}(y_2) = \epsilon_2$ let
\begin{displaymath}
\begin{split}
 K^{\epsilon_1, \epsilon_2}_{w_6}(y; \mu)  =    & \int_{-i\infty}^{i\infty} \int_{-i\infty}^{i\infty}  |4\pi^2 y_1|^{-s_1} |4\pi^2 y_2|^{-s_2}  G(s, \mu) S^{\epsilon_1, \epsilon_2}(s, \mu)\frac{ds_1\, ds_2}{(2\pi i)^2}. 
\end{split}
\end{displaymath}
 We generally follow the Barnes integral convention that the contour should pass to the right of all of the poles of the Gamma functions in the form $\Gamma(s_j + a)$  and to the left of all of the poles of the Gamma functions\footnote{Such Gamma functions may occur through the functional equation $ (\sin(\pi(s_1 + s_2))\Gamma(s_1 + s_2) )^{-1}  =  \frac{1}{\pi} \Gamma(1 - s_1 - s_2)$.} in the form $\Gamma(a - s_j)$. 
 Moreover, we choose the contour such that all integrals are absolutely convergent, which can always be arranged by shifting the unbounded part appropriately. \end{defi}

To substantiate this last claim we observe that the integral kernels have no exponential increase in any of the variables. This is obvious for $K_{w_4}$, and in the case of $K^{\epsilon_1, \epsilon_2}_{w_6}$ the exponential behaviour is given by $\exp(-\frac{\pi}{2} h^{\epsilon_1, \epsilon_2}(\Im s, \Im \mu))$, where $h^{\epsilon_1, \epsilon_2}(t, r)$ is the non-negative function
\begin{equation}\label{exp}
\begin{split}
h^{\epsilon_1, \epsilon_2}(t, r) = & -\epsilon_2 |r_1-r_2| - \epsilon_1\epsilon_2|r_1 - r_3| - \epsilon_1|r_2 - r_3| - \epsilon_1\epsilon_2|t_1+t_2| \\
&+ \epsilon_1\epsilon_2 |t_1 - r_1| + \epsilon_1 |t_1 - r_2 |  + |t_1 - r_3| + |t_2 + r_1| + \epsilon_2 |t_2 + r_2| + \epsilon_1\epsilon_2 |t_2 + r_3|. 
\end{split}
\end{equation}
If $\mu \in \Lambda_0$ for instance, then the unbounded part of the integral  for $K_{w_4}$ must satisfy $\Re s \leq 1/6-\delta$, while the unbounded part of the integral  for $K^{\pm, \pm}_{w_6}$ must satisfy $\Re s_1, \Re s_2 \leq -\delta$ for some $\delta > 0$.

\subsection{The Kuznetsov formula}\label{kuz} We define the spectral measure by
\begin{equation*}
\text{spec}(\mu) d\mu,  \quad  \text{spec}(\mu) := \prod_{j=1}^3 \left(3\nu_j \tan\left(\frac{3\pi}{2} \nu_j\right)\right),
\end{equation*}
where $d\mu = d\mu_1 d\mu_2= d\mu_1d\mu_3 = d\mu_2 d\mu_3$ is the standard measure on the hyperplane $\mu_1 + \mu_2 + \mu_3 = 0$. 
 
We can now state the Kuznetsov formula in the version of \cite[Theorems 2, 3, 4]{Bu}. Let  $n_1$, $n_2$, $m_1$, $m_2 \in \Bbb{N}$ and let $h$ be a function  that is holomorphic on $\Lambda_{1/2 + \delta}$ for some $\delta > 0$, symmetric under the Weyl group, rapidly decaying as $|\Im \mu_j| \rightarrow \infty$  and satisfies 
\begin{equation}\label{zeros}
  h(3\nu_j \pm 1)  = 0, \quad j = 1, 2, 3.
\end{equation}
Then we have 
\begin{displaymath}
\begin{split}
&  \int \overline{A_{\pi}(m_1, m_2)} A_{\pi}(n_1, n_2) \frac{h(\nu_{\pi}) }{\mathcal{N}(\pi)} d\pi = \Delta + \Sigma_4 + \Sigma_5 + \Sigma_6,
\end{split}
\end{displaymath}
where 
\begin{displaymath}
\begin{split}
  \Delta& = \delta_{n_1, m_1} \delta_{n_2, m_2}  \frac{1}{192\pi^5} \int_{\Re \mu = 0} h(\mu) \text{spec}(\mu) d \mu,\\
   \Sigma_{4}& = \sum_{\epsilon  = \pm 1} \sum_{\substack{D_2 \mid D_1\\  m_2 D_1= n_1 D_2^2}}\frac{ \tilde{S}(-\epsilon n_2, m_2, m_1; D_2, D_1)}{D_1D_2} \Phi_{w_4}\left(  \frac{\epsilon m_1m_2n_2}{D_1 D_2} \right),  \\ 
 \Sigma_{5} &= \sum_{\epsilon  = \pm 1} \sum_{\substack{     D_1 \mid D_2\\ m_1 D_2 = n_2 D_1^2}} \frac{ \tilde{S}(\epsilon n_1, m_1, m_2; D_1, D_2) }{D_1D_2}\Phi_{w_5}\left( \frac{\epsilon n_1m_1m_2}{D_1 D_2}\right),\\
   \Sigma_6 &= \sum_{\epsilon_1, \epsilon_2 = \pm 1} \sum_{D_1,  D_2  } \frac{S(\epsilon_2 n_2, \epsilon_1 n_1, m_1, m_2; D_1, D_2)}{D_1D_2} \Phi_{w_6}  \left( - \frac{\epsilon_2 m_1n_2D_2}{D_1^2}, - \frac{\epsilon_1 m_2n_1D_1}{ D_2^2}\right)
\end{split}
\end{displaymath} 
and
\begin{equation}\label{defPhi}
\begin{split}
& \Phi_{w_4}(y) =  \int_{\Re \mu = 0} h(\mu) K_{w_4}(y; \mu )\, \text{spec}(\mu) d \mu,\\
& \Phi_{w_5}(y) = \int_{\Re \mu = 0} h(\mu) K_{w_4}(-y; -\mu )\, \text{spec}(\mu) d \mu ,\\
& \Phi_{w_6}(y_1, y_2) = \int_{\Re \mu = 0} h(\mu) K^{\text{sgn}(y_1), \text{sgn}(y_2)}_{w_6}((y_1, y_2) ;  \mu )\, \text{spec}(\mu) d \mu.
\end{split}
\end{equation}

\subsection{Choice of test function}\label{hh} We now specify a test function $h$ (depending on $\mu_0$) that satisfies the properties required for the Kuznetsov formula, is non-negative on $\Lambda'_{1/2}$, satisfies $h(\mu_0) \gg 1$ and is negligibly small outside    $O(T^{\varepsilon})$-balls about $w(\mu_0)$ for $w\in \mathcal{W}$. To this end let $\psi$ be a fixed holomorphic function on $\Lambda_{\infty}$ that is non-negative,    rapidly decaying as $|\Im \mu_j| \rightarrow \infty$ and bounded from below at the origin;  we choose  $$\psi(\mu) = \exp\left(\mu_1^2 +\mu_2^2  + \mu_3^2\right).$$ 
Let
$$P(\mu) :=  \prod_{0 \leq n \leq A} \prod_{j=1}^3 \frac{(\nu_j - \frac{1}{3}(1 + 2n))(\nu_j + \frac{1}{3}(1 + 2n))}{|\nu_{0, j}|^2} $$
for some large fixed constant $A$. This polynomial has zeros at the poles of the spectral measure,  which turns out to be convenient for later contour shifts. 
Now we  choose
\begin{equation}\label{defh}
h(\mu) :=  P(\mu)^2 \Bigl(\sum_{w \in \mathcal{W}}\psi\Bigl(\frac{w(\mu)  -  \mu_0}{T^{\varepsilon}}\Bigr)\Bigr)^2
\end{equation}
for some very small $\varepsilon$. 
This function localizes at a ball of radius $T^{\varepsilon}$ about $w(\mu_0)$ for each $w \in \mathcal{W}$.  The $T^{\varepsilon}$-radius gives us a bit of elbow room that is convenient in later estimations. In particular, we have
\begin{equation}\label{partial}
\mathscr{D}_j h(\mu) \ll_j T^{-j \varepsilon}
\end{equation}
for any differential operator $\mathscr{D}_j$ of   order $j$, which we will use frequently when we integrate by parts, as sufficiently many differentiations can save arbitrarily many powers of $T$. Moreover, we have  (not applying $\varepsilon$-convention)
\begin{equation}\label{vol}
  \int_{\Re \mu = 0}  h(\mu) \text{spec}(\mu) d\mu \ll T^{3 + 2\varepsilon}. 
\end{equation}

By construction, $h$ is symmetric, holomorphic, rapidly decaying as $|\Im \mu_j| \rightarrow \infty$ and satisfies \eqref{zeros}.  Since $\overline{\psi(\mu)} = \psi(-\bar{\mu})$ and $\overline{P(\mu)} = P(-\bar{\mu})$, it follows from unitarity that $h$ is the square of a real number for $\mu \in \Lambda'_{1/2}$, so that $h(\mu) \geq 0$ for $\mu \in \Lambda'_{1/2}$. Finally it is clear that $h(\mu_0) \gg 1$.

\subsection{Absolute convergence}\label{absconv} We will show now that holomorphicity of $h$ on $\Lambda_{1/2 + \delta}$ together with the vanishing condition \eqref{zeros} yields the general bounds  $\Phi_{w_4}(y)$, $\Phi_{w_5}(y) \ll |y|^{1/10}$ (say) and $\Phi_{w_6}(y) \ll |y_1y_2|^{1/2 + \delta'}$ for $0 < \delta' < \min(1/2, \delta)$. Together with \eqref{tildeS} and \eqref{S} this implies that  the Kloosterman terms  $\Sigma_4, \Sigma_5, \Sigma_6$ are absolutely convergent.   In fact, the following lemma shows   more quantitatively, that with our particular choice of $h$ we can truncate   the $D_1, D_2$-sums 
   at $D_1, D_2 \ll T^{C}$ for some sufficiently large $C$ at the cost of an error $O(T^{-1000})$, provided that $n_1, n_2, m_1, m_2 \ll T^{10}$, say.  
   
\begin{lemma}\label{lemma1} With $h$  as in \eqref{defh} we have
$$\Phi_{w_4}(y) \ll |y|^{1/10} T^{O(1)}, \quad \Phi_{w_5}(y) \ll |y|^{1/10} T^{O(1)}, \quad \Phi_{w_6}(y) \ll |y_1y_2|^{3/5} T^{O(1)}.$$
\end{lemma}
   
 \textbf{Proof.} We use the Mellin-Barnes representations of Definition \ref{MB}.  
  

 In the integral of $K_{w_4}(y; \mu)$   we shift the contour to $\Re s = -1/10$. The remaining integral satisfies the crude bound $\ll |y|^{1/10} \| \mu \|^{O(1)}$. We pick up poles at $s = \mu_j$. (By marginal contour shifts in $\mu$, e.g.\ $\Re \mu_1 = -\varepsilon$, $\Re \mu_2 = 0$, $\Re \mu_3 = \varepsilon$, we can make sure that none of these poles coincide.) For each of the residues we shift the $\mu_j$-contour in \eqref{defPhi} to $\Re \mu_j = -1/10$. This crosses no poles, since the spectral measure $\text{spec}(\mu)$ vanishes at $\mu_i = \mu_j$ for $i \not= j$ and  cancels the poles of the Gamma factors. 
 In this way we obtain $\Phi_{w_4}(y) \ll |y|^{1/10} T^{O(1)}.$ The same bound  holds for $\Phi_{w_5}$.\\  

In the  integral of $K^{\epsilon_1, \epsilon_2}_{w_6}(y; \mu)$ we shift the contour to  $\Re s_1, \Re s_2 = -3/5$. The remaining integral satisfies the crude bound $\ll |y_1y_2|^{3/5}  \| \mu \|^{O(1)}$. There are now two sources of (possible) poles.

a) There are  pure residues at $(s_1, s_2) = (\mu_i, -\mu_j)$ for $i \not= j$. It is easy to see that these have at most simple poles at $\mu_{\ell} - \mu_k \in \Bbb{Z}$ for $\ell \not = k$. Here we shift $\Re \mu_i= -3/5$ and $\Re \mu_j = 3/5$. Again this crosses no poles, but this time this requires in addition to the vanishing of the spectral measure at $\mu_i - \mu_j =0$ also the vanishing of   $h$ at $\mu_i - \mu_j  = \pm 1$ since $6/5 > 1$. It is at this point where \eqref{zeros} is needed.  

b) There are  mixed terms at $s_1 = \mu_j$ and $\Re s_2 = -3/5$  (and the same with exchanged indices). It is easy to see that these have at most simple poles at $\mu_{\ell} - \mu_k \in \Bbb{Z}$ for $\ell \not = k$ and $\mu_{\ell} +  s_2 \in \Bbb{Z}$ for $\ell \not= j$. Here we shift  $\Re\mu_j = -3/5$, and the other two $\mu$-coordinates go to real part $3/10$.  By the properties of the spectral measure, this crosses no poles. 

In all cases we obtain the bound $\Phi_{w_6}(y) \ll |y_1y_2|^{3/5} T^{O(1)}.$ \hfill $\square$  


\section{Analytic preliminaries}

In this section we compile various auxiliary analytic results for future reference.

\subsection{Oscillatory integrals} We will frequently show that oscillatory integrals are very small using integration by parts. For convenience we quote here a useful lemma from \cite{BKY} that can be applied in all situations.

\begin{lemma} \label{integrationbyparts}
 Let $Y \geq 1$, $X, Q, U, R > 0$, 
and suppose that $w$ 
 is a smooth function with support on some interval $[\alpha, \beta]$, satisfying
\begin{equation*} 
w^{(j)}(t) \ll_j X U^{-j}.
\end{equation*}
Suppose $H$ 
  is a smooth function on $[\alpha, \beta]$ such that
\begin{equation}\label{condi1}
 |H'(t)| \gg R, \quad 
H^{(j)}(t) \ll_j Y Q^{-j} \,\, \text{for } j=2, 3, \dots.
\end{equation}
Then 
\begin{equation}\label{intbound}
I = \int_{\Bbb{R}} w(t) e^{i H(t)} dt
 \ll_B (\beta - \alpha) X \left[(QR/\sqrt{Y})^{-B} + (RU)^{-B}\right]
\end{equation}
for any $B \geq 0$. 
\end{lemma}

This lemma is proved by repeated integration by parts, and we remark that in order to prove \eqref{intbound} for some fixed $B \geq 0$, \eqref{condi1} is needed only for  $j \leq j_0 = j_0(B)$. We will use this observation later in Section \ref{secprelude}. 

In the following special case we record a more precise asymptotic evaluation. 
\begin{lemma}\label {statphas} Let $x, t\in \Bbb{R}$, and let $W$ be a fixed, smooth  function with compact support on $\Bbb{R}_{>0}$. Let
$$I := \int_0^{\infty} W(y) e(xy) y^{-it} dy$$
and let $ A > 0$. 
There exists a smooth function $\tilde{W}_x(y)$ with compact support on $\Bbb{R}_{>0}$   satisfying $ \tilde{W}^{(j)}_x(y) \ll_j 1$ with the following property: if $|x| + |t| \geq 100$, then 
$$I = |x|^{-1/2} \Bigl|\frac{t}{2\pi e x}\Bigr|^{-it}\tilde{W}_x\Bigl(\frac{t}{x}\Bigr) + O_{  B}\big((|x| + |t|)^{-B}\big).$$ 
\end{lemma}

\textbf{Proof.}  
Assume that $|x| + |t| \geq 100$. If $\frac{t}{x} \not\in [c_1, c_2]$ for some suitable constants $c_2 > c_1 > 0$ (depending only on the support of $W$), we can use Lemma \ref{integrationbyparts} with $X = U = Q = 1$, $Y = |t|$, $R = |x| + |t|$ to show that $I \ll (|x| + |t|)^{-B}$. Otherwise we use \cite[Proposition 8.2]{BKY}  with $Q = V = V_1 = X = 1$, $Y = |t|$ and the unique stationary point $y_0 = t/(2 \pi x)$. With the notation of that result we have $h^{(n)}(y) = (-1)^{n-1}(n-1)! \cdot t/y^n$ for $n \geq 2$, so that the functions $p_n(y_0)$ are functions in $x$ and $t/x$. Combining them along with $|x|^{1/2} \cdot |h''(y_0)|^{-1/2}$ to $\tilde{W}_x(t/x)$ gives the result.   \hfill $\square$

\subsection{The Gamma  function} We will frequently use the functional equation and the duplication formula of the Gamma-function:
$$\Gamma(s) \Gamma(1-s) = \frac{\pi}{\sin(\pi s)}, \quad \Gamma(s)\Gamma(s+\textstyle\frac{1}{2}) = \sqrt{\pi} 2^{1-2s} \Gamma(2s).$$
For fixed $\sigma \in \Bbb{R}$, real $|t| \geq 10$ and any $M > 0$ we have Stirling's formula
\begin{equation}\label{stir}
  \Gamma(\sigma + it) = e^{-\frac{\pi}{2}|t|} |t|^{\sigma-\frac{1}{2}} \exp\left(i t \log \frac{|t|}{e}\right)g_{\sigma, M}(t) + O_{\sigma, M}(|t|^{-M}),
\end{equation}
where
$$t^j \frac{\partial^j}{\partial t^j} g_{\sigma, M}(t) \ll_{j, \sigma, M} 1.$$
for all fixed $j \in \Bbb{N}_0$. 
 
\subsection{Mellin formulae} It is useful to define the following functions. For $x > 0$, $\alpha \in \Bbb{C}$ let
\begin{equation}\label{jplus}
J_{\alpha}^+(x) := \frac{\pi}{2} \frac{J_{-\alpha}(2x) + J_{\alpha}(2x)}{\cos(\pi \alpha/2)}, \quad J_{\alpha}^-(x) := \frac{\pi}{2} \frac{J_{-\alpha}(2x) - J_{\alpha}(2x)}{\sin(\pi \alpha/2)}, \quad \tilde{K}_{\alpha}(x) = 2\cos\left(\frac{\pi}{2} \alpha\right) K_{\alpha}(2x),
\end{equation}
where $J_{\alpha}$ and $K_{\alpha}$ are the usual Bessel functions. 
We shall need the Mellin formulas
\begin{equation}\label{0}
  \left\{\begin{array}{l} \cos(x) \\ \sin(x) \end{array} \right\} = \int_{-i\infty}^{i\infty}   \left\{\begin{array}{l} \cos(\pi s/2 ) \\ \sin(\pi s/2) \end{array} \right\} \Gamma(s) x^{-s} \frac{ds}{2\pi i}, 
\end{equation}
\begin{equation}\label{2}
  J^{\pm}_\alpha(x) = \int_{-i\infty}^{i\infty} \Gamma(s+\tfrac{\alpha}{2})\Gamma(s-\tfrac{\alpha}{2})   \left\{\begin{array}{l} \sin(\pi s) \\ \cos(\pi s)\end{array} \right\}x^{-2s} \frac{ds}{2\pi i},
\end{equation}
\begin{equation}\label{3}
 \tilde{K}_\alpha(x) = \int_{-i\infty}^{i\infty} \Gamma(s+\tfrac{\alpha}{2})\Gamma(s-\tfrac{\alpha}{2}) \cos(\tfrac{\pi}{2}\alpha) x^{-2s} \frac{ds}{2\pi i},
\end{equation}
 	\begin{equation}\label{4}
\int_0^{\infty} (1+u^2)^{-s_1} (1+u^{-2})^{-s_2} u^{\alpha} \frac{du}{u} = \frac{1}{2} B\left(s_1 - \frac{\alpha}{2}, s_2 + \frac{\alpha}{2}\right), 
\end{equation}
\begin{equation}\label{5}
 \int_1^\infty (u^2-1)^{-s_1}(1-u^{-2})^{-s_2} u^{\alpha}  \frac{du}{u}  = \frac{1}{2} B\left(1-s_1-s_2,s_1-\frac{\alpha }{2}\right),
\end{equation}
\begin{equation}\label{6}
 \int_0^1 (1-u^2)^{-s_1}(u^{-2}-1)^{-s_2} u^{\alpha}  \frac{du}{u}  = \frac{1}{2} B\left(1-s_1-s_2,s_2+\frac{\alpha}{2}\right),
\end{equation}
cf.\ \cite[17.43.3/4, 17.43.16 along with functional equation of the Gamma function, 17.43.18, 17.43.7, 3.191.2 with $u=1$, 3.191.3]{GR}. In \eqref{2} (and in \eqref{intrep1} below), the upper sign ($+$) belongs to the upper trigonometric function ($\sin$) and the lower sign ($-$) to the lower trigonometric function ($\cos$). Moreover,  $B$ is the Euler Beta function, and  we recall the Barnes convention from Definition \ref{MB}. In \eqref{0} -- \eqref{3} it is understood that $x > 0$, and in \eqref{4} -- \eqref{6} the two arguments of the Beta function must have positive real part to make the integrals absolutely convergent.

\subsection{Bessel functions} We start with the integral representations \cite[8.432.4]{GR}   and \cite[8.421.1/2, 8.405]{GR}
\begin{equation}\label{intrep1}
\begin{split} 
 \tilde{K}_{it}(x) =   \int_{-\infty}^{\infty} \cos(2x \sinh v)\exp( it v) dv, \quad J^{\pm}_{it}(x) =   \int_{-\infty}^{\infty} \left\{\begin{array}{l} \sin\\ \cos
\end{array}\right\}(2x \cosh v)\exp( it v) dv
\end{split}
 \end{equation}
for $t \in \Bbb{R}$, $x > 0$.  The integrals are not absolutely convergent, but integration by parts shows that the tail is  very small, so that the conditional convergence causes no extra difficulty. We can use these representations to obtain the uniform bounds
\begin{equation}\label{besselbound}
  \frac{\partial^j}{\partial x^j} \tilde{K}_{it}(x), \quad \frac{\partial^j}{\partial x^j} J^{\pm}_{ it}(x) \ll_j \left( 1 + \frac{|t|}{x} \right)^{j}
\end{equation}
for $|t|, x \geq 1$ and $j \in \Bbb{N}_0$.  Indeed, if $|t|/x \geq 100$ we cut the integrals \eqref{intrep1} smoothly (using a smooth partition of unity)  into the region 
\begin{equation}\label{first}
  |v| \leq (\log |t|/x) - 10,
\end{equation}
the region 
\begin{equation}\label{second}
(\log |t|/x) -10 \leq |v| \leq (\log |t|/x) +10
\end{equation}
 and 
\begin{equation}\label{third} 
 (\log |t|/x) +10\kappa \leq |v| \leq (\log |t|/x) +10(\kappa+1), \quad \kappa = 1, 2, \ldots.
 \end{equation}
 In each region we differentiate $j$ times with respect to $x$. In \eqref{first} we integrate by parts using Lemma \ref{integrationbyparts}  with 
 $$H(v) = tv \pm 2x \left\{\begin{array}{l} \sinh v\\ \cosh v
\end{array}\right\},$$
$X = (|t|/x)^j$, $\beta - \alpha \ll \log |t|/x$, $U = \log |t|/x$, $R =Y = |t|$, $Q = 1$ to see that this portion is negligible. The region \eqref{second} contains a possible  stationary point, and here we estimate trivially.  In \eqref{third} we integrate by parts using Lemma \ref{integrationbyparts}  with $X = (e^{10\kappa} |t|/x)^j$, $\beta - \alpha = 10$, $U = Q = 1$, $R =Y = x e^{10\kappa}$ which is again a negligible contribution. 

If $|t|/x \leq 100$, we   estimate trivially the range $|v| \leq 100$ and show as above that the contribution of each interval $100 + 10\kappa \leq |v| \leq 100 + 10(\kappa+1)$ is negligible. This proves \eqref{besselbound}. \\
 
We proceed  with the following  uniform asymptotic formulae (which can in principle be obtained from \eqref{intrep1} by a careful stationary phase argument). We have 
\begin{equation}\label{unifasymp}
\tilde{K}_{it}(x/2) =  \Re \left( e^{i   \omega(x, t)}  f_M(x, t)\right)  
+ O(|t|^{-M}), \quad \omega(x, t) =  |t| \cdot \text{arccosh} \frac{|t|}{x} - \sqrt{t^2 - x^2},
\end{equation}
for $t \in \Bbb{R}$, $|t| > 1$, $\frac{1}{10} |t|\geq  x > 0$ and fixed $M > 0$ with
$$t^j \frac{\partial^j}{\partial t^j} f_M(t, x) \ll_{j, M}  |t|^{-1/2}$$ 
for any $j \in \Bbb{N}_0$, see \cite[7.13.2(19)]{EMOT}.  The error term there is only $O(x^{-M})$, but for  $x \leq |t|^{1/10}$, say, the formula \eqref{unifasymp} follows from the power series expansion for 
\begin{equation}\label{I}
  \tilde{K}_{\alpha}(x) = \frac{\pi}{2} \frac{1}{\sin(\pi \alpha/2)} \left(I_{-\alpha}(2x) - I_{\alpha}(2x)\right),
\end{equation}
\begin{equation}\label{I1}
I_{\alpha}(2x) = \sum_{k=0}^{\infty} \frac{1}{k! \Gamma(k+1+\alpha)}x^{\alpha+2k}. 
\end{equation}

Analogously, we have
\begin{equation}\label{unifasymp1}
 J^{\pm}_{it}(x/2) = 
 \Re \left( e^{i  \tilde{\omega}(x, t)} \tilde{f}_M^{\pm}(x, t)\right) 
 + O(|t|^{-M}),  
\quad \tilde{\omega}(x, t) =  |t| \cdot \text{arcsinh} \frac{|t|}{x} - \sqrt{t^2 + x^2},
\end{equation}
for $t \in \Bbb{R}$,  $|t| > 1$, $x > 0$ and fixed $M > 0$ with
$$t^j \frac{\partial^j}{\partial t^j} \tilde{f}^{\pm}_M(t, x) \ll_{j, M}  \frac{1}{x^{1/2} + |t|^{1/2}}$$ 
for any $j \in \Bbb{N}_0$, see \cite[7.13.2(17)]{EMOT}. Notice that \eqref{unifasymp1}  holds without the restriction $x \leq \frac{1}{10} |t|$ (there is no ``transitional range'').   Again the error term in \cite{EMOT} is $O(x^{-M})$, but for small $x$ the error term $O(|t|^{-M})$ follows from the power series expansion
\begin{equation}\label{powerJ}
J_{\alpha}(2x) = \sum_{k=0}^{\infty} \frac{(-1)^k}{k! \Gamma(k+1+\alpha)}x^{\alpha+2k}. 
\end{equation}


\section{Integral representations} In this section we establish alternative expressions for the kernel functions given  in Definition \ref{MB} in terms of the Bessel functions $J^{\pm}$ and $\tilde{K}$ defined in \eqref{jplus}. These representations will play an important role later,  but are also of independent interest. 

 
\begin{lemma}\label{stade0} For $y\in \Bbb{R} \setminus \{0\}$ and $\mu \in \Lambda_0$ we have
\begin{equation*}
\begin{split}
K_{w_4}(y; \mu) = &\frac{1}{3072 \pi^5} \int_0^{\infty}  J^{-}_{\mu_1 - \mu_2}(2\sqrt{u})  \left(\frac{\pi^3 |y|}{u^{3/2}}\right)^{-\mu_3 } \exp\left(-\frac{2i \pi^3 y}{u} \right) \frac{du}{u}\\
& +  \frac{1}{3072\pi^5}  \int_0^{\infty}     \tilde{K}_{\mu_1 - \mu_2}(2\sqrt{u}) \left(\frac{\pi^3 |y|}{u^{3/2}}\right)^{-\mu_3} \exp\left(\frac{2i \pi^3 y}{u}\right) \frac{du}{u}. 
\end{split}
\end{equation*}
\end{lemma}

 \textbf{Remarks.} The integrals     just fail to  be absolutely convergent at $0$, but since $\exp(\pm 2i\pi^3y/u)$ is highly oscillating in a neighbourhood of  $u=0$, the integrals exist in a Riemann sense, and the portion $0 < u < 1$ can be made absolutely convergent after partial integration.  It follows from the definition that $K_{w_4}(y; \mu)$ is Weyl-group invariant. This is not easily visible from the above formula.  \\

 \textbf{Proof.} By \eqref{2} and \eqref{3} along with the duplication formula and the functional equation of the Gamma function and the fact that $\mu_1 + \mu_2 + \mu_3 = 0$, we have 
$$\frac{2}{\pi} \int_0^{\infty} u^{\frac{1}{2}\mu_3 } \left(J_{\mu_1 - \mu_2}^{-}(2 \sqrt{u})  + \tilde{K}_{\mu_1 - \mu_2}(2 \sqrt{u}) \right) u^{s-1} du =    \frac{\Gamma(\frac{1}{2}(s - \mu_1))\Gamma(\frac{1}{2}(s - \mu_2))}{\Gamma(\frac{1}{2}(1-s + \mu_1))\Gamma(\frac{1}{2}(1-s + \mu_2))}$$
(the left hand side is absolutely convergent in $0 < \Re s < 1/4$), and by \eqref{0} we have
\begin{equation}\label{notL1}
 \frac{2}{\sqrt{\pi}} \int_0^{\infty} u^{-\mu_3}\cos(2u)  u^{s-1} du =  \frac{\Gamma(\frac{1}{2}(s - \mu_3)) }{\Gamma(\frac{1}{2}(1-s + \mu_3)) }
 \end{equation}
(the left hand side is conditionally convergent in $0 < \Re s < 1$). By a formal application of Parseval's identity  we conclude 
\begin{equation}\label{parseval}
\begin{split}
&\int_{-i\infty}^{i\infty}  \frac{\Gamma(\frac{1}{2}(s - \mu_1))\Gamma(\frac{1}{2}(s - \mu_2))}{\Gamma(\frac{1}{2}(1-s + \mu_1))\Gamma(\frac{1}{2}(1-s + \mu_2))} \frac{\Gamma(\frac{1}{2}(s - \mu_3)) }{\Gamma(\frac{1}{2}(1-s + \mu_3)) } y^{-s}  \frac{ds}{2\pi i} \\
&=  \frac{4}{\pi^{3/2}} \int_0^{\infty} u^{\frac{1}{2}\mu_3 } \left(J_{\mu_1 - \mu_2}^{-}(2 \sqrt{u})  + \tilde{K}_{\mu_1- \mu_2}(2 \sqrt{u}) \right)  \left(\frac{y}{u}\right)^{-\mu_3}\cos\left(\frac{2y}{u}\right) \frac{du}{u}
\end{split}
\end{equation}
for $y > 0$. Since the integrand in  \eqref{notL1} is not in $L^1$, this formal argument needs some justification. One way is to work instead with the Mellin pair
$$ 
 \frac{2}{\sqrt{\pi}} \int_0^{\infty} e^{-\varepsilon u} u^{-\mu_3}\cos(2u)  u^{s-1} du =  \frac{\Gamma(\frac{1}{2}(s - \mu_3)) }{\Gamma(\frac{1}{2}(1-s + \mu_3)) }  \frac{2^{s-1}((2i +\varepsilon)^s + (\varepsilon - 2i)^s)}{(4+\varepsilon^2)^s \cos(\pi s/2)}  $$
 for small $\varepsilon > 0$, which can also be derived from \eqref{0}. Then on the left hand side of a correspondingly modified version of \eqref{parseval} we can let $\varepsilon \rightarrow 0$ inside the integral, e.g.\ by dominated convergence (recall the Barnes integral convention). On the right hand side we cannot use $L^1$-theory directly since the pointwise limit is not in $L^1$. However, we can split the integral into two pieces $u< 1$ and $u \geq 1$. The latter is in $L^1$, and in the former we can first integrate by parts (using that $\cos(2y/u)$ is highly oscillating) to obtain an $L^1$-integrand, then interchange limit and integration, and finally integrate by parts backwards in a Riemann sense. This proves \eqref{parseval}.

Similarly,  
$$\frac{2}{\pi} \int_0^{\infty} u^{\frac{1}{2}\mu_3  } \left(-J_{\mu_1 - \mu_2}^{-}(2 \sqrt{u})  + \tilde{K}_{\mu_1- \mu_2}(2 \sqrt{u}) \right) u^{s-1} du =   \frac{\Gamma(\frac{1}{2}(1+s - \mu_1))\Gamma(\frac{1}{2}(1+s - \mu_2))}{\Gamma(\frac{1}{2}(2-s + \mu_1))\Gamma(\frac{1}{2}(2-s + \mu_2))}.$$
and 
$$
 \frac{2}{\sqrt{\pi}} \int_0^{\infty} u^{-\mu_3}\sin(2u)  u^{s-1} du =  \frac{\Gamma(\frac{1}{2}(1+s - \mu_3)) }{\Gamma(\frac{1}{2}(2-s + \mu_3)) },$$
 so that by the same argument 
 \begin{displaymath}
\begin{split}
&\int_{-i\infty}^{i\infty}  \frac{\Gamma(\frac{1}{2}(1+s - \mu_1))\Gamma(\frac{1}{2}(1+s - \mu_2))}{\Gamma(\frac{1}{2}(2-s + \mu_1))\Gamma(\frac{1}{2}(2-s + \mu_2))} \frac{\Gamma(\frac{1}{2}(1+s - \mu_3)) }{\Gamma(\frac{1}{2}(2-s + \mu_3) } y^{-s}  \frac{ds}{2\pi i} \\
&=  \frac{4}{\pi^{3/2}} \int_0^{\infty} u^{\frac{1}{2}\mu_3 } \left(-J_{\mu_1- \mu_2}^{-}(2 \sqrt{u})  + \tilde{K}_{\mu_1 - \mu_2}(2 \sqrt{u}) \right)  \left(\frac{y}{u}\right)^{-\mu_3}\sin(2y/u) \frac{du}{u}.
\end{split}
\end{displaymath}
The integral formula for $K_{w_4}$ from Definition \ref{MB} follows now easily.  \hfill $\square$ \\

The function $K_{w_6}^{++}(y; \mu)$ is essentially a ${\rm GL}(3)$-Whittaker function, and the integral kernels with different signs are close relatives. For such functions one can derive a nice integral representation in terms of standard Bessel functions in the spirit of  \cite[(6.1.3)]{Go}. 
For $y_1, y_2 \in \Bbb{R} \setminus \{0\}$ and $\mu \in \Lambda_0$ we define the absolutely convergent expressions 
\begin{equation}\label{BesselJ1}
 \mathcal{J}^{\pm}_{1}(y; \mu) = 
  \Bigl| \frac{y_1}{y_2}\Bigr|^{\frac{1}{2}\mu_2}      \int_0^\infty J^{\pm}_{3\nu_3 }\left(2\pi |y_1|^{1/2}\sqrt{1+u^2}\right) J^{\pm}_{3\nu_3}\left(2\pi |y_2|^{1/2}\sqrt{1+u^{-2}}\right) u^{3\mu_2} \frac{du}{u},
 \end{equation}
\begin{equation}\label{BesselJ2}
 \mathcal{J}_{2}(y; \mu) = 
 \Bigl| \frac{y_1}{y_2}\Bigr|^{\frac{1}{2}\mu_2}      \int_1^\infty J^{-}_{3\nu_3 }\left(2\pi |y_1|^{1/2}\sqrt{u^2-1}\right) J^{-}_{3\nu_3}\left(2\pi |y_2|^{1/2}\sqrt{1-u^{-2}}\right) u^{3\mu_2} \frac{du}{u},
 \end{equation}
\begin{equation}\label{BesselJ3}
 \mathcal{J}_{3}(y; \mu) = 
 \Bigl| \frac{y_1}{y_2}\Bigr|^{\frac{1}{2}\mu_2}      \int_0^\infty \tilde{K}_{3\nu_3 }\left(2\pi |y_1|^{1/2}\sqrt{1+u^2}\right) J^{-}_{3\nu_3}\left(2\pi |y_2|^{1/2}\sqrt{1+u^{-2}}\right) u^{3\mu_2} \frac{du}{u},
 \end{equation}
\begin{equation}\label{BesselJ4}
 \mathcal{J}_{4}(y; \mu) = 
  \Bigl| \frac{y_1}{y_2}\Bigr|^{\frac{1}{2}\mu_2}      \int_0^1 \tilde{K}_{3\nu_3 }\left(2\pi |y_1|^{1/2}\sqrt{1-u^2}\right) \tilde{K}_{3\nu_3}\left(2\pi |y_2|^{1/2}\sqrt{u^{-2}-1}\right) u^{3\mu_2} \frac{du}{u},
 \end{equation}
   \begin{equation}\label{BesselK}
 \mathcal{J}_{5}(y;\mu) =  \Bigl|\frac{y_1 }{y_2 }\Bigr|^{\frac{1}{2}\mu_2}   \int_0^\infty \tilde{K}_{3\nu_3 }\left(2\pi |y_1|^{1/2}\sqrt{1+u^2}\right) \tilde{K}_{3\nu_3}\left(2\pi |y_2|^{1/2}\sqrt{1+u^{-2}}\right) u^{3\mu_2} \frac{du}{u}.
\end{equation}

\begin{lemma}\label{stade}  For $y_1, y_2 > 0$ we have
\begin{equation}\label{lem++}
   K_{w_6}^{++}(y; \mu)    = \frac{1}{12\pi^2}\frac{ \cos\left(\frac{3}{2}\pi \nu_1\right) \cos\left(\frac{3}{2}\pi \nu_2\right)}{ \cos\left(\frac{3}{2}\pi \nu_3\right)} \mathcal{J}_{5}(y;\mu);
   \end{equation}
for $y_1 > 0 > y_2$ we have
   \begin{equation}\label{lem+-}
   \sum_{w \in \{I, w_4, w_5\}}  K_{w_6}^{+-}(y; w(\mu))    = \frac{1}{24\pi^2} \sum_{w \in \{I, w_4, w_5\}} \Bigl( \mathcal{J}_2(y; w(\mu)) +  \mathcal{J}_3(y; w(\mu)) +  \mathcal{J}_4(y; w(\mu))  \Bigr) ; 
      \end{equation}
 for $y_2 > 0 > y_1$ we have
   \begin{equation}\label{lem-+}
   \changed{K_{w_6}^{-+}((y_1, y_2); \mu)    =    K_{w_6}^{+-}((y_2, y_1); w_4(-\mu))};
         \end{equation}
 and for $y_1, y_2 < 0$ we have
   \begin{equation}\label{lem--}
   \sum_{w \in \{I, w_4, w_5\}} K_{w_6}^{--}(y; w(\mu))  =   \frac{1}{48\pi^2} \sum_{w \in \{I, w_4, w_5\}} \Bigl(4 \mathcal{J}_1^{-}(y; w(\mu)) +2  \mathcal{J}_1^+(y; w(\mu))  \Bigr).
\end{equation}
\end{lemma}

 \textbf{Remark.}  It is a very challenging exercise find these identities, but once they are given, it is  a straightforward (but  tedious) exercise in trigonometry and Mellin inversion to prove them.\\ 



\textbf{Proof.} To prove \eqref{lem++}, we insert the Mellin formula  \eqref{3} for both $\tilde{K}$-factors in \eqref{BesselK}
 and compute the $u$-integral using \eqref{4}. After changing variables $s_1 - \frac{1}{2} \mu_2 \mapsto s_1$, $s_2 + \frac{1}{2} \mu_2 \mapsto s_2$ and using $\mu_1 + \mu_2 + \mu_3= 0$, this produces the Mellin-Barnes integral for $K_{w_6}^{++}$   from Definition \ref{MB}. \\
 
The proof of \eqref{lem--} uses the same argument, followed by the 
trigonometric identity
\begin{displaymath}
\begin{split}
32\pi^2  \sum_{w \in \{I, w_4, w_5\}}    S^{--}(s,  w(\mu)) 
= & \frac{1}{3}  \sum_{w \in \{I, w_4, w_5\}}  \Bigl[4 \cos\left(\pi(s_1 + \textstyle\frac{1}{2}  w(\mu)_2)\right)  \cos\left(\pi(s_2 - \textstyle\frac{1}{2}  w(\mu)_2)\right)   \\
& \quad\quad\quad \quad\quad  
+ 2 \sin\left(\pi(s_1 + \textstyle\frac{1}{2}  w(\mu)_2)\right)  \sin\left(\pi(s_2 - \textstyle\frac{1}{2} w(\mu)_2)\right)\Bigr]
\end{split}
\end{displaymath}
for $\mu_1 + \mu_2 + \mu_3 = 0$, which can be verified by brute force. \\

Equation \eqref{lem-+} follows directly from Definition \ref{MB}, see also \eqref{3cycle}. \\

Finally, for   the   proof of \eqref{lem+-} we  insert \eqref{2} and \eqref{3}  in \eqref{BesselJ2} -- \eqref{BesselJ4},   compute the $u$-integral with \eqref{4} -- \eqref{6}, and apply the functional equation of the Gamma function to obtain
 \begin{displaymath}
 \begin{split}
&    \mathcal{J}_j(y; \mu ) =  \frac{1}{2}   \int_{-i\infty}^{i\infty} \int_{-i\infty}^{i\infty} \frac{G(s, \mu)}{\sin(\pi(s_1+s_2))} R_j(s, \mu)
 |4\pi y_1|^{-s_1}|4\pi y_2|^{-s_2}  \frac{ds_1 \, ds_2}{(2\pi i)^2} 
  \end{split}
 \end{displaymath}
 for $j = 2, 3, 4$, where
 \begin{displaymath}
 \begin{split}
& R_2(s,  \mu)  = \cos\left(\pi (s_1 +  \tfrac{1}{2} \mu_2  )\right)\cos\left(\pi (s_2 - \tfrac{1}{2} \mu_2)\right) \sin\left(\pi(s_2 +  \mu_2)\right),\\
& R_3(s,  \mu) = \sin(\pi(s_1+s_2)) \cos\left(\tfrac{\pi}{2}( \mu_3 -  \mu_1)\right) \cos \left(\pi(s_2 - \tfrac{ \mu_2}{2})\right),\\
& R_4(s,  \mu) = \cos\left(\tfrac{\pi}{2}( \mu_3 -  \mu_1)\right)^2 \sin\left(\pi(s_1 -  \mu_2)\right).
 \end{split}
 \end{displaymath}
 Using the trigonometric identity
 \begin{displaymath}
 \begin{split}
  -32\pi^2& \sum_{w \in \{I, w_4, w_5\}}   S^{+-}(s,  w(\mu)) \sin(\pi(s_1+s_2))\\
  & =-\frac{2}{3} \sum_{w \in \{I, w_4, w_5\}}\big(R_2(s,  w(\mu)) + R_3(s,  w(\mu)) + R_4(s,  w(\mu))\big),
  \end{split}
  \end{displaymath}
 which again can be verified by brute force (but is very challenging to find), we see that   \eqref{lem+-} follows  from Definition \ref{MB}. \hfill $\square$

 \section{Some finite Fourier transforms}

This section contains  bounds for multiple Fourier transforms of the Kloosterman sums in the Kuznetsov formula. 

\begin{lemma}\label{lem1} Let $s_1, s_2, r_1, n_2 \in \Bbb{Z} \setminus \{0\}$. Let $D, \delta \in \Bbb{N}$ and $x, y \in \Bbb{Z}$. Then
\begin{equation}\label{Klo5}
\begin{split}
\Bigl|\frac{1}{D\delta}&\sum_{ n_1 \, (\text{{\rm mod }}D)}\sum_{ m_1 \, (\text{{\rm mod }} \delta)} \tilde{S}(n_1r_1, n_2s_2, m_1s_1; D, D\delta ) e\left( - \frac{xn_1}{D} - \frac{ym_1}{\delta}\right)\Bigr| \leq D (r_1,  D)(s_1,  \delta),  
\end{split}
\end{equation}
and the left hand side vanishes unless $(D, x) = (r_1, x)$,  $(\delta, y) = (s_1, y)$ and $D \mid n_2s_2$. 
\end{lemma}

\textbf{Proof.} 
The left hand side of \eqref{Klo5} equals
\begin{displaymath}
\begin{split}
& \frac{1}{D\delta} \sum_{\substack{C_1\, (\text{{\rm mod }} D)\\ (C_1, D) = 1}} \sum_{\substack{C_2 \, (\text{{\rm mod }}D\delta)\\ (C_2, \delta) = 1}}\sum_{ \substack{n_1 \, (\text{{\rm mod }}D)\\ m_1 \, (\text{{\rm mod }}\delta)} }  e\left(s_2n_2\frac{\bar{C}_1 C_2}{D} + s_1m_1 \frac{\bar{C}_2}{\delta} + r_1n_1 \frac{C_1}{D}\right)  e\left( - \frac{xn_1}{D} - \frac{ym_1}{\delta}\right)\\
& = \sum_{\substack{C_1\, (\text{{\rm mod }}D)\\ (C_1, D) = 1\\ r_1C_1 \equiv x \, (\text{{\rm mod }}D)}} \sum_{\substack{C_2 \, (\text{{\rm mod }}D\delta)\\ (C_2, \delta) = 1\\ s_1 \bar{C}_2 \equiv y \, (\text{{\rm mod }}\delta)}}e\left(s_2n_2\frac{\bar{C}_1 C_2}{D} \right),   
\end{split}
\end{displaymath}
and the result follows. \hfill $\square$\\

For $r_1, r_2, s_1, s_2 \in \Bbb{Z} \setminus \{0\}$, $x_1, x_2, y_1, y_2 \in \Bbb{Z}$ and $D_1, D_2 \in \Bbb{N}$ we define
\begin{displaymath}
\begin{split}
&\widehat{S}_{r_1, s_1, r_2, s_2}(x_1, y_1, x_2, y_2; D_1, D_2) := \\
&\frac{1}{D_1^2D_2^2} \sum_{\substack{n_1, m_1\, (\text{{\rm mod }} D_1)\\ n_2, m_2 \, (\text{\rm{mod  }} D_2)}} S(n_1r_1, m_2r_2, m_1s_1, n_2s_2; D_1, D_2)e\left( - \frac{x_1m_1 + y_1 n_1}{D_1} - \frac{x_2m_2 + y_2 n_2}{D_2}\right).
\end{split}
\end{displaymath}

For the following lemma we introduce some notation. As usual we denote Euler's function by $\phi$.  For a prime $\ell$ we write $r \mid \ell^{\infty}$ if $r$ is a power of $\ell$, and we denote by $(\ell^{\infty}, r)$ the highest power of $\ell$ dividing $r$.

\begin{lemma}\label{lem2} {\rm (a)} We have the general bound
 \begin{equation*}\label{Klo6}
 \begin{split}
|\widehat{S}_{r_1, s_1, r_2, s_2}(x_1, y_1, x_2, y_2; D_1, D_2)| \leq  (r_1, D_1)  (r_2, D_2)  (D_1, D_2),  
  \end{split}
  \end{equation*}
  and the left hand side vanishes unless $x_1y_1 \equiv r_1s_1 D_2 \, (\text{{\rm mod }}D_1)$ and $x_2y_2 \equiv r_2s_2 D_1 \, (\text{{\rm mod  }}D_2)$. \\
  {\rm (b)}  We have $\widehat{S}_{r_1, s_1, r_2, s_2}(0, 0, x_2, y_2; D_1, D_2)  = 0$  unless
  $$\frac{D_1}{(D_1, r_1s_1)} \mid (x_2, y_2).$$
     Similarly,  
  $\widehat{S}_{r_1, s_1, r_2, s_2}(x_1, y_1, 0, 0; D_1, D_2)  = 0$ unless $\frac{D_2}{(D_2, r_2s_2)} \mid (x_1, y_1).$\\
  {\rm (c)} If $(r_1r_2, s_1s_2) = 1$, then $\widehat{S}_{r_1, s_1, r_2, s_2}(0, 0, 0, 0; D_1, D_2)  = 0$ unless $D_1 = D_2$, in which case it equals $\phi(D)$. \\
  {\rm (d)} Let $\ell$ be a prime and assume that $r_1r_2s_1s_2 \mid \ell^{\infty}$. Then
 \begin{equation}\label{partd}
 |\widehat{S}_{r_1, s_1, r_2, s_2}(0, 0, 0, 0; D_1, D_2)| \leq (D_1, D_2) (\ell^{\infty}, [D_1, D_2])r_2.
 \end{equation}
   \end{lemma}
 
\textbf{Remark.} Parts (c) and (d) will be used for the treatment of the central Poisson term. The key point is that there is no cancellation in the sum
$$\frac{1}{D^4} \sum_{\substack{n_1, m_1\, (\text{{\rm mod }} D_1)\\ n_2, m_2 \, (\text{\rm{mod  }} D_2)}} S(n_1r_1, m_2r_2, m_1s_1, n_2s_2; D, D),$$
see also \cite[Property 4.10]{BFG}. This  is very different from the ${\rm GL}(2)$ case, where always
$$\frac{1}{c^2} \sum_{n, m \, (\text{mod } c)} S(m, n, c) = 0.$$
The  right hand side of \eqref{partd} could be made more symmetric (and slightly sharper), but the present form suffices for our needs. \\
 
\textbf{Proof.} By twisted multiplicativity of Kloosterman  sums \cite[Property 4.7]{BFG} the quantities  $\widehat{S}_{r_1, s_1, r_2, s_2}(x_1, y_1, x_2, y_2; D_1, D_2)$ also enjoy  twisted multiplicativity, and we have
\begin{displaymath}
\begin{split}
 & \widehat{S}_{r_1, s_1, r_2, s_2}(x_1, y_1, x_2, y_2; D_1D_1', D_2D_2')  \\
  &= \widehat{S}_{r_1, s_1, r_2, s_2}(x_1\bar{D}'_1, y_1D'_1\bar{D}'_2, x_2D_2'\bar{D}_1', y_2\bar{D}_2'; D_1 , D_2 ) \widehat{S}_{r_1, s_1, r_2, s_2}(x_1\bar{D}_1, y_1D_1\bar{D}_2, x_2D_2\bar{D}_1, y_2\bar{D}_2;  D_1',  D_2')
\end{split}
\end{displaymath}
for $(D_1D_2, D_1'D_2') = 1$. Hence  it suffices to assume   that $D_j = q^{\alpha_j}$ are  powers of a prime $q$. By orthogonality of additive characters, $\widehat{S}_{r_1, s_1, r_2, s_2}(x_1, y_1, x_2, y_2; q^{\alpha_1}, q^{\alpha_2})$    equals the number of solutions $B_j, C_j \, (\text{mod } q^{\alpha_j})$ with $(B_j, C_j, q^{\alpha_j}) = 1$ satisfying 
\begin{equation}\label{count1}
 r_1B_1 \equiv  y_1 \, ({\rm mod } \, q^{\alpha_1}),
 \end{equation}
 \begin{equation}\label{count2}
  r_2B_2 \equiv  x_2 \, ({\rm mod } \, q^{\alpha_2}),
 \end{equation}
 \begin{equation}\label{count3}
 s_1(Y_1q^{\alpha_2} - Z_1B_2) \equiv x_1 \, ({\rm mod } \, q^{\alpha_1}),
 \end{equation}
\begin{equation}\label{count4}
s_2(Y_2q^{\alpha_1} - Z_2B_1) \equiv  y_2 \, ({\rm mod } \, q^{\alpha_2}),
\end{equation}
\begin{equation}\label{equa}
B_1B_2 + C_1 q^{\alpha_2} + C_2 q^{\alpha_1} \equiv 0 \, ({\rm mod } \, q^{\alpha_1+\alpha_2}),
  \end{equation}
where  $Y_j B_j+ Z_j C_j \equiv 1 \, ({\rm mod } \, q^{\alpha_j})$.\\

 \noindent   (a) If $q \nmid B_1$, then  we can choose $Y_1 = \bar{B}_1$, $Z_1 = 0$ obtaining $ s_1r_1 q^{\alpha_2} \equiv x_1y_1 \, (\text{mod }q^{\alpha_1})$ from  \eqref{count1} and \eqref{count3}.  
 On the other hand, if $q \nmid C_1$, we can choose $Z_1 = \bar{C}_1$, $Y_1 = 0$ obtaining  $-s_1 \bar{C}_1B_2   \equiv x_1 \, ({\rm mod } \, q^{\alpha_1})$. Multiplying \eqref{equa} with $\bar{C}_1$, this implies, in connection with   \eqref{count1}, also $x_1y_1 \equiv  r_1s_1 q^{\alpha_2} \, ({\rm mod } \, q^{\alpha_1})$. The same argument works with exchanged indices, and we conclude that the sum in question vanishes unless
 $$x_1y_1 \equiv r_1s_1 q^{\alpha_2} \, (\text{mod }q^{\alpha_1}), \quad x_2y_2 \equiv r_2s_2 q^{\alpha_1} \, (\text{mod }q^{\alpha_2}).$$
From \eqref{count1} and \eqref{count2},  the number of choices for $B_1, B_2 $ is at most $(r_1,  q^{\alpha_1})(r_2,   q^{\alpha_2}) $. Having these two fixed, we can now choose $C_1$ or $C_2$ freely, then the other variable is determined by \eqref{equa}. Hence we conclude that 
the sum in question is
$$\leq (r_1, q^{\alpha_1})  (r_2, q^{\alpha_1})  q^{\min(\alpha_1, \alpha_2)} $$ 
as desired. \\
(b) Let $x_1 = y_1 = 0$. We denote by $v_q$ the usual $q$-adic valuation. There is nothing to prove unless  
\begin{equation}\label{max}
\max(v_q(r_1), v_q(s_1)) < \alpha_1
\end{equation}
which we assume from now on. Then \eqref{count1} reads  
\begin{equation}\label{previous}
q^{\alpha_1} \mid r_1B_1,
\end{equation} which implies $q \mid B_1$, so $q \nmid C_1$. With $Y_1 = 0$, $Z_1 = \bar{C}_1$ we obtain $q^{\alpha_1} \mid s_1B_2$ from  \eqref{count3}.  From this and \eqref{count2} we therefore conclude $v_q(x_2) \geq \min(\alpha_2, \alpha_1 + v_q(r_2) - v_q(s_1))$. From the second conclusion of part (a)  with $x_1 = y_1 = 0$ we know that $\alpha_2 \geq \alpha_1 - v_q(r_1) - v_q(s_1)$, and the desired divisibility condition for $x_2$ follows.  We have already   proved $q^{\alpha_1} \mid s_1B_2$,  which together with  \eqref{max}  implies 
 $q\mid B_2$, $q \nmid C_2$, so that with $Z_2 = \bar{C}_2$, $Y_2 = 0$ we obtain $-s_2 \bar{C}_2 B_1 \equiv y_2 \, (\text{mod } q^{\alpha_2})$ from \eqref{count4}. This together with \eqref{previous} proves the divisibility condition for $y_2$ by the same argument. \\
(c) The $w_6$ Kloosterman sum satisfies the symmetry property \cite[Property 4.4]{BFG}
\begin{equation*}
S(n_1, m_2, m_1, n_2; D_1, D_2) = S(m_1, n_2, n_1, m_2; D_1, D_2),
\end{equation*}
so that without loss of generality we can assume that $q \nmid r_1r_2$ (otherwise we interchange $s_1, s_2$ with $r_1, r_2$). By \eqref{count1} and \eqref{count2}, we have $B_1 = B_2 = 0$ and so $q \nmid C_1C_2$, so that the condition \eqref{equa}  implies $\alpha_1 = \alpha_2$, and $C_2 = -C_1$ can be chosen freely, but coprime to $q$. \\
(d) For $D_1, D_2$ coprime to $\ell$, the desired bound follows from  the conclusion of part (a), so by twisted multiplicativity it suffices to consider the case $D_1 = \ell^{\alpha_1}$, $D_2 = \ell^{\alpha_2}$. Let $r_j = \ell^{\rho_j}$, $s_j = \ell^{\sigma_j}$ ($j = 1, 2$). We need to show
\begin{equation*}
|\widehat{S}_{r_1, s_1, r_2, s_2}(0, 0, 0, 0; \ell^{\alpha_1}, \ell^{\alpha_2})| \leq \ell^{\alpha_1+\alpha_2 + \rho_2}. 
\end{equation*}
We have at most $\ell^{\rho_2}$ choices for $B_2$ from \eqref{count2} and trivially $\ell^{\alpha_1}$ choices for $B_1$.  Once they are determined, we conclude as above that there are at most $\ell^{\min(\alpha_1, \alpha_2)}$ choices for the pair $(C_1, C_2)$. This proves the claim and completes the proof of the lemma. 
  \hfill $\square$

\section{Key lemmas on the weight functions}\label{key}

In this section we summarize   important properties of the functions $\Phi_w$ defined in \eqref{defPhi}. The proofs are long and difficult, and we postpone them to the end of the paper.  We have not aimed for the greatest possible generality, but rather for a compact presentation of the necessary bounds needed in our particular application. 
We recall that the test function $h$ (depending on $\varepsilon$,  $A$ and $\mu_0$) was specified in Section \ref{hh}, 
and $\mu_0$ satisfies \eqref{munu}. 

The following two lemmas will be used to truncate various sums. They show on the one hand that $\Phi_w$ is negligibly small for small arguments (thereby quantifying and improving Lemma \ref{lemma1}), and bound on the other hand the derivatives of $\Phi_w$ which in connection with Lemma \ref{integrationbyparts} can be used to truncate sums with Fourier integrals containing $\Phi_w$. Statements and proofs of Lemma \ref{lem3} and Lemma \ref{lem4} do not use $\varepsilon$-convention. 

\begin{lemma}\label{lem3} 
Let $0 < |y| \leq T^{3 - \varepsilon}$. Then for any constant $B \geq 0$ one has
\begin{equation}\label{w41}
\Phi_{w_4} (y) \ll_{\varepsilon, B} T^{-B}.
\end{equation}
If $T^{3-\varepsilon} < |y|$, then  
\begin{equation}\label{w42}
 |y|^j   \Phi^{(j)}_{w_4}(y)   \ll_{j, \varepsilon}   T^{3+2\varepsilon} (T + |y|^{1/3})^j
 \end{equation}
for any $j \in \Bbb{N}_0$. 
\end{lemma}

The true order of magnitude of $\Phi_{w_4}$ is essentially $T^{2}$, but this is of little relevance here.  The focus of this lemma is on the cut-off point $y \gg T^{3+o(1)}$ and the size of the oscillation, and these bounds are sharp (and need to be sharp for our purposes). 

 
\begin{lemma}\label{lem4} 
Let $\Upsilon := \min(|y_1|^{1/3} |y_2|^{1/6}, |y_1|^{1/6} |y_2|^{1/3})$. If $\Upsilon \leq T^{1-\varepsilon}$, then
\begin{equation}\label{J1}
\Phi_{w_6}(y_1, y_2)   \ll_{B, \varepsilon}   T^{-B}
 \end{equation}
for any fixed constant $B \geq 0$. If $\Upsilon \geq T^{1-\varepsilon}$, then 
 \begin{equation}\label{J2}
\begin{split}
|y_1|^{j_1} |y_2|^{j_2} &\frac{\partial^{j_1}}{\partial y_1^{j_1} }\frac{\partial^{j_2}}{\partial y_2^{j_2} } \Phi_{w_6}(y_1, y_2) \\
& \ll_{j_1, j_2, \varepsilon } T^{3+2\varepsilon}  \left(T+|y_1|^{1/2} + |y_1|^{1/3}|y_2|^{1/6}   \right)^{j_1}  \left(T+ |y_2|^{1/2} + |y_2|^{1/3}|y_1|^{1/6}   \right)^{ j_2} 
\end{split}
 \end{equation}
 for all $j_1, j_2 \in \Bbb{N}_0$. 
\end{lemma}

Again the cut-off point $T^{1+o(1)}  \ll \Upsilon$ and the size of the oscillation are  sharp; the true order of magnitude of $\Phi_{w_6}$ is roughly $T^{3/2}$.\\  

We continue with strong bounds  for multiple Fourier transforms. The following two lemmas feature the smooth, compactly supported weight function $W$ from Section \ref{stage}. For  $\mu \in \Lambda_0$, and $\Xi$, $\Xi_1$, $\Xi_2$, $U$, $U_1$, $U_2$, $V$, $V_1$, $V_2 \in \Bbb{R}$ it is convenient to define the quantities
\begin{equation}\label{deftildeK}
\begin{split}
 &\tilde{\mathcal{K}}_{\mu}(\Xi, U, V):=   \int_{0}^{\infty}\int_{0}^{\infty} K_{w_4}( \xi \eta \Xi; \mu) e(\xi U+\eta V) W(\xi) W(\eta) d\xi\, d\eta
\end{split}
\end{equation}
and
\begin{equation}\label{defK}
\begin{split}
\mathcal{K}_{\mu}(\Xi_1, \Xi_2; & U_1, V_1; U_2, V_2):= \int_{\Bbb{R}^4} K^{\text{sgn}(\Xi_1), \text{sgn}(\Xi_2)}_{w_6}(   \xi_1 \eta_1 \Xi_1,  \xi_2 \eta_2 \Xi_2; \mu) \\
&\times e(\xi_1 U_1+\eta_1 V_1+\xi_2 U_2+\eta_2 V_2)  \   W(\xi_1) W(\eta_1)\overline{W}(\xi_2) \overline{W}(\eta_2)  d\xi_1\, d\xi_2\, d\eta_1\, d\eta_2.\end{split}
\end{equation}

\begin{lemma}\label{lem5}  Let   $U, V, \Xi  \in \Bbb{R}$, $|\Xi|, T > 1$.  Let $\mu \in \Lambda_0$  satisfy \eqref{munu}  
and let $\varepsilon, B > 0$.  \\ 
{\rm (a)} We have
$$\tilde{\mathcal{K}}_{\mu}(\Xi, U, 0) \ll_{\varepsilon, B} T^{-B}$$
unless $U \leq T^{\varepsilon}$. Similarly, $U = 0$ requires $V \leq T^{\varepsilon}$ for a non-negligible contribution.\\
{\rm (b)} We have
$$\tilde{\mathcal{K}}_{\mu}(\Xi, 0, 0) \ll_{\varepsilon, B} \begin{cases} T^{-B}, & |\Xi| \geq T^{3+\varepsilon},\\
T^{-3/2 + \varepsilon}, & \text{always.}\\
\end{cases}$$
{\rm (c)} If $|U|, |V| \geq T^{\varepsilon}$, then $\tilde{\mathcal{K}}_{\mu}(\Xi, U, V) \ll_{\varepsilon, B}  T^{-B}$ unless $U \asymp V$, in which case
$$\tilde{\mathcal{K}}_{\mu}(\Xi, U, V) \ll_{\varepsilon, B}  (|UV|^{1/2} + T)^{-3/2+\varepsilon}.$$ 
\end{lemma}


\begin{lemma}\label{lem6} Let  $U_1, V_1, U_2, V_2  \in \Bbb{R}$, $T, |\Xi_1|, |\Xi_2| > 1$, and assume that $\Xi_1, \Xi_2, U_1, U_2, V_1, V_2 \changed{\ll} T^{O(1)}$.   Let $\mu \in \Lambda_0$  satisfy \eqref{munu}. 
Let $\varepsilon, B > 0$. \\
{\rm (a)} We have
 $$\mathcal{K}_{\mu}(\Xi_1, \Xi_2; 0, V_1; U_2, V_2) \ll_{\varepsilon, B} T^{-B}$$
unless $V_1 \leq T^{\varepsilon}$. Similarly, $V_1 = 0$ requires $U_1 \leq T^{\varepsilon}$, $U_2 = 0$ requires $V_2 \leq T^{\varepsilon}$, and $V_2 = 0$ requires $U_2 \leq T^{\varepsilon}$ for a non-negligible contribution. \\
 {\rm (b)} If $|U_1|, |V_1| \geq T^{\varepsilon}$, we have   
  $$\mathcal{K}_{\mu}(\Xi_1, \Xi_2; U_1, V_1; 0, 0) \ll_{\varepsilon} (|U_1V_1|^{1/4} |\Xi_2|)^{-1+\varepsilon}. $$
  Similarly, if $|U_2|, |V_2| \geq T^{\varepsilon}$,  then  
  $\mathcal{K}_{\mu}(\Xi_1, \Xi_2; 0, 0; U_2, V_2) \ll (|U_2V_2|^{1/4} |\Xi_1|)^{-1+\varepsilon}. $\\
{\rm (c)} We have $\mathcal{K}_{\mu}(\Xi_1, \Xi_2; 0, 0, 0, 0) \ll_{\varepsilon, B} T^{-B}$ unless $\min(|\Xi_1|, |\Xi_2|) \geq T^{3-\varepsilon}$. In this case we have
\begin{equation}\label{c1}
\mathcal{K}_{ \mu}\left(\frac{  \Xi_1}{D}, \frac{  \Xi_2}{D}; 0, 0; 0, 0\right) \ll_{\varepsilon} T^{-3+\varepsilon}
\end{equation}
and 
\begin{equation}\label{c2}
\sum_{\epsilon \in \{\pm 1\}^2} \sum_{w \in \mathcal{W}} \sum_D \frac{\phi(D)}{D^2} \mathcal{K}_{w(\mu)}\left(\frac{\epsilon_1 \Xi_1}{D}, \frac{\epsilon_2 \Xi_2}{D}; 0, 0; 0, 0\right) \ll_{\varepsilon} |\Xi_1 \Xi_2|^{-1/2 + \varepsilon}.
\end{equation}
{\rm (d)} For $|U_1|, |V_1|, |U_2|, |V_2| \geq T^{\varepsilon}$ we define
\begin{equation}\label{defw}
\Upsilon_1 = \frac{|\Xi_1|}{ e^2 |U_1V_1|}, \quad \Upsilon_2 = \frac{|\Xi_2|}{  e^2 |U_2V_2|}.
\end{equation}
Then we 
have  
\begin{equation}\label{11cbound1}
|U_1V_1U_2V_2|^{1/2} \mathcal{K}_{\mu}(\Xi_1, \Xi_2; U_1, V_1; U_2, V_2) \ll_{\varepsilon}   T^{-1 +\varepsilon}
\end{equation}
if
$$ |\Upsilon_1 - 1| + \Upsilon_2 \ll T^{-\Halfbfrac} \ll   |U_2V_2|T^{-2}    \quad    \text{or} \quad  |\Upsilon_2 - 1| + \Upsilon_1 \ll T^{-\Halfbfrac} \ll  |U_1V_1|  T^{-2},$$
and 
\begin{equation}\label{11cbound2}
|U_1V_1U_2V_2|^{1/2} \mathcal{K}_{\mu}(\Xi_1, \Xi_2; U_1, V_1; U_2, V_2) \ll_{\varepsilon}   T^{-1 -\Eighthbfrac + \varepsilon}
\end{equation}
otherwise. \\
{\rm (e)} If $|U_1|, |V_1|, |U_2|, |V_2| \geq T^{\varepsilon}$ and in addition
$$|U_2V_2| \geq (|U_1V_1| + T^2)T^{\varepsilon},$$
then  $\mathcal{K}_{\mu}(\Xi_1, \Xi_2; U_1, V_1; U_2, V_2) \ll_{\varepsilon, B} T^{-B}$ unless 
\begin{equation}\label{cond11c}
|\Upsilon_2 - 1| \leq T^{\varepsilon} \left(\frac{1}{|U_2V_2|^{1/4}} + \frac{|U_1 V_1|^{1/2}+T}{|U_2V_2|^{1/2}}\right).
\end{equation}
Similarly,  $|U_1V_1| \geq (|U_2V_2| + T^2)T^{\varepsilon}$ requires
$|\Upsilon_1 - 1| \leq T^{\varepsilon} (|U_1V_1|^{-1/4} + (|U_2 V_2|^{1/2}+T)|U_1V_1|^{-1/2})$ 
for a non-negligible contribution. 
\end{lemma} 

\textbf{Remark.} Parts (a) -- (c) and (e) treat special configurations; in particular, part (c) will be used for the central term in the Poisson summation formula. Part (d) treats the generic situation. Its proof is by far the longest. The constant $\Eighthb$ could be improved at the cost of increasing the length of the paper. 
 
\section{The diagonal term}

We are now prepared to start with the proof of Theorem \ref{thm1}.  We   return to \eqref{basic} and apply the Kuznetsov formula as described in Section \ref{kuz} with  test function   $h$ as in \eqref{defh} and parameters
$$m_1s_1 \mapsto m_1, \quad n_2s_2 \mapsto m_2, \quad m_2r_2 \mapsto n_1, \quad n_1r_1 \mapsto n_2.$$
 In the following sections we estimate each of the four terms $\Delta$, $\Sigma_4$, $\Sigma_5$ and $\Sigma_6$ on the arithmetic side of the Kuznetsov formula. We recall the size conditions \eqref{lambda} and \eqref{sizeM} on $L$ and $M$. We start with the diagonal term.\\

By \eqref{vol} we have trivially
$$\Delta \ll \frac{T^{3+\varepsilon}}{M^2L^2}    \sum_{j=1}^3 \sum_{\substack{\ell_1, \ell_2 \asymp L\\ \ell_1, \ell_2 \text{ prime}}}  \sum_{\substack{r_0r_1r_2 = \ell_1^j \\ s_0s_1s_2 = \ell_2^j }} \sum_{\substack{r_2n_1 \asymp M\\ s_0m_1 \asymp M\\s_1n_2 \asymp M\\ r_0m_2 \asymp M}} \delta_{\substack{m_1s_1 = m_2r_2\\ n_2s_2 = n_1 r_1}}.$$
The diagonal contribution $\ell_1 = \ell_2 = \ell$, say,  is at most
\begin{displaymath}
\begin{split}
& \ll \frac{T^{3+\varepsilon}}{M^2L^2}   \sum_{j=1}^3 \sum_{ \ell  \asymp L }  \sum_{\substack{r_0r_1r_2 = \ell^j \\ s_0s_1s_2 = \ell^j }} 
  M^2 \ll \frac{T^{3+\varepsilon}}{L}.
 \end{split}
\end{displaymath}
For the off-diagonal term $\ell_1 \not= \ell_2$ we notice that $m_1$ must be divisible by $r_2$ and $m_2$ must be divisible by $s_1$. Hence we have $\min ( M/(s_0r_2),  M/(r_0s_1) )$ choices for  the pair $(m_1, m_2)$. Similarly we have   $\min (M/(r_1s_1),  M/(s_2r_2 ))$ choices for the pair  $(n_1, n_2)$. Using $\min(A, B) \leq \sqrt{AB}$, we see that 
  the off-diagonal contribution is at most
\begin{displaymath}
\begin{split}
\ll  \frac{T^{3+\varepsilon}}{M^2L^2}   \sum_{j=1}^3\sum_{ \ell_1, \ell_2 \asymp L }  \sum_{\substack{r_0r_1r_2 = \ell_1^j \\ s_0s_1s_2 = \ell_2^j }} 
 \frac{M^2}{(r_0r_1r_2^2s_0s_1^2s_2)^{1/2}} \ll \frac{T^{3+\varepsilon}}{L},
 \end{split}
\end{displaymath}
so that in total 
\begin{equation}\label{Delta}
\Delta \ll \frac{T^{3+\varepsilon}}{L}.
\end{equation}

\section{The $w_4$ and $w_5$ terms}
 By symmetry it is enough bound the $w_4$-term 
 $$\Sigma_4 \leq \frac{1}{L^2}\sum_{\epsilon = \pm 1} \sum_{j=1}^3\sum_{\substack{\ell_1, \ell_2 \asymp L\\ \ell_1, \ell_2 \text{ prime}}} \sum_{\substack{r_0r_1r_2 = \ell_1^j \\ s_0s_1s_2 = \ell_2^j }} |\Sigma_4(r, s)|,$$
 where 
\begin{displaymath}
\begin{split}
\Sigma_4(r, s) := \frac{T^{\varepsilon}}{M^2 }   &  \sum_{m_1, m_2, n_1, n_2}  W\Bigl(\frac{r_2n_1}{M}\Bigr) W\Bigl(\frac{s_0m_1}{M}\Bigr) \overline{W\Bigl(\frac{s_1n_2}{M}\Bigr) W\Bigl(\frac{r_0m_2}{M}\Bigr)} \\
& \times 
\sum_{\substack{D, \delta\\ n_2s_2 \delta = m_2r_2 D}} \frac{\tilde{S}(-\epsilon n_1r_1, n_2s_2, m_1s_1; D, D\delta)}{D^2\delta} \Phi_{w_4}\left( \frac{\epsilon n_1n_2m_1 r_1s_1s_2}{D^2\delta} \right). 
\end{split}
\end{displaymath}
By the argument of Section \ref{absconv} we can truncate the $D, \delta$ sum at some $T^B$ for some sufficiently large $B$ at the cost of a negligible error. Then by  
Lemma \ref{lem3}  we can truncate the sum, again with a  negligible error, at
$$\left( \frac{  n_1n_2m_1 r_1s_1s_2}{D^2\delta}\right)^{1/3} \geq T^{1-\varepsilon},$$
or in other words
$$ D^{2/3}\delta^{1/3} \leq \frac{M}{T^{1-\varepsilon}} \left(\frac{r_1s_2}{r_2s_0}\right)^{1/3}.$$
We apply Poisson summation in the $m_1$, $n_1$ variables and estimate the remaining sums trivially.  By \eqref{Klo5}  this gives
\begin{displaymath}
\begin{split}
\Sigma_4(r, s) \ll &\frac{T^{\varepsilon}}{M^2 }   
  \sum_{D^2 \delta \leq T^{\varepsilon} \frac{M^3r_1s_2}{T^3r_2s_0}}  \frac{1}{D^2 \delta} \sum_{\substack{ r_0m_2 \asymp M\\ s_1n_2 \asymp M\\ n_2s_2 \delta = m_2r_2 D}} 
  \sum_{\substack{x, y \in \Bbb{Z}\\ (D, x) = (r_1, x)\\ (\delta, y) = (s_1, y)}} D(r_1,  D) (s_1,   \delta) \bigl|\widehat{\Phi}_{D, \delta, n_2, r, s}(x, y)\bigr|,
\end{split} 
\end{displaymath}
where
$$\widehat{\Phi}_{D, \delta, n_2, r, s}(x, y) = \frac{M^2}{r_2s_0}\int_{0}^{\infty} \int_0^{\infty} \Phi_{w_4}\left(\frac{\epsilon \xi \eta n_2  r_1s_1s_2  M^2}{D^2 \delta r_2s_0 }\right) e\left(\frac{\xi xM}{r_2D}\right) e\left(\frac{\eta yM}{s_0 \delta}\right) W(\xi) W(\eta)  d\xi \,d \eta.$$
Integration by parts  using \eqref{w42} and the condition $s_1n_2 \asymp M$ as well as the bound on $D^2\delta$ shows  that this is negligible unless
\begin{displaymath}
\begin{split}
&|x| \leq X := T^{\varepsilon} \left(\frac{M^3 r_1s_2}{r_2s_0D^2 \delta}\right)^{1/3} \frac{r_2D}{M} = T^{\varepsilon}\left(\frac{r_1r_2^2s_2D}{ s_0\delta}\right)^{1/3},\\
& |y| \leq Y := T^{\varepsilon} \left(\frac{M^3 r_1s_2}{r_2s_0D^2 \delta}\right)^{1/3} \frac{s_0\delta}{M} = T^{\varepsilon}\left(\frac{r_1s_2s_0^2\delta^2}{r_2D^2}\right)^{1/3}.
\end{split}
\end{displaymath} 
We insert the definition \eqref{defPhi} for $\Phi_{w_4}$ and pull the $\mu$-integral outside. By the properties of $h$ and \eqref{vol}  we are left with bounding
 \begin{equation*}
\begin{split}
\Sigma_4(r, s; \mu) :=  \frac{T^{3+\varepsilon}}{r_2s_0} \sum_{D^2 \delta \leq  T^{\varepsilon} \frac{M^3r_1s_2}{T^3r_2s_0}}  \frac{Dr_1s_1  }{D^2 \delta}  \sum_{\substack{ r_0m_2 \asymp M\\ s_1n_2 \asymp M\\ n_2s_2 \delta = m_2r_2 D}} \sum_{\substack{|x| \leq X, |y| \leq Y\\ (D, x) = (r_1, x)\\ (\delta, y) = (s_1, y)}} \Big|\tilde{\mathcal{K}}_{\mu}\left(\frac{\epsilon   n_2  r_1s_1s_2  M^2}{D^2 \delta r_2s_0 }, \frac{ xM}{r_2D}, \frac{  yM}{s_0 \delta} \right)\Big|
\end{split}
\end{equation*}
with $\tilde{\mathcal{K}}_{\mu}(\Xi, U, V)$ as in \eqref{deftildeK}, where $\mu \in \Lambda_0$ satisfies \eqref{munu}.  Notice that our summation conditions imply 
$$\frac{  n_2  r_1s_1s_2  M^2}{D^2 \delta r_2s_0 } \gg T^{3-\varepsilon},$$
so that  the condition $|\Xi|  \geq 1$ of Lemma \ref{lem5} is satisfied.  As usual  in the Poisson summation formula, the central terms need special treatment. The summation conditions imply
\begin{displaymath}
\begin{split}
\min \left(\frac{M}{r_2 D}, \frac{M}{s_0\delta}\right)  & = \min \left(\frac{M(m_2r_2)^{1/3}}{r_2 D^{2/3} (n_2s_2\delta)^{1/3}}, \frac{M(n_2s_2)^{2/3}}{s_0\delta^{1/3}(m_2r_2D)^{2/3}}\right) \\
&\gg T^{1-\varepsilon} \min\left(\frac{ s_0s_1}{r_0r_1r_2s_2^2}, \frac{s_2r_0^2}{s_0^2s_1^2r_1r_2}\right)^{1/3} \gg \frac{T^{1-\varepsilon}}{L^{j}} \geq T^{4/5}
\end{split}
\end{displaymath}
by \eqref{lambda}. 
Thus we conclude from Lemma \ref{lem5}(a) that $xy = 0$ implies $x = y = 0$, up to a negligible error. 

 We start with the contribution $\Sigma_4^{0}(r, s; \mu)$ of the  central term $x = y = 0$ to $\Sigma_4(r, s; \mu)$.  The summation conditions imply $D = r_1$, $\delta = s_1$. The condition $\Xi \leq T^{3 + \varepsilon}$ from Lemma \ref{lem5}(b) implies
  $$\frac{  n_2  s_2  M^2}{r_1  r_2s_0 } \leq T^{3+\varepsilon}.$$
 Using again   Lemma \ref{lem5}(b)   together with a well-known divisor argument, we obtain
 \begin{equation}\label{w4central}
\begin{split}
\Sigma_4^{0}(r, s; \mu) & \ll  \frac{T^{3+\varepsilon} }{r_2s_0}       \sum_{\substack{ s_2 n_2 \leq r_1r_2s_0  T^{3+\varepsilon}M^{-2}    \\ n_2s_2 s_1 = m_2r_2 r_1}}  T^{-\frac{3}{2}}  \ll \frac{T^{\frac{9}{2}+\varepsilon}r_1}{M^2} \ll T^{\frac{3}{2} + 3\lambda + 2\eta + \varepsilon}. 
\end{split}
\end{equation}

We proceed to bound the terms $xy \not= 0$ in $\Sigma_4(r, s; \mu)$, say $\Sigma_4^{\ast}(r, s; \mu)$.  First we observe that the summation condition implies 
 \begin{equation}\label{boundD}
\frac{s_1r_2D}{r_0s_2} \asymp \frac{m_2 r_2D}{n_2s_2} \asymp \delta \quad \text{so that} \quad D^3 \frac{s_1r_2}{r_0s_2}  \asymp D^2 \delta \ll T^{\varepsilon} \frac{M^3 r_1s_2}{T^3 r_2s_0}. 
\end{equation}
 Moreover, fixing $D$ and $m_2$ determines $n_2$ and $\delta$, up to a divisor function.  According to Lemma \ref{lem5}(c) and \eqref{boundD} we 
 now obtain
 \begin{equation}\label{w4main}
 \begin{split}
   \Sigma^{\ast}_4(r, s; \mu) & \ll  T^{3+\varepsilon} 
 \frac{r_1s_1 }{r_2s_0} \sum_{D^2 \delta \leq T^{\varepsilon} \frac{M^3r_1s_2}{T^3r_2s_0}}  \frac{1 }{D \delta} \sum_{\substack{ r_0m_2 \asymp M\\ s_1n_2 \asymp M\\ n_2s_2 \delta = m_2r_2 D}} \sum_{\substack{0 < |x| \leq X\\ 0 <  |y| \leq Y }}\left(\frac{D \delta r_2s_0}{|xy|M^2}\right)^{\frac{3}{4}}\\
   & \ll  T^{3+\varepsilon} \frac{r_1s_1 }{r_2s_0} \sum_{D^2 \delta \leq T^{\varepsilon} \frac{M^3r_1s_2}{T^3r_2s_0}}  \frac{(r_2s_0)^{\frac{3}{4}} }{(D \delta)^{\frac{1}{4}}M^{\frac{3}{2}}} \sum_{\substack{ r_0m_2 \asymp M\\ s_1n_2 \asymp M\\ n_2s_2 \delta = m_2r_2 D}}    \left(\frac{ \delta r_1^2 r_2s_0 s_2^2  }{ D}\right)^{\frac{1}{12}}\\
     & \ll   T^{3+\varepsilon}  \frac{r_1s_1 }{r_2s_0} \sum_{D^3  \ll   T^{\varepsilon} \frac{M^3 r_0r_1s_2^2}{T^3 r_2^2s_0s_1}}  \frac{(r_2s_0)^{\frac{3}{4}} }{D^{\frac{1}{3}} (D\frac{s_1r_2}{r_0s_2})^{\frac{1}{6}} M^{\frac{3}{2}}} \sum_{ m_2 \asymp   M/r_0}   (r_1^2r_2s_0s_2^2)^{\frac{1}{12}}\\
   & \ll   T^{3+\varepsilon} 
 \frac{r_1s_1 }{r_2s_0} \frac{(r_2s_0)^{3/4}}{(\frac{s_1r_2}{r_0s_2})^{\frac{1}{6}}M^{\frac{3}{2}}}   (r_1^2r_2s_0s_2^2)^{\frac{1}{12}}  \left( \frac{M^3 r_0r_1s_2^2}{T^3 r_2^2s_0s_1} \right)^{\frac{1}{6}}  \frac{M}{r_0}  =  T^{\frac{5}{2}+\varepsilon}  \left(\frac{r_0^2r_1^4 s_1^2 s_2^2}{r_2^2 s_0}\right)^{1/3},
      \end{split}
 \end{equation}
 which is at most $T^{5/2+\varepsilon} L^6$ for $r_0r_1r_2, s_0s_1s_2 \ll L^3$.  Combining \eqref{w4central} and \eqref{w4main}, we obtain
\begin{equation}\label{w4final}
  \Sigma_4 \ll T^{\frac{3}{2} + 3\lambda + 2\eta + \varepsilon} + T^{\frac{5}{2}+ 6\lambda + \varepsilon}. \end{equation}

\section{The $w_6$-term}


\subsection{Truncation and Poisson summation} We have
\begin{displaymath}
\begin{split}
\Sigma_6 = &\frac{T^{\varepsilon}}{M^2L^2}   \sum_{j=1}^3\sum_{\substack{\ell_1, \ell_2 \asymp L\\ \ell_1, \ell_2 \text{ prime}}} \sum_{\substack{r_0r_1r_2 = \ell_1^j \\ s_0s_1s_2 = \ell_2^j }}\Bigl| \sum_{\epsilon \in \{\pm 1\}^2} \sum_{m_1, m_2, n_1, n_2}  W\Bigl(\frac{r_2n_1}{M}\Bigr) W\Bigl(\frac{s_0m_1}{M}\Bigr) \overline{W\Bigl(\frac{s_1n_2}{M}\Bigr) W\Bigl(\frac{r_0m_2}{M}\Bigr)}  \\
&  \times  \sum_{D_1, D_2} 
\frac{S(\epsilon_2 n_1r_1, \epsilon_1 m_2r_2, m_1s_1, n_2s_2; D_1, D_2)}{D_1D_2} \Phi_{w_6} \left(\frac{-\epsilon_2 n_1m_1s_1r_1D_2}{D_1^2}, \frac{ -\epsilon_1n_2m_2D_1s_2r_2}{D_2^2} \right)\Bigr|.
\end{split}
\end{displaymath}
It is absolutely crucial to keep the $\epsilon$-sum inside the absolute values.
By the argument of Section \ref{absconv} we can truncate the $D_1, D_2$-sum at some $T^B$ for some sufficiently large $B$ at the cost of a negligible error. Then 
by \eqref{J1} we can truncate the sum further at
$$\frac{(n_2m_2s_2r_2)^{1/3} (n_1m_1s_1r_1)^{1/6}}{D_2^{1/2}} \geq T^{1-\varepsilon}, \quad \frac{(n_1m_1s_1r_1)^{1/3} (n_2m_2s_2r_2)^{1/6}}{D_1^{1/2}} \geq T^{1-\varepsilon},$$
or in other words
\begin{equation}\label{D}
D_1 \leq T^{\varepsilon} \frac{M^2}{T^2} \left(\frac{r_1^2s_1s_2}{r_0r_2 s_0^2}\right)^{1/3} \ll \Deltavar := \frac{M^2L^j}{T^{2-\varepsilon}}, \quad D_2 \leq T^{\varepsilon} \frac{M^2}{T^2} \left(\frac{r_1r_2s_2^2}{r_0^2 s_0s_1}\right)^{1/3}  \ll \Deltavar.
\end{equation}
We apply Poisson summation to all four variables $n_1, n_2, m_1, m_2$. This gives
\begin{displaymath}
\begin{split}
\Sigma_6 \ll &\frac{T^{\varepsilon}}{L^2}    \sum_{j=1}^3\sum_{\substack{\ell_1, \ell_2 \asymp L\\ \ell_1, \ell_2 \text{ prime}}} \sum_{\substack{r_0r_1r_2 = \ell_1^j \\ s_0s_1s_2 = \ell_2^j }} \frac{M^2}{r_0r_2s_0s_1} \Bigl| \sum_{\epsilon \in \{\pm 1\}^2}  \sum_{\changed{D_1, D_1 \leq \Deltavar} }  \sum_{\substack{x_1, x_2  \in \Bbb{Z}\\ y_1, y_2 \in \Bbb{Z}}} \widehat{S}_{\epsilon_2r_1, s_1, \epsilon_1 r_2, s_2}(x_1, y_1, x_2, y_2; D_1, D_2)\\
&  \times \int_{\Bbb{R}^4} \Phi_{w_6} \left( \frac{-\epsilon_2 M^2 \eta_1\xi_1s_1r_1D_2}{D_1^2 r_2s_0}, \frac{-\epsilon_1 M^2 \eta_2\xi_2D_1s_2r_2}{D_2^2 s_1r_0} \right) e\left(\frac{x_1\xi_1M}{D_1s_0} + \frac{y_1 \eta_1M}{D_1r_2} + \frac{x_2\xi_2M}{D_2r_0}  + \frac{y_2 \eta_2M}{D_2s_1}\right) \\
& \quad\quad\quad \times W(\xi_1) W(\eta_1)\overline{W}(\xi_2) \overline{W}(\eta_2) d\xi_1\, d\xi_2\, d\eta_1\, d\eta_2 \Bigr|. 
\end{split}
\end{displaymath}
 Integration by parts in connection with \eqref{J2} shows that the integral is negligible unless
\begin{displaymath}
  \changed{|x_1|}  \leq T^{\varepsilon} \left(\frac{M\sqrt{ s_1r_1D_2}}{D_1\sqrt{r_2s_0}} + \frac{M(s_1s_2r_1^2)^{1/6}}{D_1^{1/2}(r_0r_2s_0^2)^{1/6}} \right)  \cdot \frac{D_1s_0}{M} \ll T^{\varepsilon} L^j(D_1 + D_2)^{1/2} =:   X,
\end{displaymath}
and similarly, $|y_1|, |x_2|, |y_2| \leq T^{\varepsilon} \Xvar$. 
 We insert the definition \eqref{defPhi} and sort the integration over $\mu$ by Weyl chambers. Recalling that  $h$ is Weyl-group invariant, we   pull the integration over one Weyl  chamber outside, leaving the sum over the Weyl group inside.   By the properties of $h$ and \eqref{vol}  we are left with bounding
 \begin{displaymath}
\begin{split}
 \Sigma_6(\mu) := &\frac{T^{3+\varepsilon}M^2}{L^2}   \sum_{j=1}^3\sum_{\substack{\ell_1, \ell_2 \asymp L\\ \ell_1, \ell_2 \text{ prime}}} \sum_{\substack{r_0r_1r_2 = \ell_1^j \\ s_0s_1s_2 = \ell_2^j }}\\
 &  \Bigl|   \sum_{w \in \mathcal{W}} \sum_{\epsilon \in \{\pm 1\}^2} \sum_{\changed{D_1, D_2 \leq \Deltavar} }  \sum_{\changed{|x_1|, |x_2| ,  |y_1|, |y_2| \leq  \Xvar}} \frac{\widehat{S}_{\epsilon_2r_1, s_1, \epsilon_1 r_2, s_2}(x_1, y_1, x_2, y_2; D_1, D_2)}{r_0r_2s_0s_1D_1D_2}\\
&\times   \mathcal{K}_{w(\mu)}\left(\frac{-\epsilon_2 M^2  s_1r_1D_2}{D_1^2 r_2s_0}, \frac{-\epsilon_1 M^2 D_1s_2r_2}{D_2^2 s_1r_0};  \frac{x_1 M}{D_1s_0}, \frac{y_1  M}{D_1r_2} ; \frac{x_2 M}{D_2r_0}, \frac{y_2 M}{D_2s_1} \right)\Big|
\end{split}
\end{displaymath}
with $\mathcal{K}_{\mu}(\Xi_1, \Xi_2; U_1, V_1; U_2, V_2)$ as in \eqref{defK}, where $\mu \in \Lambda_0$ satisfies \eqref{munu}. The first two arguments of $K_{\mu}$  satisfy
$$\frac{ M^2  s_1r_1D_2}{D_1^2 r_2s_0}, \frac{ M^2 D_1s_2r_2}{D_2^2 s_1r_0}  \gg \frac{T^{1-\varepsilon}}{L^{4j}}$$
by \eqref{D} and \eqref{sizeM}, so that the condition $|\Xi_1|, |\Xi_2| \geq 1$ of Lemma \ref{lem6} is satisfied by \eqref{lambda}. 
Since
$$\min\left(\frac{M}{D_1s_0}, \frac{M}{D_1r_2}\right) \gg \frac{T^{1/2 - \varepsilon}}{L^{2j}} \geq T^{1/6}$$
 by \eqref{lambda},   
 we conclude similarly as in the preceding section  from Lemma \ref{lem6}(a) that $x_1y_1 = 0$ implies $x_1 = y_1 = 0$ (otherwise the contribution is negligible).  Similarly, $x_2y_2 = 0$ implies $x_2 = y_2 = 0$. 

\subsection{The central term}\label{centralsec} We start with bounding the contribution $\Sigma^{0}_6(\mu)$ of the terms $x_1 =x_2 = y_1 = y_2 = 0$. We consider first the terms $\Sigma^{0, \not=}_6(\mu)$ with $\ell_1 \not= \ell_2$, in which case in particular $(r_1r_2, s_1s_2) = 1$.  By Lemma \ref{lem2}(c)  we obtain
\begin{displaymath}
\begin{split}
 \Sigma^{0, \not=}_{6}(\mu) = &\frac{T^{3+\varepsilon}M^2}{L^2}   \sum_{j=1}^3\sum_{\substack{\ell_1, \ell_2 \asymp L\\ \ell_1 \not= \ell_2 \text{ prime}}} \sum_{\substack{r_0r_1r_2 = \ell_1^j \\ s_0s_1s_2 = \ell_2^j }} \Bigl|   \sum_{w \in \mathcal{W}} \sum_{\epsilon \in \{\pm 1\}^2} \sum_{D  \leq \Deltavar}  \frac{\phi(D)}{r_0r_2s_0s_1D^2}\\
&\times   \mathcal{K}_{w(\mu)}\left(\frac{-\epsilon_2 M^2  s_1r_1}{D r_2s_0}, \frac{-\epsilon_1 M^2 s_2r_2}{D  s_1r_0}; 0, 0; 0, 0 \right)\Big|. 
\end{split}
\end{displaymath}
We conclude from Lemma \ref{lem6}(c)  and \eqref{sizeM} that   $\mathcal{K}_{\mu}$ is negligible for $D \gg L^{6+\varepsilon}$. In particular, by   \eqref{lambda}  and \eqref{D} we can complete the $D$-sum at the cost of a negligible error. Applying \eqref{c2}, we obtain
\begin{displaymath}
\begin{split}
 \Sigma^{0, \not=}_{6}(\mu) \ll &\frac{T^{3+\varepsilon}M^2}{L^2}   \sum_{j=1}^3\sum_{\substack{\ell_1, \ell_2 \asymp L\\ \ell_1 \not= \ell_2 \text{ prime}}} \sum_{\substack{r_0r_1r_2 = \ell_1^j \\ s_0s_1s_2 = \ell_2^j }}\frac{1}{r_0r_2s_0s_1 }  \left(\frac{ M^2  s_1r_1}{ r_2s_0} \cdot \frac{M^2 s_2r_2}{  s_1r_0}\right)^{-1/2} \ll \frac{T^{3+\varepsilon}}{L}.
\end{split}
\end{displaymath}
We proceed to bound the   contribution  $\Sigma^{0, =}_{6}(\mu) $ of the terms $\ell_1 = \ell_2  = \ell$, say.  We write $D_j = \ell^{\alpha_j} D_j'$ for $j= 1, 2$ with $\ell \nmid D_j'$. By Lemma \ref{lem2}(d) and \eqref{c1} we obtain
\begin{displaymath}
\begin{split}
 \Sigma^{0, \not=}_{6}(\mu)  & \ll \frac{T^{3+\varepsilon}M^2}{L^2}   \sum_{j=1}^3\sum_{\ell \asymp L \text{ prime}} \sum_{\substack{r_0r_1r_2 = \ell^j \\ s_0s_1s_2 = \ell^j }}\frac{1}{r_0r_2s_0s_1 }   \sum_{\alpha_1, \alpha_2 \leq T^{\varepsilon}} \sum_{\changed{D'_1,D'_2 \leq \Deltavar}} \frac{(D_1', D_2') \ell^{\alpha_1 + \alpha_2}r_2}{\ell^{\alpha_1+\alpha_2}D'_1D'_2T^{3-\varepsilon}} 
 \ll \frac{T^{3+\varepsilon}}{L}.
 \end{split}
\end{displaymath}
Combining the previous two displays we obtain
\begin{equation}\label{central}
   \Sigma^{0}_6(\mu) \ll \frac{T^{3+\varepsilon}}{L}.  
\end{equation}

\subsection{The mixed terms} Next we consider the contribution $\Sigma^{\text{mix}}_6(\mu)$  of the terms $x_1 = y_1 = 0 \not= x_2y_2$. (By symmetry, the same argument works for $x_2 = y_2 = 0 \not= x_1 y_1$.) From now on we can sum trivially over $\epsilon \in \{\pm 1\}^2$ and $w \in \mathcal{W}$. We conclude from Lemma \ref{lem2}(a) and (b) and from Lemma \ref{lem6}(b) that 
\begin{displaymath}
\begin{split}
  \Sigma^{\text{mix}}_6(\mu) & \ll \frac{T^{3+\varepsilon}M^2}{L^2}   \sum_{j=1}^3  \sum_{\ell_1, \ell_2 \asymp L }  \sum_{\substack{r_0r_1r_2 = \ell_1^j \\ s_0s_1s_2 = \ell_2^j }} 
  \sum_{\substack{\changed{D_1, D_2 \leq \Deltavar} \\ D_1 \mid r_1s_1 D_2} }  \sum_{\substack{ 0 < |x_2y_2| \leq  \changed{\Xvar^2}\\ x_2 y_2 \equiv r_2s_2 D_1 \, (\text{mod} D_2)\\ D_1/(D_1, r_1s_1) \mid (x_2, y_2)}} \frac{(r_1, D_1)(r_2, D_2)(D_1, D_2) }{r_0r_2s_0s_1D_1D_2}\\
  & \times \frac{(r_0s_1)^{1/4}D_2^{1/2}}{|x_2y_2|^{1/4} M^{1/2}} \cdot \frac{D_1^2 r_2s_0}{M^2 s_1r_1D_2}. 
\end{split}
\end{displaymath}
We write $F := (D_1, r_1s_1) \leq r_1s_1$,   introduce the variable $$z := x_2y_2F^2D_1^{-2} \in \Bbb{Z} \setminus\{0\},$$
which, up to a divisor function, determines $x_2, y_2$, and we write 
$D_1 D = D_2F$ 
with $D \in \Bbb{N}$. In this way we obtain that $\Sigma_6^{\text{mix}}(\mu)$ is at most
\begin{displaymath}
\begin{split}
 \frac{T^{3+\varepsilon}M^2}{L^2}   \sum_{j=1}^3 \sum_{ \ell_1, \ell_2 \asymp L }  \sum_{\substack{r_0r_1r_2 = \ell_1^j \\ s_0s_1s_2 = \ell_2^j }} 
  \sum_{ DD_1 \leq \Deltavar r_1s_1 }  \sum_{\substack{ 0 < |z| \leq  \changed{\Xvar^2} (F/D_1)^2\\ zD_1/F \equiv r_2s_2 F \, (\text{mod} D)}} \frac{r_1r_2F}{r_0r_2s_0s_1DD_1} \cdot \frac{r_0^{\frac{1}{4}} r_2s_0D_1F}{s_1^{\frac{3}{4}} r_1M^{\frac{5}{2}} D^{\frac{1}{2}}|z|^{\frac{1}{4}}}.\\
\end{split}
\end{displaymath}
We can afford to   drop the factor $|z|^{-1/4}$. The summation condition on $z$ implies
$$1 \leq  \changed{\Xvar^2} \left(\frac{F}{D_1}\right)^2 \ll T^{\varepsilon}\frac{L^{2j}}{D_1}  (DF + F^2)\ll T^{\varepsilon} \frac{L^{6j}D}{D_1}  , $$
and $z$ is determined modulo $D/(D, r_2s_2s_1r_1)$, so that there are at most $1 + T^{\varepsilon} L^{8j} D_1^{-1}$ choices for $z$.  We summarize
\begin{displaymath}
\begin{split}
  \Sigma^{\text{mix}}_6(\mu) & \ll \frac{T^{3+\varepsilon}M^2}{L^2}   \sum_{j=1}^3 \sum_{ \ell_1, \ell_2 \asymp L }  \sum_{\substack{r_0r_1r_2 = \ell_1^j \\ s_0s_1s_2 = \ell_2^j }} 
 \sum_{ \substack{DD_1 \leq \Deltavar r_1s_1\\ D_1 \ll T^{\varepsilon} D L^{6j}}} \left( 1 + \frac{L^{8j}}{D_1}\right) \frac{r_2 F^2}{r_0^{3/4}s_1^{7/4} D^{3/2} M^{5/2}}  \\
  & \ll \frac{T^{3+\varepsilon}M^2}{L^2}   \sum_{j=1}^3 \sum_{ \ell_1, \ell_2 \asymp L }  \sum_{\substack{r_0r_1r_2 = \ell_1^j \\ s_0s_1s_2 = \ell_2^j }} 
 \frac{r_2 r_1^2s_1^{1/4}}{r_0^{3/4}   M^{5/2}}  \left(L^{8j} + (\Deltavar r_1s_1)^{1/4} L^{\frac{3}{4} \cdot 6j}\right).
 \end{split}
 \end{displaymath}
 In the last term we summed over $D_1 \ll T^{\varepsilon}\min(\Deltavar r_1s_1/D, D L^{6j})$ first and then over $D$. This gives 
 \begin{equation}\label{mix}
 \begin{split}
 \Sigma^{\text{mix}}_6(\mu)  & \ll \sum_{j=1}^3 \left( \frac{T^{3+\varepsilon} L^{41j/4} }{M^{1/2}}  + T^{5/2 + \varepsilon} L^{29j/4} \right) \ll T^{\frac{9}{4} + \frac{\eta}{2} + 31\lambda +\varepsilon} + T^{\frac{5}{2} +   22\lambda+\varepsilon}. 
\end{split}
\end{equation}

 
\subsection{The generic terms} Finally, we bound the contribution $\Sigma_6^{\text{gen}}(\mu)$ of the terms $x_1 y_1 x_2 y_2 \not= 0$.  
From Lemma \ref{lem2}(a)  we obtain the congruences $x_1y_1 \equiv  r_1s_1 D_2 \, (\text{mod } D_1)$, $x_2y_2 \equiv  r_2s_2 D_1 \, (\text{mod } D_2)$, which we re-write as
$$ x_1y_1=   r_1s_1 D_2 + c_1D_1, \quad x_2y_2=   r_2s_2 D_1 + c_2D_2$$
with $c_1, c_2 \in \Bbb{Z}$. This gives
$$\Sigma_6^{\text{gen}}(\mu) \ll  \frac{1}{L^2}   \sum_{\epsilon \in \{\pm 1\}^2}  \sum_{j=1}^3\sum_{ \ell_1, \ell_2 \asymp L }  \sum_{\substack{r_0r_1r_2 = \ell_1^j \\ s_0s_1s_2 = \ell_2^j }} |\Sigma_6^{\text{gen}}(r, s; \mu)|,$$
where
 \begin{displaymath}
\begin{split}
 \Sigma_6^{\text{gen}}(r, s; \mu) := & T^{3+\varepsilon}M^2 
\underset{\substack{x_1y_1=   r_1s_1 D_2 + c_1D_1\\ x_2y_2=   r_2s_2 D_1 + c_2D_2}}{\sum_{c_1, c_2 \in \Bbb{Z}} \sum_{\changed{D_1, D_2 \leq \Deltavar} }  \sum_{\substack{0 < |x_1y_1| \leq  \changed{\Xvar^2}\\  0<  |x_2y_2| \leq   \changed{\Xvar^2}}}} \frac{r_1r_2(D_1, D_2)}{r_0r_2s_0s_1D_1D_2}\\
&\times  \Big|\mathcal{K}_{\mu}\left(\frac{-\epsilon_2 M^2  s_1r_1D_2}{D_1^2 r_2s_0}, \frac{-\epsilon_1 M^2 D_1s_2r_2}{D_2^2 s_1r_0} ;  \frac{x_1 M}{D_1s_0}, \frac{y_1  M}{D_1r_2}; \frac{x_2 M}{D_2r_0}, \frac{y_2 M}{D_2s_1}  \right)\Big|.
\end{split}
\end{displaymath}
 Introducing the variables $D = (D_1, D_2)$ and $z_1D = x_1y_1$, $z_2 D = x_2y_2$, we obtain 
 \begin{displaymath}
\begin{split}
 \Sigma_6^{\text{gen}}(r, s; & \mu)\ll T^{3+\varepsilon}M^2 \frac{r_1  }{r_0 s_0s_1 }  \sum_{D \leq \Deltavar} \underset{\substack{z_1=   r_1s_1 D_2 + c_1D_1\\ z_2=   r_2s_2 D_1 + c_2D_2}}{\sum_{c_1, c_2 \in \Bbb{Z}} \sum_{\changed{D_1,D_2 \leq \Deltavar/D} }  \sum_{\substack{0 < |z_1| \changed{\ll} T^{\varepsilon}  (D_1 + D_2) L^{2j}\\  0<  |z_2| \changed{\ll} T^{\varepsilon}  (D_1+ D_2)L^{2j}}}} \frac{1 }{DD_1D_2}\\
&\times  \sum_{\substack{x_1y_1 = Dz_1\\ x_2y_2 = Dz_2}} \Big|\mathcal{K}_{\mu}\left( \frac{-\epsilon_2 M^2  s_1r_1D_2}{DD_1^2 r_2s_0}, \frac{-\epsilon_1 M^2 D_1s_2r_2}{DD_2^2 s_1r_0};   \frac{x_1 M}{DD_1s_0}, \frac{y_1  M}{DD_1r_2}; \frac{x_2 M}{DD_2r_0}, \frac{y_2 M}{DD_2s_1} \right)\Big|,
\end{split}
\end{displaymath}
and we want to apply Lemma \ref{lem6}(d) and (e)  with
\begin{equation}\label{variables}
\begin{split}
 |U_1V_1|   = \frac{|z_1| M^2}{DD_1^2s_0r_2}, \quad |U_2V_2|   = \frac{|z_2| M^2}{DD_2^2r_0s_1} , \quad   \Upsilon_1 = \frac{r_1s_1D_2}{ e^2|z_1|}, \quad \Upsilon_2 = \frac{r_2s_2D_1}{  e^2|z_2|}.
\end{split}
\end{equation}
We now distinguish the three cases $c_1 = c_2 = 0$, $c_1 c_2 = 0$ but $(c_1, c_2) \not= 0$, and $c_1 c_2 \not= 0$. We call the corresponding contributions $\Sigma_6^{\text{gen}, 0}(r, s; \mu)$, $\Sigma_6^{\text{gen}, \text{mix}}(r, s; \mu)$ and $\Sigma_6^{\text{gen}, \ast}(r, s; \mu)$, respectively. \\

\emph{Case I:} If $c_1 = c_2 = 0$,  then $\Upsilon_1 = \Upsilon_2 = e^{-2}  $, so that the condition $|\Upsilon_1 - 1| + \Upsilon_2 \ll T^{-\Halfb}$ or $|\Upsilon_2- 1| + \Upsilon_1 \ll T^{-\Halfb}$ from Lemma \ref{lem6}(d) can never happen.  
By \eqref{11cbound2} and \eqref{D} we conclude 
\begin{equation}\label{gen0}
\begin{split}
\Sigma_6^{\text{gen}, 0}(r, s; \mu) 
& \ll   T^{2 - \Eighthbfrac+\varepsilon} 
\frac{r_1}{r_0s_0s_1} \sum_{D \leq \Deltavar}\sum_{\changed{D_1,D_2 \leq \Deltavar/D} }   \frac{   \sqrt{r_0s_1s_0r_2}}{  \sqrt{r_1s_1D_2 r_2s_2D_1}} \\
& = T^{2 - \Eighthbfrac+\varepsilon} 
\frac{r_1^{1/2}}{(r_0s_1)^{1/2} (s_0s_1s_2)^{1/2} } \sum_{D \leq \Deltavar}\sum_{\changed{D_1,D_2 \leq \Deltavar/D} }   \frac{   1}{  \sqrt{ D_2  D_1}}
 \ll T^{3 - \Eighthbfrac + \varepsilon} L^3. 
\end{split} 
\end{equation}\\

\emph{Case  II:} If $c_1c_2 \not= 0$, we     insert  again the bound from Lemma \ref{lem6}(d) getting
\begin{displaymath}
\begin{split}
\Sigma_6^{\text{gen}, \ast}(r, s; \mu) \ll   T^{2+\varepsilon}   &
\frac{r_1r_2^{1/2} }{(r_0s_0s_1)^{1/2} }  \sum_{D \leq \Deltavar}    \underset{\substack{z_1=   r_1s_1 D_2 + c_1D_1\\ z_2=   r_2s_2 D_1 + c_2D_2}}{ \sum_{\substack{\changed{D_1,D_2 \leq \Deltavar/D}\\ c_1c_2 \not= 0}}  \sum_{\substack{0 < |z_1| \leq T^{\varepsilon}\max(D_1, D_2) L^{2j} \\ 0 < |z_2| \changed{\leq} T^{\varepsilon} \max(D_1, D_2) L^{2j} }}} \frac{1}{|z_1z_2|^{1/2}}\\
& \times \left(T^{-\Eighthbfrac} + \delta_{\substack{|\Upsilon_2 -1| + \Upsilon_1 \ll T^{-\Halfbfrac}\\ |U_1V_1| \gg T^{ 2- \Halfbfrac}}}+ \delta_{\substack{|\Upsilon_1 - 1| + \Upsilon_2 \ll T^{-\Halfbfrac}\\ |U_2V_2| \gg T^{2- \Halfbfrac}}} \right). 
 \end{split} 
\end{displaymath}
We split the variables into dyadic intervals
$$D_1 \asymp \mathcal{D}_1, \quad D_2 \asymp \mathcal{D}_2, \quad |z_1| \asymp Z_1, \quad |z_2| \asymp Z_2,$$
where $\changed{\mathcal{D}_1, \mathcal{D}_2 \ll \Deltavar/D}$ and $Z_1, Z_2 \ll T^{\varepsilon}L^{2j} \Deltavar/D.$ 
If $Z_1 \leq Z_2$, we choose first $z_1$ and $D_2$. This determines, up to a divisor function, $c_1$ and $D_1$. Having all of these fixed, we choose $c_2$, which determines the last variable $z_2$. The number of choices for $z_1$, $D_2$ and $c_2$ is   $$\ll Z_1 \mathcal{D}_2 \left(1 + \frac{Z_2  }{\mathcal{D}_2}\right) \ll \sqrt{Z_1Z_2} \left(\mathcal{D}_2 + \sqrt{Z_1Z_2}  \right)  \ll  T^{\varepsilon} \sqrt{Z_1Z_2} \frac{\Deltavar}{D} L^{j}  .$$
If $Z_2 \leq Z_1$, we choose $z_2$, $D_1$ and $c_1$ getting the same bound.

Under the additional conditions  $|\Upsilon_2 -  1| + \Upsilon_1 \ll T^{-\Halfb} \ll |U_1V_1|T^{-2} $  we can estimate more efficiently: using first that $|\Upsilon_2 -  1|  \ll T^{-\Halfb}$ and then $|U_1V_1| \gg T^{2- \Halfb}$, we pick  as above $z_2, D_1$ and $c_1$ in
\begin{displaymath}
\begin{split}
& \ll \mathcal{D}_1 \left(1 + \frac{Z_2}{T^{\Halfb}}\right)\left(1 + \frac{Z_1 }{\mathcal{D}_1}\right)  =  \mathcal{D}_1 + Z_1 +  \frac{\mathcal{D}_1 Z_2}{T^{\Halfb}} + \frac{Z_1Z_2}{T^{\Halfb} } \\
& \ll \frac{Z_2}{r_2s_2} + Z_1 +   \frac{Z_1^{1/2} M Z_2 T^{-1+\Quarterbfrac}}{(Ds_0r_2)^{1/2}T^{\Halfb}} + \frac{Z_1Z_2}{T^{\Halfb}} \\
&\ll T^{\varepsilon} (Z_1Z_2)^{1/2} \left(\Bigl(\frac{\Deltavar}{D}\Bigr)^{1/2} L^{j} + \frac{M(\Deltavar/D)^{1/2} L^{j}}{(D s_0r_2)^{1/2} T^{1 + \Quarterbfrac}}  +  \frac{\Deltavar  L^{2j}}{DT^{\Halfbfrac}} 
\right)
 \end{split}
\end{displaymath}
ways which then determine the other variables (up to a divisor function). We also notice that $\Upsilon_1 \ll T^{-\Halfb}$ implies $Z_1\gg r_1s_1 T^{\Halfb}$ and hence $D \ll L^{2j} \Deltavar T^{-\Halfb}.$ 

Similar bounds hold under the dual conditions $|\Upsilon_1 -  1| + \Upsilon_2 \ll T^{-\Halfb}$, $|U_2V_2| \gg T^{2- \Halfb}$.  

Altogether we conclude that $\Sigma_6^{\text{gen}, \ast}(r, s; \mu)$ is at most 
\begin{equation}\label{genast}
\begin{split}
& T^{2+\varepsilon} 
\frac{r_1r_2^{1/2} }{(r_0 s_0s_1)^{1/2} }\Biggl\{  \sum_{D \leq \Deltavar}   T^{-\Eighthbfrac} \frac{\Deltavar}{D} L^{2j} + \sum_{D \ll \frac{L^{2j} \Deltavar}{ T^{\Halfb}}} \left( \changed{\frac{\Deltavar^{\frac{1}{2}}}{D^{\frac{1}{2}}}} L^{j} + \frac{M L^{j}  }{(Ds_0r_2)^{\frac{1}{2}} T^{1+\Quarterbfrac}} \changed{\frac{\Deltavar^{\frac{1}{2}}}{D^{\frac{1}{2}}}}\right)\Biggr\}\\
\ll &  T^{2+\varepsilon}   \frac{r_1r_2^{1/2} }{(r_0 s_0s_1)^{1/2} } \cdot \frac{M^2 L^{3j}}{T^{2+\Eighthbfrac}}  \ll T^{3-\Eighthbfrac+\varepsilon} L^{12} 
  \end{split} 
\end{equation}
by \eqref{D}. \\

\emph{Case III:}  We finally turn to the estimation of $\Sigma_6^{\text{gen}, \text{mix}}(r, s; \mu)$ and assume without loss of generality $c_2 = 0 \not= c_1$ (the other case is analogous). Here we have $\Upsilon_2 =   e^{-2}$, so that the condition $|\Upsilon_2-1| + \Upsilon_1 \ll T^{-\Halfb}$ or $|\Upsilon_1 - 1| + \Upsilon_1 \ll T^{-\Halfb}$ from Lemma \ref{lem6}(d) can never be satisfied. 
Cutting into dyadic ranges, we are left with bounding 
\begin{displaymath}
\begin{split}
 \Sigma^{\text{gen}, \text{mix}}_6&(r, s; \mu; \mathcal{D}, \mathcal{D}_1, \mathcal{D}_2, Z_1) \ll T^{3+\varepsilon}M^2 
\frac{r_1   }{r_0 s_0s_1}\sum_{D\asymp \mathcal{D}}\underset{z_1 = r_1s_1D_2 + c_1D_1}{\sum_{c_1\not= 0} \sum_{\substack{D_1 \asymp \mathcal{D}_1\\ D_2 \asymp \mathcal{D}_2} }\sum_{\substack{  |z_1| \asymp Z_1 \\ z_2 = r_2s_2D_1 }}}\\
 & \sum_{\substack{x_1y_1 = Dz_1\\ x_2y_2 = Dz_2}}\frac{1  }{ DD_1D_2} \Big|\mathcal{K}_{\mu}\left(\frac{-\epsilon_2 M^2  s_1r_1D_2}{DD_1^2 r_2s_0}, \frac{-\epsilon_1 M^2 D_1s_2r_2}{DD_2^2 s_1r_0};  \frac{x_1 M}{DD_1s_0}, \frac{y_1  M}{DD_1r_2}; \frac{x_2 M}{DD_2r_0}, \frac{y_2 M}{DD_2s_1}  \right)\Big|
\end{split}
\end{displaymath}
for $\mathcal{D} \ll \Deltavar$, $\changed{\mathcal{D}_1, \mathcal{D}_2 \ll \Deltavar/\mathcal{D}}$ and  $Z_1 \ll T^{\varepsilon} L^{2j}\Deltavar/\mathcal{D}.$  Picking $z_1$ and $D_2$ determines  (up to a divisor function) $c_1$ and $D_1$. Using \eqref{11cbound2}, we obtain  similarly as before 
\begin{equation}\label{prelim}
\begin{split}
 \Sigma^{\text{gen}, \text{mix}}_6(r, s; \mu; \mathcal{D}, \mathcal{D}_1, \mathcal{D}_2, Z_1)  &  \ll  T^{2-\Eighthbfrac+\varepsilon}  \frac{r_1   }{r_0 s_0s_1  }  \sum_{D \asymp \mathcal{D}} \underset{z_1 = r_1s_1D_2 + c_1D_1}{\sum_{c_1\not= 0} \sum_{\substack{D_1 \asymp \mathcal{D}_1\\ D_2 \asymp \mathcal{D}_2} }\sum_{\substack{  |z_1| \asymp Z_1 \\ z_2 = r_2s_2D_1 }}}  \frac{   \sqrt{r_0s_1s_0r_2}}{\sqrt{ |z_1| r_2s_2  D_1} }\\
& \ll T^{2-\Eighthbfrac+\varepsilon} 
\frac{r_1   }{\sqrt{r_0 s_0s_1 s_2} }  \sum_{D \asymp \mathcal{D}}   \frac{1}{\sqrt{Z_1 \mathcal{D}_1}} \cdot Z_1 \mathcal{D}_2. 
\end{split}
\end{equation}


This bound is acceptable unless
\begin{equation}\label{accept}
  \frac{Z_1M^2}{\mathcal{D} \mathcal{D}_1^2 s_0r_2}  \gg \left(  \frac{r_2s_2\mathcal{D}_1M^2}{\mathcal{D} \mathcal{D}_2^2 r_0s_1} + T^2\right)T^{\varepsilon}.
 \end{equation}
Indeed, if $Z_1 \mathcal{D}_2^2 r_0s_1\ll r_2s_2 \mathcal{D}_1^3 T^{\varepsilon}$, we can continue the preceding calculation with the estimate
\begin{equation}\label{genmix1}
\begin{split}
 \Sigma^{\text{gen}, \text{mix}}_6&(r, s; \mu; \mathcal{D}, \mathcal{D}_1, \mathcal{D}_2, Z_1)  \ll  T^{2-\Eighthbfrac+\varepsilon} \frac{r_1   }{\sqrt{r_0 s_0s_1 s_2} }  \sum_{D \asymp \mathcal{D}}   \mathcal{D}_1 \left(\frac{r_2s_2}{r_0s_1}\right)^{\frac{1}{2}} \ll T^{3-\Eighthbfrac + \varepsilon}L^{6}; 
 \end{split}
 \end{equation}
and if $Z_1M^2 \ll \mathcal{D} \mathcal{D}_1^2 s_0r_2 T^{2+\varepsilon}$, we obtain
 \begin{equation}\label{genmix2}
\begin{split}
 \Sigma^{\text{gen}, \text{mix}}_6&(r, s; \mu; \mathcal{D}, \mathcal{D}_1, \mathcal{D}_2, Z_1) \\
 &  \ll  T^{2-\Eighthbfrac+\varepsilon}  
 \frac{r_1   }{\sqrt{r_0 s_0s_1 s_2} }  \sum_{D \asymp \mathcal{D}}   \frac{\mathcal{D}_2(\mathcal{D} \mathcal{D}_1s_0r_2)^{1/2}T}{ M} \ll T^{3-\Eighthbfrac + \varepsilon} L^{\frac{15}{2}}
  \end{split}
 \end{equation}
by \eqref{D}. From now on we assume   \eqref{accept}. In this case we have
$$|U_1V_1| \gg (|U_2V_2| + T^2)T^{\varepsilon}$$
 with the notation from \eqref{variables}, so that Lemma \ref{lem6}(e) gives us the additional information $$\Upsilon_1 = 1 + O\left(\frac{T^{\varepsilon}}{|U_1V_1|^{1/4}} + T^{\varepsilon} \frac{\abs{U_2 V_2}^{1/2} + T}{|U_1V_1|^{1/2}}\right)$$
(up to a negligible error). Having picked $z_1$, we conclude that
$$D_2 = \frac{e^2 z_1}{r_1s_1}\left( 1 + O\left(\frac{T^{\varepsilon}}{|U_1V_1|^{1/4}} + T^{\varepsilon} \frac{\abs{U_2 V_2}^{1/2} + T}{|U_1V_1|^{1/2}}\right)\right)$$
and in particular 
$Z_1/(r_1s_1) \asymp \mathcal{D}_2$. Hence  the number of pairs $(z_1, D_2)$ is 
\begin{displaymath}
\begin{split}
& \ll Z_1 \left(1 + \mathcal{D}_2 T^{\varepsilon}\left(\frac{\mathcal{D}^{1/4} \mathcal{D}_1^{1/2} (s_0r_2)^{1/4}}{Z_1^{1/4} M^{1/2}} + \frac{T\mathcal{D}^{1/2} \mathcal{D}_1 (s_0r_2)^{1/2}}{Z_1^{1/2} M}\right)\right).\\
\end{split}
\end{displaymath}
Replacing the last factor $Z_1\mathcal{D}_2$ on the right hand side of  \eqref{prelim} with this quantity, we obtain under the present assumption \eqref{accept} that
\begin{equation}\label{genmix3}
\begin{split}
&   \Sigma^{\text{gen}, \text{mix}}_6 (r, s; \mu; \mathcal{D}, \mathcal{D}_1, \mathcal{D}_2, Z_1)\\
   & \ll \frac{ T^{2-\Eighthbfrac+\varepsilon}   
r_1   }{\sqrt{r_0 s_0s_1 s_2} }  \sum_{D \asymp \mathcal{D}}  \left(Z_1^{\frac{1}{2}} + \frac{Z_1^{\frac{1}{4}}\mathcal{D}_2(\mathcal{D} s_0r_2)^{\frac{1}{4}}}{M^{1/2}} + \frac{\mathcal{D}_2 T (\mathcal{D} \mathcal{D}_1 s_0r_2)^{\frac{1}{2}}}{M}\right) \ll T^{3-\Eighthbfrac + \varepsilon} L^{\frac{15}{2}}.
   \end{split}
\end{equation}
By \eqref{genmix1}, \eqref{genmix2} and \eqref{genmix3} we obtain
\begin{equation}\label{genmixfin}
   \Sigma^{\text{gen}, \text{mix}}_6(\mu) \ll T^{3 - \Eighthbfrac + \varepsilon} L^{\frac{15}{2}}
\end{equation} 
in all cases. \\
 
Combining \eqref{gen0}, \eqref{genast} and \eqref{genmixfin}, we obtain
\begin{equation*} 
\Sigma_6^{\text{gen}}(\mu) \ll T^{3 - \Eighthbfrac + \varepsilon} L^{12}.
\end{equation*}

Together with \eqref{central}  and \eqref{mix} we obtain finally
\begin{equation}\label{6finalbound}
  \Sigma_6(\mu) \ll T^{\varepsilon}\left(T^{3 -\lambda}+  T^{3 - \Eighthbfrac  + 12\lambda}  + T^{\frac{9}{4} + \frac{\eta}{2} + 31\lambda} + T^{\frac{5}{2} + 22\lambda}\right). 
\end{equation}

\section{Proof of Theorem \ref{thm1}}

Collecting the bounds \eqref{Delta}, \eqref{w4final} and \eqref{6finalbound}, we see that we can bound $\mathcal{L}_M(\pi_0)$ in \eqref{basic} by
$$\mathcal{L}_M(\pi_0)^2 \ll T^{\varepsilon} \left(T^{3 - \Eighthbfrac  + 12\lambda}  + T^{\frac{9}{4} + \frac{\eta}{2} + 31\lambda} + T^{\frac{5}{2} + 22\lambda} + T^{\frac{3}{2} + 3\lambda + 2\eta} + T^{3 - \lambda}\right)$$
and we recall the  trivial bound 
$$\mathcal{L}_M(\pi_0)^2 \ll T^{3 - 2\eta+\varepsilon},$$
see \eqref{etasmall}. Recalling \eqref{basic1}, we choose the trivial bound if $\eta \leq 1/100$, otherwise we choose $\lambda = \deltaFinal$  to complete the proof of Theorem \ref{thm1}. \\

The rest of the paper is concerned with the proofs of the bounds from Section \ref{key}.

\section{Proof of Lemma \ref{lem3}}
We start by inserting the integral representation of Lemma \ref{stade0} into the definition of $\Phi_{w_4}$ in \eqref{defPhi}. In the $\mu$-plane we introduce the new variables $\rho = (\rho_1, \rho_2)$, where $\rho_1 = \Im(\mu_1 + \mu_2)$ and $\rho_2 = \Im(\mu_1 - \mu_2)$. We start with the integral involving the   $\tilde{K}$-function:
$$\int_{\Bbb{R}^2}   h(\mu) \int_0^{\infty} 
\exp\left(i \rho_1 \log \frac{\pi^3 |y|}{u^{3/2}}\right) 
\tilde{K}_{i\rho_2}(2\sqrt{u}) \exp\left(\pm \frac{2i \pi^3y}{u}\right) \frac{du}{u} \text{spec}(\mu) d\rho.$$
We recall that the $u$-integral is not absolutely convergent at $0$, but as in the proof of Lemma \ref{stade0} this causes no substantial difficulties (we can temporarily integrate by parts in the region $0 < u < 1$, for instance).  
We can  at the cost of a negligible error   replace $h$ by a real-analytic function that still satisfies \eqref{partial} and \eqref{vol} and is compactly supported in $\min_{w \in \mathcal{W}} |\mu - w(\mu_0)| \leq T^{2\varepsilon}$, say.  In particular, by our assumptions on $\mu_0$ we can assume from now on that  $\mu$ satisfies \eqref{munu}.    
It is easy to see that the $\rho_1$-integral is negligible  unless 
\begin{equation}\label{uw4}
  u \asymp  |y|^{2/3},  
\end{equation}
and (at the cost of a negligible error) we restrict the integral to  this  interval. 

Assume first $ |y| \leq T^{3-\varepsilon}.$ 
Then we conclude 
\begin{equation}\label{uw4a}
  u^{1/2} \ll T^{1 -\varepsilon/3},
  \end{equation} so that   we can insert the uniform asymptotic formula \eqref{unifasymp}. 
The error term saves arbitrarily many powers of $T$ (choosing $M$ large enough) and is therefore admissible. 
Moreover, with $H(\rho_2) = \omega(4\sqrt{u}, \rho_2)$ we have
$$|H'(\rho_2)| = 
 \text{arccosh}\left(\frac{|\rho_2|}{4\sqrt{u}}\right) \gg_{\varepsilon} \log T \gg 1, \quad H^{(j)}(\rho_2)  \ll T^{1-j} \quad (j \geq 2)$$
 by \eqref{uw4a} (recall that $|\rho_2| \asymp T$).  
Integrating by parts sufficiently often  by means of  \eqref{partial} and Lemma \ref{integrationbyparts} with $R = 1$, $Y = Q = T$, $U = T^{\varepsilon}$, we see that the $\rho_2$-integral is negligible as desired. This proves \eqref{w41}. 

The integral involving the $J^{\pm}$-function can be treated in the same way using the analogous formula \eqref{unifasymp1}. 

On the other hand, if $|y| \geq T^{3-\varepsilon}$, we differentiate $j$ times with respect to $y$ under the integral sign.  Keeping in mind that $u$ is restricted to the range \eqref{uw4}, each such differentiation produces a factor $T|y|^{-1} + u^{-1} \asymp (T +  |y|^{1/3})|y|^{-1}$, and a trivial estimate using \eqref{defPhi}, \eqref{vol}, Lemma \ref{stade0} and \eqref{besselbound} (with $j = 0$) 
completes the proof of \eqref{w42}.  \hfill $\square$

\section{Proof of Lemma \ref{lem4}}

   The strategy of the proof of is similar to the preceding one, but the details are more involved. We will have to play off some $\varepsilon$'s against others, therefore need to be very careful with the value of $\varepsilon$ and again do not use $\varepsilon$-convention in this proof. 

We start with the discussion of the $(+ +)$ case where $y_1, y_2 > 0$.  
Without loss of generality we can assume $y_1 \geq y_2$.   As in the previous proof we  replace   $h$, at the cost of a negligible error,  by a real analytic function $\tilde{h}$  that  satisfies \eqref{partial} and is compactly supported in $\min_{w \in \mathcal{W}} |\mu - w(\mu_0)| \leq T^{2\varepsilon}$. In particular, we can assume from now on that $\mu$ satisfies \eqref{munu}. It follows directly from Definition \ref{MB}  that   $K^{++}_{w_6} $ is symmetric in $\mu$, hence we can and will  assume without loss of generality that $\tilde{h}$ is supported only in the positive Weyl chamber $\Im \nu_1, \Im \nu_2 \geq 0$, so that $\Im \nu_3 \leq 0$.  There we have
\begin{equation*}
 \frac{ \cos\left(\frac{3}{2}\pi \nu_1\right) \cos\left(\frac{3}{2}\pi \nu_2\right)}{ \cos\left(\frac{3}{2}\pi \nu_3\right)} 
= \frac{1}{2} + O(e^{-3\pi \min(|\nu_2|, |\nu_1|)}) = \frac{1}{2} + O_B(T^{-B})
\end{equation*}
for any $B \geq 0$. 

Now we insert the integral representation \eqref{lem++} and \eqref{BesselK}. We consider first  the $\mu_2$-integral 
\begin{equation*}
\int_{(0)} \tilde{h}(\mu)  \left(\frac{y_1}{y_2}\right)^{\frac{1}{2}\mu_2} u^{3\mu_2}  \text{{\rm spec}}(\mu) d\mu_2.  
\end{equation*}
 Integrating by parts in combination with \eqref{partial}, we can save arbitrarily many powers of $T$ unless
\begin{equation}\label{u}
u = \frac{y_2^{1/6}}{y_1^{1/6}}\left(1 + O(T^{-\varepsilon/2})\right).
\end{equation}
It follows in particular that  
\begin{equation}\label{sizes}
y_1^{1/2} \sqrt{1 + u^2} \asymp y_1^{1/2}, \quad y_2^{1/2}\sqrt{1 + u^{-2}} \asymp  y_2^{1/3} y_1^{1/6}.
\end{equation}
At this point we see the significance of the somewhat technical looking cut-off point $\Upsilon$: if at least one of the  arguments in the Bessel functions in \eqref{BesselK} is significantly smaller than the index $3\nu_3 \asymp T$ of the Bessel function, then we claim that the integral is negligible.   

If  $\Upsilon \geq T^{1-\varepsilon}$,  we differentiate $j_1$ times with respect to $y_1$ and $j_2$ times with respect to $y_2$ under the integral sign and estimate trivially, using \eqref{defPhi}, \eqref{vol}, \eqref{besselbound} and \eqref{sizes}.  This proves \eqref{J2}. 

From now on we assume 
\begin{equation}\label{T}
\Upsilon \leq T^{1-\varepsilon}
\end{equation} 
and aim at proving \eqref{J1}. By \eqref{sizes} and \eqref{T}, we 
can insert into the second Bessel function in \eqref{BesselK} the uniform asymptotic formula \eqref{unifasymp} (with $t = 3 \Im  \nu_3$, $x = 4\pi y_2^{1/2} \sqrt{1+u^{-2}} \ll |t|^{1-\varepsilon}$ and $M$ sufficiently large so that the error term is negligible). 
The treatment of the other Bessel function in \eqref{BesselK} depends on the  relative size of $y_1$ and $y_2$.\\

\emph{Case I:} Let us first assume that $y_1$ and $y_2$ are ``close'' in the sense that $y_2 \geq y_1T^{-\varepsilon}$. Then $y_1^{1/2} \leq y_2^{1/2} T^{\varepsilon/2} \leq y_2^{1/3} y_1^{1/6} T^{\varepsilon/2} \leq T^{1-\varepsilon/2}$ by \eqref{T}, hence by \eqref{sizes} we can insert  \eqref{unifasymp} also for the other Bessel function, again with a negligible error term. Now we consider the $\nu_3$-integral, where we write momentarily $\rho := |3\Im \nu_3|$ for notational simplicity: 
\begin{equation}\label{pm1}
\begin{split}
\int_{-\infty}^{\infty} & \exp\left(\pm i \omega(4\pi y_2^{1/2} \sqrt{1+u^{-2}}, \rho )\right)  \exp\left( \pm i \omega(4\pi y_1^{1/2} \sqrt{1+u^{2}}, \rho)\right)  \\
& \times   \tilde{h}(\mu)  f^{\pm}_M\left(4\pi y_2^{1/2} \sqrt{1+u^{-2}}, \rho\right)  f^{\pm}_M\left(4\pi y_1^{1/2} \sqrt{1+u^{2}}, \rho \right)   \text{spec}(\mu) d\rho,
\end{split}
\end{equation}
where $f_M^+ = f_M$, $f_M^- = \bar{f}_M$ as in \eqref{unifasymp}. 
With $H(\rho) = \pm \omega(4\pi y_2^{1/2} \sqrt{1+u^{-2}}, \rho ) \pm \omega(4\pi y_1^{1/2} \sqrt{1+u^{2}}, \rho)$ we have 
\begin{equation}\label{pm2}
H'(\rho) = \pm  \, \text{arccosh}\left(\frac{\rho}{4\pi y_2^{1/2} \sqrt{1+u^{-2}}}\right) \pm \, \text{arccosh}\left(\frac{\rho}{4\pi y_1^{1/2} \sqrt{1+u^{2}}}\right)
\end{equation}
and $H^{(j)}(\rho) \ll T^{1-j}$ for $j \geq 2$ whenever  $\rho \asymp T$ and \eqref{sizes} and \eqref{T} hold.  
Moreover, in this case the first term in \eqref{pm2} is $\gg_{\varepsilon} \log T$, hence recalling \eqref{partial} and the asymptotic formula $\text{arccosh}(x) = \log(2x) + O(x^{-2})$, we can save arbitrarily many powers of $T$ by integrating by parts unless  the two $\pm$ signs are different and 
$$\frac{y_2^{1/2}\sqrt{1 + u^{-2}}}{y_1^{1/2}\sqrt{1 + u^{2}}} = 1 + O(T^{-\varepsilon/2}),$$
or in other words
\begin{equation}\label{pm3}
u = \frac{y_2^{1/2}}{y_1^{1/2}}\left( 1 + O(T^{-\varepsilon/2})\right).
\end{equation}
(Here we used again Lemma \ref{integrationbyparts} with $U = T^{\varepsilon}$, $Y= Q = T$, $R = T^{\varepsilon/2}$.) 
Combining this with \eqref{u} 
we see that 
we can localize $u$ at $|u -1| \ll T^{-\varepsilon/2}$ at the cost of a negligible error.  We insert a corresponding smooth cut-off function $\psi(u)$ with $\psi^{(j)}(u) \ll_j T^{\frac{\varepsilon}{2} j}$ for $j \in \Bbb{N}_0$ into the integral and consider
\begin{equation}\label{715}
\begin{split}
 \int_0^{\infty} &\exp \left(\epsilon  i \omega(4\pi y_2^{1/2} \sqrt{1+u^{-2}},3\rho )- \epsilon  i \omega(4\pi y_1^{1/2} \sqrt{1+u^{2}}, 3\rho )\right)  u^{3\mu_2} \\
&\times \psi(u)  f^{\epsilon}_M\left(4\pi y_2^{1/2} \sqrt{1+u^{-2}}, 3\rho\right)  f^{-\epsilon}_M\left(4\pi y_1^{1/2} \sqrt{1+u^{2}}, 3\rho\right)  \frac{du}{u}
\end{split}
\end{equation}
for $\epsilon  \in \{\pm \}$ with $f_M^{\pm}$ as in \eqref{pm1}. 
With $$H(u) = \epsilon \omega(4\pi y_2^{1/2} \sqrt{1+u^{-2}},3\rho )- \epsilon  \omega(4\pi y_1^{1/2} \sqrt{1+u^{2}}, 3\rho ) +  3\Im \mu_2 \log u$$ and recalling $\rho = 3|\Im \nu_3|$ we have
\begin{displaymath}
\begin{split}
H'(u) &= \frac{3\Im \mu_2}{u} +\epsilon  \frac{\sqrt{( 3\Im \nu_3 )^2 u^2 - 16\pi^2y_2 (1 + u^2)}}{u^2(1+u^2)}+ \epsilon  \frac{u \sqrt{( 3\Im \nu_3 )^2  - 16\pi^2y_1 (1 + u^2)}}{ (1+u^2)} \\
& = \left(3\Im\mu_2 + \epsilon \Bigl |\frac{3}{2}\Im \nu_3\Bigr| + \epsilon \Bigl|\frac{3}{2}\Im \nu_3\Bigr|\right) \left(1 + O(T^{-\varepsilon/2})\right) = \left\{ \begin{array}{ll}  6\Im \nu_2, & \epsilon = 1\\ -6\Im \nu_1, & \epsilon = -1\end{array} \right\} \left(1 + O(T^{-\varepsilon/2})\right) \asymp T
\end{split}
\end{displaymath}
  and $H^{(j)}(u) \ll T^{1-j}$ for $j \geq 2$. 
Hence by partial integration (using Lemma \ref{integrationbyparts} with $R = Y = T$, $U = T^{-\varepsilon/2}$, $Q = 1$), the $u$-integral is negligible.\\

\emph{Case II:} On the other hand, if $y_2 \leq y_1T^{-\varepsilon}$,  we substitute for the first Bessel function in \eqref{BesselK} the integral representation \eqref{intrep1}
with $t = 3 \Im  \nu_3$, $x = 4\pi y_1^{1/2} \sqrt{1+u^{2}}$. 
We consider first the $\nu_3$-integral (again with $\rho := |3\Im \nu_3|$)
$$\int_{(0)} \exp\left(\pm i \omega(4\pi y_2^{1/2} \sqrt{1+u^{-2}}, \rho )\right)  \exp(i \rho v) \tilde{h}(\mu)  f^{\pm}_M\left(4\pi y_2^{1/2} \sqrt{1+u^{-2}}, \rho\right) d\rho.$$
With $H(\rho) = \pm \omega(4\pi y_2^{1/2} \sqrt{1+u^{-2}}, \rho ) + \rho v$  we have
$$ H'(\rho) = v  \pm  \, \text{arccosh}\left(\frac{\rho}{4\pi y_2^{1/2} \sqrt{1+u^{-2}}}\right), \quad H^{(j)}(\rho) \ll T^{1-j} \quad (j \geq 2)$$
whenever  $\rho \asymp T$ and \eqref{sizes} and \eqref{T} hold.  
Integration by parts in connection with \eqref{partial} implies as above that the $\nu_3$-integral is negligible unless
\begin{equation}\label{condv}
|v| = \, \text{arccosh}\left(\frac{\rho}{4\pi y_2^{1/2} \sqrt{1+u^{-2}}}\right) + O(T^{-\varepsilon/2}), 
\end{equation}
which implies 
$$\cosh v = \frac{\rho}{4\pi y_2^{1/2} \sqrt{1+u^{-2}}}(1 + O(T^{-\varepsilon/2})) \asymp \frac{T}{y_2^{1/3} y_1^{1/6}}$$
by \eqref{u}. As before, we can remember \eqref{condv} by inserting a smooth cut-off function $\tilde{\psi}(v)$ with $\tilde{\psi}^{(j)} \ll_j T^{\frac{\varepsilon}{2} j}$ for $j \in \Bbb{N}_0$ at the cost of a negligible error. Then the $v$-integral becomes
$$\int_{\Bbb{R}} \tilde{\psi}(v) \cos\left(4\pi y_1^{1/2} \sqrt{1+u^2}  \sinh v\right) \exp(i \rho v)  dv. $$
With $H(v) = \rho v \pm 4\pi y_1^{1/2} \sqrt{1+u^2}  \sinh v$ we have
$$H'(v) = \rho  \pm 4\pi y_1^{1/2} \sqrt{1+u^{2}} \cosh v \asymp  O(T) + T \left(\frac{y_1}{y_2}\right)^{1/3} \asymp T\left(\frac{y_1}{y_2}\right)^{1/3}$$
and $H^{(j)} (v) \ll T(y_1/y_2)^{1/3}$ for $j \geq 2$.  Applying Lemma \ref{integrationbyparts} with $Y = R = T(y_1/y_2)^{1/3}$, $Q = 1$, $R = T^{-\varepsilon/2}$, we see that the $v$-integral is negligible.\ This completes the proof of \eqref{J1} in the (++) case. \\

The $(--)$ case is  similar. By symmetry we have
$$\Phi_{w_6}(y) =  \frac{1}{3}\int_{\Re \mu = 0} h(\mu) \sum_{w \in \{I, w_4, w_5\}} K^{--}_{w_6}(y ;  w(\mu) )\, \text{spec}(\mu) d \mu$$
for $y_1, y_2 < 0$. 
Here we use \eqref{lem--} and \eqref{BesselJ1} in place of \eqref{lem++} and \eqref{BesselK}\changed{, and} argue as before.  
The analogue of \eqref{unifasymp} is \eqref{unifasymp1}. 
This case is a little simpler because  it is not necessary to distinguish between $|y_2| \geq |y_1|T^{-\varepsilon}$ and $|y_2| \leq |y_1|T^{-\varepsilon}$, the  method of Case I works regardless of the relative size of $y_1$ and $y_2$, since \eqref{unifasymp1} has no transitional range.  This completes the proof of the lemma in the $(--)$ case. \\

Finally we treat the $(+-)$ case (and the $(-+)$ follows by \eqref{lem-+}). Here we use the integral representation \eqref{lem+-} together with \eqref{BesselJ2}--\eqref{BesselJ4}. The treatment of the term involving $\mathcal{J}_3$ is identical to the preceding argument. The $u$-integrals in $\mathcal{J}_2$ and $\mathcal{J}_4$ are slightly different and require a  variation of the argument. The analogue of \eqref{u} is again
\begin{equation}\label{unew}
u = \frac{|y_2|^{1/6}}{|y_1|^{1/6}}(1 + O(T^{-\varepsilon/2})),
\end{equation}
but we have only a weaker version of \eqref{sizes}. In the case of $\mathcal{J}_2$, we have
$$|y_1|^{1/2} \sqrt{u^2 - 1} \ll |y_1|^{1/3} |y_2|^{1/6}, \quad |y_2|^{1/2} \sqrt{1 - u^2} \ll |y_2|^{1/2},$$
and in the case of $\mathcal{J}_4$, we have 
$$|y_1|^{1/2} \sqrt{1-u^2} \ll |y_1|^{1/2}, \quad  |y_2|^{1/2} \sqrt{  u^{-2}-1} \ll |y_2|^{1/3} |y_1|^{1/6}.$$
This is still sufficient   to prove \eqref{J2}.   From now on we assume \eqref{T}. Let us first assume that  $|y_1|$, $|y_2|$ are not close, i.e.\ 
\begin{equation}\label{notclose}
  \min(|y_1|, |y_2|) \leq \max(|y_1|, |y_2|)T^{-\varepsilon}.
\end{equation}  
   Then for $\mathcal{J}_4$ the argument of Case II goes through with minor notational changes (note that $u\leq 1$ implies automatically $|y_1| \geq |y_2|$ in this case). For $\mathcal{J}_2$ we can copy the argument of Case I above. Note that the assumption $y_2 \geq y_1 T^{-\varepsilon}$ was only needed to insert the uniform asymptotic asymptotic expansion, which in the present case of the Bessel $J$-function can   be done even under the assumption \eqref{notclose};  moreover, \eqref{notclose} in connection with \eqref{unew} implies that $u$ is supported away from 1. 

It remains to treat the case when $|y_1|$ and $|y_2|$ are close, i.e.\   $\min(|y_1|, |y_2|) \geq \max(|y_1|, |y_2|)T^{-\varepsilon}$.  As  in \eqref{pm3} we conclude
$$u = \frac{|y_2|^{1/2}}{|y_1|^{1/6}} (1 + O(T^{-\varepsilon/2}))$$
 so that we can localize $u$ at $|u-1| \ll T^{-\varepsilon/2}$, and we have automatically 
 \begin{equation}\label{y1y2}
    |y_1| = |y_2| (1+O(T^{-\varepsilon/2}))
 \end{equation}   
     (otherwise the contribution is negligible).  We are now facing the small technical problem that the integrals $\mathcal{J}_2$ and $\mathcal{J}_4$ have a sharp cut-off at $u=1$ which prevents partial integration with respect to $u$. Therefore we first extract smoothly the region $|u- 1| \ll T^{-1 + \varepsilon}$ and insert a smooth weight function $\psi(u)$ with support in $T^{-\varepsilon/2} \gg |u-1| \gg T^{-1+\varepsilon}$ and $\psi^{(j)} \ll_j T^{(1 - \varepsilon)j}$. We can then apply the argument following  \eqref{715} except that we apply Lemma \ref{integrationbyparts} with $U = T^{-1+\varepsilon}$ instead of $T^{-\varepsilon/2}$. 
     
For the remaining region $|u-1| \ll T^{-1+\varepsilon}$ we need to glue together the two integrals $\mathcal{J}_2$ and $\mathcal{J}_4$ and consider the integral 
\begin{equation*}
  \int_0^\infty  \tilde{\psi}(u) \left(\mathcal{A}_{3\nu_3 }(u) - \mathcal{A}_{-3\nu_3 }(u) \right) \left(\mathcal{B}_{3\nu_3}(u) - \mathcal{B}_{-3\nu_3}(u) \right)
  u^{3\mu_2} \frac{du}{u},
\end{equation*}     
 where $\tilde{\psi}(u)$ has  support in $|u-1| \ll T^{-1+\varepsilon}$ and satisfies $\tilde{\psi}^{(j)} \ll_j  T^{(1-\varepsilon)j}$, and
 \begin{equation}\label{B1}
   \mathcal{A}_{\alpha}(u) = \sin(\pi \alpha/2)^{-1}  \begin{cases} J_{\alpha}(4\pi |y_1|^{1/2} \sqrt{u^2 - 1}), & u > 1,\\ I_{\alpha}(4\pi |y_1|^{1/2} \sqrt{1-u^2}), & u < 1,\end{cases}
 \end{equation}
 and 
  \begin{equation}\label{B2}
   \mathcal{B}_{\alpha}(u) = \sin(\pi \alpha/2)^{-1} \begin{cases} J_{\alpha}(4\pi |y_2|^{1/2} \sqrt{1-u^{-2} }), & u > 1,\\ I_{\alpha}(4\pi |y_2|^{1/2} \sqrt{ u^{-2}-1}), & u < 1.\end{cases}
 \end{equation}
 Here we used \eqref{I}. Notice that the spectral measure has a double zero  at $\nu_j = 0$ which cancels the poles of the two $\sin(\pi \alpha/2)^{-1}$ at $\alpha = 3\nu_3 = 0$. 
By \eqref{T}, \eqref{y1y2} and the support of $\tilde{\psi}$, the arguments of the Bessel functions in \eqref{B1} and \eqref{B2} are $\ll T^{(1-\varepsilon)/2}$. This is important, because now the power series expansions \eqref{I1} and \eqref{powerJ} are decreasing term by term, and truncating at $k < M$ gives an error of size $$O\left(T^{2 \cdot \frac{1-\varepsilon}{2} \cdot M} T^{-M-1/2}\right) = O(T^{-\frac{1}{2} - \varepsilon M})$$
 for $|\nu_3| \asymp T$ and  $|u-1| \ll T^{-1+\varepsilon}$. 
 This is negligible for sufficiently large $M$. 
 
The remaining terms of the power series expansions produce integrals of the type    
 \begin{displaymath}
 \begin{split}
  \int_0^{\infty} &  \left\{\begin{array}{ll} (-1)^{k_1+k_2}, & u > 1\\ 1, & u < 1 
\end{array}\right\}  \frac{(2\pi|y_1|^{1/2}\sqrt{|1 - u^2|})^{ \epsilon_1  3\nu_3 + 2k_1}(2\pi|y_2|^{1/2}\sqrt{ |u^{-2}-1|})^{\epsilon_2  3\nu_3 + 2k_2} }{\sin(\frac{3}{2} \pi \nu_3)^2 k_1! k_2! \Gamma(k_1 + 1 + \epsilon_1  3\nu_3)\Gamma(k_2 + 1 + \epsilon_2  3\nu_3) }   \tilde{\psi}(u) u^{3\mu_2} \frac{du}{u} 
\end{split}
\end{displaymath}
 for $\epsilon_1, \epsilon_2 \in \{\pm 1\}$.  
 
 If $\epsilon_1 = \epsilon_2$, 
 we keep $u$ fixed and  pull the $\nu_3$-integral inside. We replace the modified function $\tilde{h}(\mu)$ by the original holomorphic function $h(\mu)$ at the cost of a negligible error. Now we  
 can  shift the $\nu_3$ contour to the far right if $\epsilon_1 = 1$ and to the far left if $\epsilon_1 = -1$, thereby saving arbitrarily many powers of $T$. The decay properties of $h$ ensure that the poles  of $1/\sin(\frac{3}{2} \pi \nu_3)$  contribute negligibly.  
 
 If $\epsilon_1 = - \epsilon_2$, then we are left with a $u$-integral of the form
 \begin{displaymath}
 \begin{split}
 & \int_0^{\infty} \tilde{\psi}(u) \left\{\begin{array}{ll} (-1)^{k_1+k_2}, & u > 1\\ 1, & u < 1 
\end{array}\right\}   |1 - u^2|^{k_1}  |u^{-2}-1|^{ k_2} u^{3\mu_2 + \epsilon_1 3\nu_3}  \frac{du}{u} \\
= &  \int_0^{\infty} \tilde{\psi}(u) \frac{(u^2 - 1)^{k_1+k_2}}{u^{2k_2}}  \left\{\begin{array}{ll} u^{-6\nu_1}, & \epsilon_1  = 1\\ u^{6\nu_2}, & \epsilon_1 = -1 
\end{array}\right\}  \frac{du}{u}. 
\end{split}
\end{displaymath}
Notice that the integrand on the right hand side, unlike its appearance on the left hand side, is smooth. 
Integrating by parts one last time with $H(u) = -6 \Im \nu_1 \log u$ or $H(u) = 6\Im \nu_2 \log u$, $R =Q = Y = T$, $U = T^{1-\varepsilon}$, we see that this integral is negligible. This completes the proof of \eqref{J1} in the $(+-)$ case.   \hfill $\square$

\section{Proof of Lemma \ref{lem5}}

 We start by inserting the Mellin-Barnes representation from Definition \ref{MB} into the right hand side of \eqref{deftildeK}. We choose the integration line $\Re s = \varepsilon$. 
 We then need to analyze
  \begin{equation}\label{start}
 \begin{split}
 \int_{\Bbb{R}^2} \int_{(\varepsilon)} & \tilde{G}^{\pm}(s, \mu) | \xi\eta\Xi|^{-s}  \frac{ds }{2\pi i } \, e(\xi U+\eta V )     W(\xi) W(\eta)  d\xi \, d\eta. 
 \end{split}
 \end{equation}
If $V = 0$, we integrate by parts with respect to $\eta$ to see that we can restrict the $s$-integral to $|\Im s| \leq T^{\varepsilon/2}$, the remaining part being negligible. But then the $\xi$-integral is negligible unless $U \leq T^{\varepsilon}$. This proves (a). 

If $U = V = 0$, we obtain
 \begin{displaymath}
 \begin{split}
  \int_{(\varepsilon)} & \tilde{G}^{\pm}(s, \mu) | \Xi|^{-s} \widehat{W}(1-s)^2 \frac{ds }{2\pi i}, 
  \end{split}
 \end{displaymath}
 where the Mellin transform $\widehat{ W}$ is entire and rapidly decaying. 
 The upper bound $\tilde{\mathcal{K}}_{\mu}(\Xi, 0, 0) \ll T^{-3/2 + \varepsilon}$ follows from Stirling's formula. If $|\Xi| \geq T^{3+\varepsilon}$, we can shift the contour to the right to see that the integral is negligible (absolute convergence is ensured by the rapid decay of $\widehat{W}$). 
This completes the proof of  (b). 

Finally we prove (c). Here we assume $|U|, |V| \geq T^{\varepsilon}$. We integrate over $\xi$ and $\eta$ using Lemma \ref{statphas} and see that $\tilde{\mathcal{K}}_{\mu}(\Xi; U, V)$ is negligible unless 
$t \asymp U \asymp V$ 
in which case \eqref{start} equals (up to a negligible error)
$$\frac{|\Xi|^{-\varepsilon}}{|UV|^{1/2}}  \int_{\Bbb{R}}  \tilde{G}^{\pm}(\epsilon + i t, \mu) \left(\frac{|\Xi t^2|}{|UV|}\right)^{-it} F\left(\frac{t}{U}\right) dt$$
for some fixed,  smooth, compactly supported function $F$. Writing $r_j = \Im \mu_j$ for $j = 1, 2, 3$ and recalling \eqref{munu}, a trivial estimate using Stirling's formula shows the bound 
$$\tilde{\mathcal{K}}_{\mu}(\Xi, U, V) \ll \frac{(|UV| + T)^{\varepsilon}}{|UV|^{1/2}} \int_{t \asymp U}  \prod_{j=1}^3 (1 + |t - r_j|)^{-1/2} dt 
\ll \frac{1}{(|UV|^{1/2} + T)^{3/2 - \varepsilon}}.$$

\section{Proof of Lemma \ref{lem6}}


\subsection{The simple parts} As in the previous section we start by inserting the Mellin-Barnes representation from Definition \ref{MB} into \eqref{defK}. In all cases we choose the contour $$[\varepsilon  , \varepsilon + i T^{B}] \cup [\varepsilon + iT^{B}, -1/10 + iT^B] \cup [-1/10 +iT^B, -1/10 + i \infty)$$ together with its reflection on the real axis, for some sufficiently large constant $B$ and some sufficiently small $\varepsilon> 0$. The third portion is negligible by Stirling's formula, the second portion is negligible by Lemma \ref{statphas} and the fact that $U_j, V_j \ll T^{O(1)}$.  
Hence we then need to analyze
  \begin{displaymath}
 \begin{split}
 \int_{\Bbb{R}^4} \int_{|t_j| \leq T^B} & G(\varepsilon + i t, \mu) S^{\epsilon_1, \epsilon_2}(\varepsilon+ it, \mu)|4\pi^2 \xi_1\eta_1\Xi_1|^{-\varepsilon - it_1}|4\pi^2 \xi_2\eta_2\Xi_2|^{-\varepsilon - it_2} \frac{dt_1\, dt_2}{(2\pi)^2}\\
 & \times e(\xi_1 U_1+\eta_1 V_1+\xi_2 U_2+\eta_2 V_2)     W(\xi_1) W(\eta_1)\overline{W}(\xi_2) \overline{W}(\eta_2)  d\xi_1\, d\xi_2\, d\eta_1\, d\eta_2,
 \end{split}
 \end{displaymath}
where $(\epsilon_1, \epsilon_2) = (\text{sgn}(\Xi_1), \text{sgn}(\Xi_2))$.  If $U_1 = 0$, then after integrating by parts in the $\xi_1$-integral we see easily that we can restrict the $t_1$-integral to $|t_1| \leq T^{\varepsilon/2}$, the remaining part being negligible.  But then the $\eta_1$-integral is negligible unless $V_1 \leq T^{\varepsilon}$. This proves part (a).\\
 
To prove part (b), we recall $U_2 = V_2 = 0$ and $|U_1|,  |V_1| \geq T^{\varepsilon}$. We integrate over $\xi_1$ and $\eta_1$ using Lemma \ref{statphas} and see that the $t_1$-integral is negligible unless 
$t_1 \asymp U_1 \asymp V_1$. We shift the $s_2$-contour to $\Re s_2 = 1-2\varepsilon$ (note that because of a possible pole at $s_1 + s_2 = 1$ of $S^{\epsilon_1, \epsilon_2}(s, \mu)$ we cannot shift much further to the right) and truncate the contour at $|\Im s_2| \leq T^{\varepsilon}$ by the rapid decay of $\widehat{W}(s_2)$. Now we estimate $\mathcal{K}_{\mu}(\Xi_1, \Xi_2; U_1, V_1; 0, 0)$ trivially  by 
\begin{displaymath}
\begin{split}
 & \frac{1}{|U_1V_1|^{1/2}}\int_{|t_2| \leq T^{\varepsilon}} \int_{t_1 \asymp U_1} \frac{|G\left((\varepsilon + i t_1, 1 - 2\varepsilon + it_2), \mu\right) S^{\epsilon_1, \epsilon_2}\left((\varepsilon+ it, 1 - 2\varepsilon + it_2), \mu\right)| }{ |\Xi_2|^{1-2\varepsilon}} dt_1\, dt_2\\
   \ll &  \frac{1}{|U_1V_1|^{1/2}} \cdot T^{\varepsilon} |U_1|   \cdot \frac{T^{3(\frac{1}{2} - 2\varepsilon)}}{(T+|U_1|)^{3(\frac{1}{2}-\varepsilon)} |U_1|^{\frac{1}{2} - \varepsilon}} |\Xi_2|^{-1+2\varepsilon} 
  \end{split}
  \end{displaymath}
by Stirling's formula.\\

For part (c) we write $${\tt W}(x) = \frac{1}{4\pi^2}  \int_{\Bbb{R}} W(\xi) W\left(\frac{x}{4\pi^2\xi}\right) \frac{d\xi}{\xi}.$$
Then 
$$\mathcal{K}_{\mu}(\Xi_1, \Xi_2; 0, 0; 0, 0) = \int_{-i\infty}^{i\infty}\int_{-i\infty}^{i\infty}  G(s, \mu) S^{\epsilon_1, \epsilon_2}(s, \mu)|  \Xi_1|^{-s_1}|  \Xi_2|^{-s_2} \widehat{ {\tt W}}(1-s_1) \widehat{\overline{ {\tt W}}}(1-s_2)  \frac{ds}{(2\pi i)^2},$$
where as before $(\epsilon_1, \epsilon_2) = (\text{sgn}(\Xi_1), \text{sgn}(\Xi_2))$. 
Here we may shift the $s_1$ and/or $s_2$ contour to the left. The poles at $s_1 = - n + \mu_j$, $s_2 = - n - \mu_j$ for $n \in \Bbb{N}_0$ contribute negligibly by the rapid decay of $\widehat{ {\tt W}}(s)$. By Stirling's formula we see that this forces $|\Xi_1|$ and $|\Xi_2|$ to be $\gg T^{3-\varepsilon}$ for a non-negligible contribution, and we obtain \eqref{c1} for the contours at $\Re s_1 = \Re s_2 = \varepsilon$. 

We proceed to prove \eqref{c2}. Shifting both contours to $\Re s_1 = \Re s_2 = -1/2$ (at the cost of a negligible error), say, we have
\begin{displaymath}
\begin{split}
&\sum_{\epsilon \in \{\pm 1\}^2} \sum_{w \in \mathcal{W}} \sum_D \frac{\phi(D)}{D^2} \mathcal{K}_{w(\mu)}\left(\frac{\Xi_1}{D}, \frac{\Xi_2}{D}; 0, 0; 0, 0\right) \\
&= \sum_{\epsilon \in \{\pm 1\}^2} \sum_{w \in \mathcal{W}} \int_{(-1/2)} \int_{(-1/2)} \frac{\zeta(1 - s_1 - s_2)}{\zeta(2 - s_1 - s_2)} G(s,  \mu) S^{\epsilon_1, \epsilon_2}(s, w(\mu))|  \Xi_1|^{-s_1}|  \Xi_2|^{-s_2} \widehat{ {\tt W}}(1-s_1) \widehat{\overline{ {\tt W}}}(1-s_2)  \frac{ds}{(2\pi i)^2}
\end{split}
\end{displaymath}
plus an error $O(T^{-B})$, where we used the factorization
$$\sum_{D} \frac{\phi(D)}{D^s} = \frac{\zeta(s-1)}{\zeta(s)}$$
for $\Re s > 2$. 
Recall that $G(s, \mu)$ is invariant under the Weyl group. 
The key point is now to shift to the right past the possible pole at $s_1 + s_2 = 0$. This is a very subtle point and sensitive to the signs $(\epsilon_1, \epsilon_2)$. First we notice that $$ G((s, -s), \mu) S^{+ +}((s, -s), \mu) = G((s, -s), \mu) S^{- -}((s, -s), \mu) = 0$$ from Definition \ref{MB}. We consider now
\begin{equation*}
\sum_{w \in \mathcal{W}}\left(S^{+ -}((s, -s), w(\mu)) + S^{- +}((s, -s), w(\mu))\right)  G((s, -s),  \mu).
\end{equation*}
We have
$$S^{+ -}((s, -s),\mu)  G((s, -s),  \mu) = \frac{\cos(\frac{\pi}{2}(\mu_2 - \mu_3))}{32\sin(\frac{\pi}{2}(\mu_1 -\mu_2))\sin(\frac{\pi}{2}(\mu_1 -\mu_3)) (\mu_1 - s)(s - \mu_2)(s - \mu_3)}.$$
This is negligible unless
$$ \Im \mu_2< \Im  \mu_1 < \Im \mu_3 \quad \text{ or } \quad \Im \mu_3 <\Im  \mu_1 <   \Im \mu_2, $$
in which case it equals, up to a negligible error, 
$$\frac{1}{ 16(\mu_1 - s)(s - \mu_2)(s - \mu_3)}.$$
Similarly,
$$S^{-+}((s, -s),\mu)  G((s, -s),  \mu) = -\frac{1}{16 (\mu_1- s)(s - \mu_2)(s - \mu_3)} + O(T^{-B}).$$
if 
$$\Im \mu_1  <  \Im \mu_3  < \Im \mu_2\quad \text{ or } \quad  \Im \mu_2 <  \Im \mu_3  <  \Im \mu_1$$
and $S^{-+}((s, -s),\mu)  G((s, -s),  \mu) = O(T^{-B})$ otherwise. 
Hence for given $\mu$, only two values of $w\in \mathcal{W}$ contribute non-negligibly to the $(+ -)$ term and two values of $w \in \mathcal{W}$ contribute non-negligibly to the $(-+)$ term.  Adding the contributions of the four relevant Weyl chambers gives 0, and we conclude from the previous discussion  that the residue of the integrand at $s_1 + s_2 = 0$ is negligible. We see now why it was important to keep the $\epsilon$-sum and the $w$-sum intact. 
 We can now shift to $\Re s_1 = \Re s_2 = 1/2 - \varepsilon$ at the cost of a negligible error and conclude the desired bound from Stirling's formula.

\subsection{Prelude}\label{secprelude} We proceed to prove (d) and (e). This requires a lot of work and will take the rest of the paper. Here we assume  $|U_1|,  |V_1|, |U_2|,  |V_2| \geq T^{\varepsilon}$. As in the proof of Lemma \ref{lem5}  we integrate over $\xi$ and $\eta$ using Lemma \ref{statphas} and see that   $\mathcal{K}_{\mu}(\Xi_1, \Xi_2; U_1, V_1; U_2, V_2)$ is negligible unless
$$t_1 \asymp U_1 \asymp V_1, \quad t_2 \asymp U_2 \asymp V_2,$$
in which case it equals (up to a negligible error)
\begin{equation}\label{prelude}
 \frac{(\Xi_1 \Xi_2)^{-\varepsilon}}{|U_1V_1U_2V_2|^{1/2}}\int_{\Bbb{R}^2} G(\varepsilon + i t, \mu) S^{\epsilon_1, \epsilon_2} (\varepsilon + i t, \mu)  (\Upsilon_1t_1^2)^{-it_1}  (\Upsilon_2t_2^2)^{-it_2} F\left(\frac{t_1}{U_1}, \frac{t_2}{U_2}\right) dt,
 \end{equation}
where $\Upsilon_j$ is as in \eqref{defw} 
and $F$ is a fixed, smooth, compactly supported function. 

The idea is now to insert Stirling's formula for the Gamma functions and analyze the integral by a two-dimensional stationary phase method. There are two difficulties. First, the quality of Stirling's approximation depends on the distance to the origin, so that we need to  insert several dyadic partitions and treat some ranges trivially. Secondly, due to the complexity of the equations defining the stationary points, we will not actually attempt to locate them, but rather argue that the derivatives of the phase function cannot be too small for too long.  We now make these ideas precise.

We write $\Im \mu = r$. The Gamma functions contained in $G(s, \mu)$ naturally split the integrals in $t_1$ and $t_2$ into intervals based on the signs of $t_1 - r_j$ and $t_2 + r_j$. Because  $|r_j - r_k | \gg T$ for $j \not= k$,  we know that $t_1$ and $-t_2$ can each be close to at most one of the $r_j$. Let $r_j$ be nearest to $t_1$ and $r_k$ nearest to $-t_2$ (possibly $j = k$). This pair of indices is kept fixed for the rest of the argument. 
We introduce a dyadic partition of unity and insert a localizing factor
$$F_1\left(\pm \frac{t_1 - r_j}{B_1}\right)F_2\left(\pm \frac{t_2 - r_k}{B_2}\right)F_3\left(\pm \frac{t_1 + t_2}{B_3}\right),$$
where $B_1, B_2, B_3 > 0$ and $F_1, F_2, F_3$ are smooth, compactly supported functions. We note already at this point that the sizes of $B_1, B_2, B_3, U_1, U_2$ are linked by the triangle inequality and not all combinations are possible. The upper bound contained in Stirling's formula shows that the integrand in \eqref{prelude} is
$$\ll T^{\varepsilon} \frac{1}{ (T+ |U_1|) (T+ |U_2|)}  \left(\frac{1+B_3}{(1+B_1)(1+B_2)}\right)^{1/2}.$$
In order to insert the precise version Stirling's formula that captures the oscillation, we need that $B_1, B_2, B_3 \geq T^{\varepsilon}$.  

We first treat some degenerate cases.  If both $B_1$ and $B_2$ are $\leq T^{\varepsilon}$, then a trivial estimate gives $K_{\mu}(\Xi_1, \Xi_2; U_1, V_1; U_2, V_2) \ll T^{-3/2 + \varepsilon} |U_1V_1U_2V_2|^{-1/2}$ which is stronger than \eqref{11cbound1} and \eqref{11cbound2}.  If $B_3 \leq T^{\varepsilon}$, then $B_1^{1-\varepsilon} \ll B_2 \ll B_1^{1+\varepsilon}$, so that again a trivial estimate gives $K_{\mu}(\Xi_1, \Xi_2; U_1, V_1; U_2, V_2) \ll T^{-2 + \varepsilon} |U_1V_1U_2V_2|^{-1/2}$, which is also stronger than \eqref{11cbound1} and \eqref{11cbound2}.  The only remaining case when, say, only $B_1 \leq T^{\varepsilon}$  is  most efficiently treated on the way, and we will indicate the necessary modifications in due course.  
 
We assume from now on that $B_1, B_2, B_3 \geq T^{\varepsilon}$ and insert Stirling's formula \eqref{stir} for all Gamma factors. In this way we transform the integral in \eqref{prelude} into 
$$\int_{\Bbb{R}^2} \exp\left(-\frac{\pi}{2}h^{\epsilon_1, \epsilon_2}(t, r) \right) e^{i g(t, r)} \mathcal{F}(t, r)  dt,$$
where $h^{\epsilon_1, \epsilon_2}(t, r)$ was defined in \eqref{exp}, 
\begin{displaymath}
\begin{split}
g(t, r) =& - (t_1+t_2) \log|t_1+t_2| - t_1 \log (\Upsilon_1t_1^2) - t_2 \log (\Upsilon_2 t_2^2) \\
&+ \sum_{l=1}^3 \big((t_1 - r_l) \log |t_1 - r_l| + (t_2 + r_l)\log|t_2 + r_l|\big), 
\end{split}
\end{displaymath}
and  the smooth function $\mathcal{F}$ has support in $t_j \asymp U_j$ and satisfies
$$\frac{\partial^n}{\partial t_1^n} \frac{\partial^m}{\partial t_2^m}  \mathcal{F}(t, r) \ll_{n, m} T^{\varepsilon} \frac{1}{ (T+ |U_1|) (T+ |U_2|)}  \left(\frac{B_3}{B_1B_2}\right)^{1/2} \frac{1}{E_1^nE_2^m} $$ 
for $n, m \in \Bbb{N}_0,$ where
$$E_i := \min(B_i, B_3, |U_i|). $$ 
The function $h^{\epsilon_1, \epsilon_2}$ is non-negative and piecewise linear with kink points only at $t_1 = r_l$, $t_2 = -r_l$, $t_1 = -t_2$ for $l \in \{1, 2, 3\}$. As the support of $\mathcal{F}$ is, by definition, away from the kink points, we have either $\exp(-\frac{\pi}{2}h^{\epsilon_1, \epsilon_2}(t, r) ) = 1$ for all $t \in \text{supp}(\mathcal{F})$, or the exponential factor is negligibly small. Hence it suffices to analyze the integral
\begin{equation}\label{integr}
\mathcal{I} = \mathcal{I}(B_1, B_2, B_3) := \int_{\Bbb{R}^2}   e^{i g(t, r)} \mathcal{F}(t, r)  dt,
\end{equation}
for all possible choices of $B_1, B_2, B_3$ to obtain an  upper bound for the quantity $$|U_1U_2V_1V_2|^{1/2} \mathcal{K}_{\mu}(\Xi_1, \Xi_2; U_1, V_1; U_2, V_2)$$
featured in Lemma \ref{lem6}(d). 
For $i = 1, 2$ we write
\begin{equation}\label{der-g}
g_i(t) = g_i(t, r) :=   \frac{\partial }{\partial t_i} g(t, r) = \log \Big| \frac{(t_i + (-1)^i r_1)(t_i + (-1)^i r_2)(t_i + (-1)^i r_3)}{\Upsilon_i t_i^2(t_1+t_2)}\Big|.
\end{equation}
We have
\begin{equation}\label{derbounds}
\begin{split}
 \frac{\partial^n}{\partial t_i^n} g_i(t) & \ll E_i^{-n}, \quad n \geq 1,  \\
 \frac{\partial^m}{\partial t_l^m} \frac{\partial^n}{\partial t_i^n} g_i(t) & \asymp B_3^{-n-m}, \quad n \geq 0, m \geq 1, i \not= l. 
\end{split}
\end{equation}
We compute the first few derivatives explicitly:  let 
$$C_1 = \frac{1}{2}(r_1^2 + r_2^2 + r_3^2) \asymp T^2, \quad C_2 = r_1r_2 r_3 \asymp -T^3,\quad g'' = \text{Hess}(g) =  \left(\begin{matrix} \frac{\partial g_1}{\partial t_1} & \frac{\partial g_1}{\partial t_2} \\ \frac{\partial g_2}{\partial t_1} & \frac{\partial g_2}{\partial t_2} \end{matrix}\right).$$ 
Then a direct computation shows
\begin{equation}\label{p1}
\begin{split}
\frac{\partial g_1}{\partial t_1}(t) = & \frac{P_1(t)}{t_1(t_1 + t_2) \prod_{l=1}^3(t_1 - r_l)},\quad P_1(t) = t_2(t_1^3 + C_1t_1 + 2C_2) + t_1(2C_1 t_1 + 3C_2);\\
\frac{\partial g_2}{\partial t_2}(t) = & \frac{\tilde{P}_1(t)}{t_2(t_1 + t_2) \prod_{l=1}^3(t_2 + r_l)},\quad \tilde{P}_1(t_1, t_2) = P_1(-t_2, -t_1);\\
\det g''(t) = &\frac{P_2(t)}{t_1t_2\prod_{l=1}^3(t_1 - r_l)(t_2 + r_l)},\\
&P_2(t) = 2C_1 t_1^2t_2^2 - 3C_2t_1t_2(t_1 - t_2) + 2C_1^2 t_1t_2 - 4C_1 C_2(t_1 - t_2) - 6C_2^2.\\
\end{split}
\end{equation}
If $\{j, k, l\} = \{1, 2, 3\}$ (in particular $j \not= k$), then
\begin{equation}\label{factP1}
 P_2(r_j, -r_k) = r_jr_k(r_j-r_l)(r_k- r_l) (r_j-r_k)^2 \asymp T^6.
\end{equation}
This can be seen by direct computation, but it is more elegant to observe that
$$\lim_{t \rightarrow (r_j, -r_k)} (t_1 - r_j)(t_2 + r_k) \det g''(t) = 1,$$
which implies \eqref{factP1} immediately. This differs from the case $j = k$ in that $P_2(r_j,-r_j) = 0$. \\
Using $r_1 + r_2 + r_3 = 0$, one also verifies that
\begin{equation}\label{factP2}
 \frac{\partial g_2}{\partial t_2}(t)  = \frac{r_j(t_2 - 2r_lr_m/r_j)}{t_2(t_2+r_l)(t_2 + r_m)} + \frac{t_1 - r_j}{(t_2 + r_j)(t_1 + t_2)} = \frac{r_j(t_2 - 2r_lr_m/r_j)}{t_2(t_2+r_l)(t_2 + r_m)} + O\left(\frac{B_1}{B_2B_3}\right), 
\end{equation}
where $\{j, l, m \} = \{1, 2, 3\}$. \\


Before we proceed, it is convenient to introduce the notation
$$A \lll B$$
to mean $A \leq \delta B$ for a sufficiently small  constant $\delta$ 
  (where ``sufficiently small'' depends on $c, C$ in \eqref{cC}, the support of the weight functions and $\varepsilon$). Similarly we write $A \ggg B$ to mean $A \geq \Delta B$ for a sufficiently large constant $\Delta$. \\

We now introduce another partition of unity that localizes the size of the derivatives $g_1(t)$ and $g_2(t)$,  and insert a factor $F_4(|g_1(t)|/B_4) F_5(|g_2(t)|/B_5)$ into \eqref{integr}.  Applying Lemma \ref{integrationbyparts} with $$Y = Q = E_1, \quad U = E_1 \min(1, B_4),  \quad R = B_4$$ 
or $Y = Q = E_2$, $U = E_2 \min(1, B_5)$,  $R = B_5$, we see that the integral is negligible unless  $$g_i(t) \ll T^{\varepsilon} E_i^{-1/2}.$$ 
In certain situations, this can be refined a little. If 
\begin{equation}\label{ggg1}
|U_1| \ggg T+|U_2|, \quad |U_1||U_2| \ggg T^2,
\end{equation}
 then we can compute explicitly 
\begin{equation}\label{varder}
\begin{split}
\frac{(-1)^n}{(n-1)!} \frac{\partial^n g_1}{\partial t_1^n}(t) &= \sum_{l=1}^3 \frac{1}{(t_1 - r_l)^n} - \frac{2}{t_1^n} - \frac{1}{(t_1 + t_2)^n}\\
& = \frac{1 + O((T^2+|U_2|^2)|U_1|^{-2})}{t_1^{n-1}} \Bigl(\sum_{l=1}^3 \frac{1}{t_1 - nr_l} - \frac{2}{t_1} - \frac{1}{t_1 + nt_2}\Bigr) \\
&= \frac{nt_1^3 t_2(1 + O((T^2 + |U_2|^2)|U_1|^{-2}+T^2|U_1 U_2|^{-1}))}{t_1^n(t_1 - nr_1)(t_1 - nr_2)(t_1 - nr_3)(t_1 + nt_2)} \asymp \frac{|U_2|}{|U_1|^{n+1}}.
\end{split}
\end{equation}
for $n \in \Bbb{N}$.  (Here we used again that $r_1 + r_2 + r_3 = 0$.) Then we can apply Lemma  \ref{integrationbyparts} similarly, but with $Y = |U_2|$, $Q = |U_1|$, and obtain that the integral is negligible unless  
\begin{equation}\label{better}
g_1(t) \ll   |U_2|^{1/2 }|U_1|^{-1+\varepsilon}.
\end{equation}
Notice that depending on the value of $\varepsilon$, \eqref{varder} is needed for some $n \leq n_0 = n_0(\varepsilon)$, see the remark after Lemma \ref{integrationbyparts}. This $n_0$ determines the implicit constants in \eqref{ggg1}. 

A similar statement holds with exchanged indices: if
\begin{equation}\label{ggg2}
|U_2| \ggg T+|U_1|, \quad |U_1||U_2| \ggg T^2,
\end{equation}
then the integral is negligible unless  
\begin{equation}\label{better2}
g_2(t) \ll   |U_1|^{1/2 }|U_2|^{-1+\varepsilon}.
\end{equation}

With this in mind, let us define
 $\mathcal{A}_{j, k} = \mathcal{A}_{j, k}(U_1, U_2, B_1, B_2, B_3)$ to be the set of all $t = (t_1, t_2)$ satisfying
\begin{equation}\label{Ajk}
  g_i(t) \ll T^{\varepsilon} E_i^{-1/2}, \quad t_i \asymp U_i, \quad |t_1 - r_j| \asymp B_1, \quad |t_2 + r_k| \asymp B_2, \quad |t_1+t_2| \asymp B_3,
\end{equation}
except in the situation where \eqref{ggg1} holds, in which case we replace the first condition for $i=1$ with \eqref{better}, and in the analogous situation where \eqref{ggg2} holds, in which case we replace the first condition for $i=2$ with \eqref{better2}. Notice that even in these exceptional cases \eqref{Ajk} still holds.

The set $\mathcal{A}_{j, k}$ is the subset of the support of $\mathcal{F}$ where stationary points can lie, and the above discussion shows
\begin{equation}\label{mainint}
 \mathcal{I}   \ll   \frac{T^{\varepsilon}}{ (T+ |U_1|) (T+ |U_2|)}  \left(\frac{B_3}{B_1B_2}\right)^{1/2} \text{meas}(\mathcal{A}_{j, k}),
\end{equation}
 up to a negligible error. We note that in particular $\mathcal{A}_{j, k}$ is empty unless the consistency relation 
\begin{equation}\label{consist}
(T + |U_i|)^2 B_i \asymp \Upsilon_i |U_i|^2 B_3, \quad i = 1, 2 
\end{equation}
holds, see \eqref{der-g}. 

\subsection{Proof of part (e)}
At this point we can already complete the proof of (e). If $|U_2| \gg (|U_1| + T) T^{\varepsilon}$, then
$$g_2(t) = - \log \Upsilon_2 - \log\left| 1 + \frac{t_1}{t_2}\right| + \sum_{l=1}^3 \log\left|1 + \frac{r_l}{t_2}\right| = - \log \Upsilon_2 + O\left(\frac{|U_1| + T}{|U_2|}\right).$$
Since   $B_2, B_3 \asymp |U_2|$, we have $E_2 = |U_2|$, so that $\mathcal{A}_{j, k}$ is empty unless 
$$|\Upsilon_2 - 1| \ll T^{\varepsilon}\left(\frac{|U_1| + T}{|U_2|} + \frac{1}{|U_2|^{1/2}}\right)$$
which in view of $|U_2| \asymp |V_2|$ is equivalent to  \eqref{cond11c}. Notice that for this argument we do not need to insert Stirling's formula for $\Gamma(s_1 - \mu_j)$, so that the argument works even in the previously excluded case $B_1 \leq T^{\varepsilon}$, but $B_2, B_3$ are automatically $\gg T^{\varepsilon}$, hence the proof of (e) is complete. 

 
\subsection{The measure of the critical set} Obviously, the crucial point in the estimation of \eqref{mainint} is the size of the set $\mathcal{A}_{j, k}$, i.e.\ an estimate how long the derivatives can be small.   The following lemma  gives a first result in this direction. 

\begin{sublemma}\label{sub1}  Let $M \subseteq \Bbb{R}$ be an interval. Then we have the following bounds: 
\begin{equation}\label{1a}
 \text{{\rm meas}}(\mathcal{A}_{j, k} \cap (M \times \Bbb{R})) \ll  T^{\varepsilon} \frac{\min(\text{{\rm meas}}(M), B_1, |U_1|)}{\sqrt{E_2} |\frac{\partial g_2}{\partial t_2} (t)|}, \quad \text{provided } \Big|\frac{\partial g_2}{\partial t_2} (t)\Big| \gg \frac{T^{\varepsilon}}{E_2^{5/4}};
\end{equation}
\begin{equation}\label{1b}
 \text{{\rm meas}}(\mathcal{A}_{j, k} )  \ll T^{\varepsilon} \frac{\min(B_1, |U_1|)B_3}{\sqrt{E_1} }; 
\end{equation}
\begin{equation}\label{2a}
 \text{{\rm meas}}(\mathcal{A}_{j, k} \cap ( \Bbb{R} \times M)) \ll T^{\varepsilon} \frac{\min(\text{{\rm meas}}(M), B_2, |U_2|)}{\sqrt{E_1} |\frac{\partial g_1}{\partial t_1} (t)|}, \quad \text{provided } \Big|\frac{\partial g_1}{\partial t_1} (t)\Big| \gg \frac{T^{\varepsilon}}{E_1^{5/4}};
\end{equation}
\begin{equation}\label{2b}
 \text{{\rm meas}}(\mathcal{A}_{j, k}  ) \ll  T^{\varepsilon} \frac{\min(  B_2, |U_2|)B_3}{\sqrt{E_2} }; 
\end{equation}
\begin{equation}\label{33}
  \text{{\rm meas}}(\mathcal{A}_{j, k} ) \ll T^{\varepsilon} \frac{1}{\sqrt{E_1E_2} |\det g''(t)|}, \quad \text{provided }  |\det g''(t)| \gg T^{\varepsilon} \frac{\max(E_1, E_2)^{1/6}}{(E_1E_2)^{7/6}}.
\end{equation}
If in addition $|U_1| \ggg T + |U_2|$ and $|U_1U_2| \ggg T^2$, then we have
\begin{equation}\label{modif1}
\text{{\rm meas}}(\mathcal{A}_{j, k}) \ll |U_2|^{-1/2  } B_2 |U_1|^{1+\varepsilon};
\end{equation}
similarly, if $|U_2| \ggg T + |U_1|$ and $|U_1U_2| \ggg T^2$, then we have
\begin{equation}\label{modif2}
 \text{{\rm meas}}(\mathcal{A}_{j, k}) \ll |U_1|^{-1/2  } B_1 |U_2|^{1+\varepsilon}.
\end{equation}
\end{sublemma}

\textbf{Proof.}  We start with the remark that the number of connected components of $\mathcal{A}_{j, k}$  is absolutely bounded. Indeed, up to changing implied constants,   $\mathcal{A}_{j, k}$ is   the intersection of the preimages of intervals (around $\pm 1$) under a pair of rational functions, since we  can re-write the first condition   in the definition \eqref{Ajk} of $\mathcal{A}_{j, k}$ as
\begin{equation}\label{rewrite}
\Bigl|\frac{(t_i + (-1)^i r_1)(t_i + (-1)^i r_2)(t_i + (-1)^i r_3)}{\Upsilon_i t_i^2(t_1+t_2)} - \alpha_i \Bigr| \ll T^{ \varepsilon} E_i^{-1/2}, \quad i = 1, 2, \quad \alpha_i \in \{\pm 1\}
\end{equation}
(and similarly with a modified right hand side under the extra condition \eqref{ggg1} or \eqref{ggg2}).  
We call the corresponding subsets $\mathcal{A}_{j, k}^{\alpha_1, \alpha_2}$. 
This in turn can  be expressed as the intersection of sets where certain polynomials are positive, since we have in general
$$\Bigl| \frac{p_1(t)}{p_2(t)} \pm 1\Bigr| < A \Longleftrightarrow A^2 p_2(t)^2 - (p_1(t) \pm p_2(t))^2 > 0.$$
A theorem of Milnor and Thom \cite[Theorem 3]{Mi} (see also \cite{Th})  gives an absolute bound for the number of connected components in terms of the degrees of the polynomials, which for us are fixed. Inspecting the defining conditions \eqref{Ajk}, we see that $\mathcal{A}_{j, k}$ is a union of four sets (characterized by two choices of $\pm 1$) each of which can be characterized 
by  two inequalities of degree 6 and $4 + 12 + 12 + 4 = 32$ linear inequalities, so that by \cite[Theorem 3]{Mi}, $\mathcal{A}_{j, k}$  has at most $4\cdot 1035$ connected components.  Therefore it suffices to prove the bounds for the measure of each connected component of $\mathcal{A}_{j, k}$. \\

We start with the proof of \eqref{33}. In this proof we keep track of the value of $\varepsilon$ and do not apply $\varepsilon$-convention. To this end we remark that we can assume that $E_i  = \min(B_i, B_3, |U_i|) \gg T^{10\varepsilon}$, say. 

 Let $t \in \mathcal{A}_{j, k}$ and consider the connected component $\mathcal{A}_{j, k}(t)$ of $\mathcal{A}_{j, k}$ containing $t$.  Let $$\mathcal{B} := g''(t)^{-1}   \text{diag}\left(E_1^{1/2}, E_2^{1/2}\right) ^{-1}(\mathcal{D}_{2\varepsilon}),$$ where $\mathcal{D}_{2\varepsilon}$ is the \changed{open} disk of radius $T^{2\varepsilon}$ centered at 0.  The idea is to show that no point of the boundary $\changed{t+}\partial \mathcal{B}$ can be contained in $\mathcal{A}_{j, k}(t)$, so that $\mathcal{A}_{j, k}(t) \subseteq \changed{t+}\mathcal{B}$. 
 
 Suppose that $t+u \in \mathcal{A}_{j, k}(t)$. Then by Taylor's theorem  and the bounds \eqref{derbounds} we have
$$g_1(t+u) - g_1(t) = u_1 \frac{\partial g_1}{\partial t_1}(t) + u_2 \frac{\partial g_1}{\partial t_2}(t) + O\left(\frac{|u_1|^2}{E_1^2} + \frac{|u_1u_2|}{B_3^2} + \frac{|u_2|^2}{B_3^2}\right).$$ 
The error term is majorized by $O(|u_1|^2 E_1^{-2} + |u_2|^2 \max(E_1, E_2)^{-2})$. An analogous formula holds for $g_2(t+u) - g_2(t)$.  
Hence if $u_1 \ll T^{\varepsilon} E_1^{2/3}$ and $u_2 \ll T^{\varepsilon} E_2^{2/3}$, we conclude from \eqref{Ajk} that 
\begin{equation}\label{contra}
T^{2\varepsilon} \gg E_1(g_1(t+u) - g_1(t))^2 + E_2(g_2(t+u) - g_2(t))^2 = \left\| \text{diag}\left(E_1^{1/2}, E_2^{1/2}\right) g''(t) u\right\|^2 + O(T^{3\varepsilon})
\end{equation}
by \eqref{derbounds}.    \changed{If also $u$ is in $\partial\mathcal{B}$, we have for suitable $0 \leq \theta < 2\pi$ that}
$$u_1 = \frac{1}{|\det g''(t)|} \left( T^{2\varepsilon} \frac{\cos\theta}{\sqrt{E_1}} \frac{\partial g_2}{\partial t_2}(t) -  T^{2\varepsilon} \frac{\sin\theta}{\sqrt{E_1}} \frac{\partial g_1}{\partial t_2}(t)\right) \ll \frac{T^{2\varepsilon}}{\sqrt{E_1} E_2 |\det g''(t)|} \ll T^{\varepsilon} E_1^{2/3},$$
provided the assumption of \eqref{33} is satisfied. Similarly, 
$u_2 \ll 
T^{\varepsilon} E_2^{2/3}$,  
so that  \eqref{contra} can be applied. On the other hand, for such $u$,  the main term on the right hand side of \eqref{contra} equals $T^{4\varepsilon}$, a contradiction.
\changed{The matrices $\text{diag}\left(E_1^{1/2}, E_2^{1/2}\right)^{-1}$ and $g''(t)^{-1}$ are clearly invertible, so that $t+\mathcal{B}$ is open and intersects $\mathcal{A}_{j, k}(t)$ at least at $t$.
By the above, $\mathcal{A}_{j, k}(t)$ does not intersect $t+\partial\mathcal{B}$, hence it is covered by $t+\mathcal{B}$ and the complement of its closure, $t+(\mathbb{R}^2\setminus\bar{\mathcal{B}})$.}
Hence $\mathcal{A}_{j, k}(t) \subseteq \changed{t+}\mathcal{B}$ \changed{by connectivity}, and the volume of $\changed{t+}\mathcal{B}$ is $\pi T^{4\varepsilon}(\sqrt{E_1E_2} |\det g''(t)|)^{-1}$. 

The proofs of the remaining statements are similar, but simpler. If $u = (u_1, u_2)$ with $u_2 = 0$ is such that  $t+u \in \mathcal{A}_{j, k}(t)$, then 
\begin{displaymath}
\begin{split}
\frac{T^{\varepsilon} }{E_1^{1/2}} \gg  g_1(t+u) - g_1(t) = u_1 \frac{\partial g_1}{\partial t_1}(t) + O\left(\frac{|u_1|^2}{E_1^2}\right),\\
\frac{T^{\varepsilon}}{ E_2^{1/2} }  \gg g_2(t+u) - g_2(t) = u_1 \frac{\partial g_2}{\partial t_1}(t) + O\left(\frac{|u_1|^2}{B_3^2}\right).
\end{split}
\end{displaymath}
Provided $|u_1| \ll E_1^{3/4}$, we see that $u_1 = \pm T^{2\varepsilon} E_1^{-1/2} 
| \frac{\partial g_1}{\partial t_1}(t)|^{-1}$ leads to a contradiction in the first inequality. Estimating the possible range for $u_2$ trivially from the conditions $t_2 \asymp U_2$, $|t_2 + r_k| \asymp B_2$, $t_2 \in M$  defining the set $\mathcal{A}_{j, k} \cap (\Bbb{R}\times M)$ proves \eqref{2a}. 

Similarly, using the second line in \eqref{derbounds} with $m=1$ and $n=0$, we arrive at a contradiction for the second inequality upon choosing $u_1 = \pm T^{2\varepsilon} B_3 E_2^{-1/2}$ (which is always less than $B_3$), and we estimate the range for $u_2$ trivially to obtain \eqref{2b}.  

If in addition \eqref{ggg1} holds, we can replace the first inequality with
$$\frac{|U_2|^{1/2 + \varepsilon} }{|U_1|} \gg  g_1(t+u) - g_1(t) = u_1 \frac{\partial g_1}{\partial t_1}(t) + O\left(\frac{|u_1|^2|U_2|}{|U_1|^3}\right)$$
by \eqref{varder}. 
Since under the present assumptions $|\frac{\partial g_1}{\partial t_1}(t)|  \asymp |U_2| |U_1|^{-2}$, we arrive at a contradiction by choosing $u_1 = |U_2|^{-1/2 + 2\varepsilon} |U_1|$ (provided $\varepsilon < 1/4$). Estimating the $u_2$ range trivially, we obtain \eqref{modif1}. 

The proofs of \eqref{1a}, \eqref{1b} and \eqref{modif2} are identical with exchanged indices. \hfill $\square$\\

We emphasize that the proofs of \eqref{2a} and \eqref{modif1}  make no use of the function $g_2$, and the proofs of \eqref{1a} and \eqref{modif2} make no use of the function $g_1$. In particular, the latter two bounds  can be used in the following section where we treat the case when $B_1$ is small and Stirling's formula cannot be inserted. 

\subsection{The case where $B_1$ is small} \label{B1small}We are now prepared to treat the remaining exceptional case where, say, $B_1 \leq T^{\varepsilon}$. This implies in particular $|U_1| \asymp T$. We distinguish several cases depending on the size of $|U_2|$. The following analysis is already a precursor to the various cases below.\\

\emph{Case 1:} Suppose that $|U_2| \ggg T$, so that $B_2 \asymp B_3 \asymp E_2 \asymp |U_2|$. In this case \eqref{modif2} with $B_1 \leq T^{\varepsilon}$ is applicable, and from \eqref{mainint} we obtain
$$\mathcal{I} \ll \frac{ T^{\varepsilon} }{ (T+ |U_1|) (T+ |U_2|)}  \left(\frac{B_3}{B_2}\right)^{1/2} \frac{|U_2|}{T^{1/2}} \ll T^{-3/2 + \varepsilon}.$$

\emph{Case 2:} Suppose that $|U_2| \asymp T$. Then $B_2, B_3 \ll T$. If $B_2 \leq T^{9/10}$ or $B_3 \leq T^{9/10}$, we can estimate trivially $\text{meas}(\mathcal{A}_{j, k}) \ll B_1B_2 \leq T^{\varepsilon} B_2$, so that 
$$\mathcal{I} \ll T^{\varepsilon} \frac{1}{T^2} \left(\frac{B_3 }{B_2 }\right)^{1/2} \cdot B_2 \leq T^{-1 - \frac{1}{20} + \varepsilon}.$$
Suppose from now on that $B_2, B_3 \geq T^{9/10}$. By \eqref{factP2} we have 
$$\frac{\partial g_2}{\partial t_2}(t) = \frac{r_j(t_2 - 2r_lr_m/r_j)}{t_2(t_2+r_l)(t_2 + r_m)} + O(T^{-9/5 + \varepsilon})$$
for $\{j, l, m\} = \{1, 2, 3\}$ (here $k$ may or may not be different from $j$).     If 
$|t_2 - 2 r_lr_m/r_j|\geq T^{9/10},$ 
then $$\Bigl|\frac{\partial g_2}{\partial t_2}(t)\Bigr| \gg  T^{-\frac{11}{10}} \geq T^{-\frac{5}{4} \cdot \frac{9}{10} + \varepsilon} \geq E_2^{-\frac{5}{4} + \varepsilon},$$ so that \eqref{1a} with $M = \Bbb{R}$ implies
$\text{meas}(\mathcal{A}_{j, k}) \ll T^{11/10+\varepsilon} \min(B_2, B_3)^{-1/2} $
and hence
$$\mathcal{I} \ll T^{\varepsilon} \frac{ T^{11/10} \min(B_2, B_3)^{-1/2} }{T^{2}} \left(\frac{B_3}{ B_2}\right)^{1/2} \ll T^{-\frac{13}{10} +  \varepsilon} .$$
In the opposite case, when $|t_2 - 2 r_lr_m/r_j|\leq T^{9/10}$,  we estimate trivially
$$\mathcal{I} \ll T^{\varepsilon} \frac{1}{T^2} \left(\frac{B_3 }{B_2 }\right)^{1/2} T^{9/10} \leq T^{-1-\frac{1}{20} + \varepsilon}.$$

\emph{Case 3:} Finally suppose $|U_2| \lll T$, so that $B_2 \asymp B_3   \asymp T$, $E_2 \asymp |U_2|$.  Then \eqref{p1} implies that $\tilde{P}_1(t) \asymp T^4$ (since $|t_1| \asymp |U_1| \asymp T$), so that $|\frac{\partial g_2}{\partial t_2}(t)|  \asymp |U_2|^{-1}$. The condition of \eqref{1a} is trivially satisfied, and we obtain
$$\mathcal{I} \ll T^{\varepsilon}  \frac{1}{T^{2}} \left(\frac{B_3}{ B_2}\right)^{1/2}  |U_2|^{1/2} \leq T^{-3/2 + \varepsilon}.$$

This proves a strong version of \eqref{11cbound1} and \eqref{11cbound2} in all  cases and completes the discussion of the case $B_1 \leq T^{\varepsilon}$.

\subsection{The nearly generic case}\label{secgeneric}

We are now prepared for the proof of \eqref{11cbound1} and \eqref{11cbound2} in the situation where $B_1, B_2, B_3 \geq T^{\varepsilon}$. 
This will be a case-by-case analysis.
We first consider the case
\begin{equation}\label{case1a}
T^{1-b} \leq |U_1|, |U_2|, B_1, B_2, B_3 \leq T^{1+b},
\end{equation}
where $0 < b < \frac{1}{50}$ is a small constant.
We will later show that the exponent of \eqref{11cbound2} results from the choice $b=\Oneb$.
Our assumption implies $T^{1-b} \leq E_1, E_2 \leq T^{1-b}$, and from the consistency relation \eqref{consist}, we also have $T^{-b} \leq \Upsilon_l \leq T^{3b}$, $l=1,2$, with the upper bound occuring at $\abs{U_l} = B_3 = T^{1-b}$, $B_l= T$ and the lower bound at $\abs{U_l} = B_3 = T$, $B_l= T^{1-b}$.

In Sublemma \ref{sub1}, every part but \eqref{33} is concerned with saving the square-root of the length in a single dimension, but this just fails to give a power saving.
Similarly, we cannot use \eqref{33}, because we do not have good control over the Hessian $\det g''(t)$.
In this section, we give a refinement of \eqref{1b}, essentially by assuming that any given component has a highly degenerate singular point.

As in the proof of Sublemma \ref{sub1} it suffices to bound the measure of the connected component $\mathcal{A}_{j, k}^{\alpha_1, \alpha_2}(t)$ containing some $t \in \mathcal{A}_{j, k}^{\alpha_1, \alpha_2}$.  The conditions \eqref{rewrite} are equivalent to (recall that $r_1+r_2 + r_3 = 0$)
\begin{displaymath}
\begin{split}
&(\Upsilon_1^{-1} - \alpha_1)t_1^3 - \Upsilon_1^{-1}C_1 t_1 - \Upsilon_1^{-1} C_2 - \alpha_1 t_1^2t_2 \ll T^{\varepsilon} B_3 |U_1|^2 E_1^{-1/2},\\
& (\Upsilon_2^{-1} - \alpha_2)t_2^3 - \Upsilon_2^{-1}C_1 t_2 - \Upsilon_2^{-1} C_2 - \alpha_2 t_2^2t_1 \ll T^{\varepsilon} B_3 |U_2|^2 E_2^{-1/2}.
\end{split}
\end{displaymath}
Rearranging the first, we have
\begin{equation}\label{t2}
t_2 = \alpha_1 \frac{(\Upsilon_1^{-1} - \alpha_1) t_1^3 - \Upsilon_1^{-1} C_1 t_1 - \Upsilon_1^{-1} C_2}{t_1^2} + O(T^{\varepsilon} B_3 E_1^{-1/2}).
\end{equation}
By assumption, $t_2 \asymp U_2$ and the error term is trivially $O(T^{\varepsilon}|U_2|)$. 
We substitute this into the second equation, and after some algebra, in particular clearing the denominators using $t_1 \asymp U_1$, we obtain have a non-trivial bound on a polynomial of the form
\begin{align*}
	\sum_{i=0}^9 a_i t_1^i  = a_9\prod_{i=1}^9 (t_1 - q_i)  \ll& T^{\varepsilon} B_3 |U_1|^6 \left(E_1^{-1/2}((1 + \Upsilon_2^{-1})|U_2|^2 + \Upsilon_2^{-1} T^2 + |U_1U_2|) + E_2^{-1/2}|U_2|^2\right) \\
	\ll & T^{\frac{17}{2}+\frac{21}{2}b+\varepsilon}
\end{align*}
for some complex numbers $a_i$, $q_i$, independent of $t$, where in particular
\begin{displaymath}
\begin{split}
a_9 &= (\Upsilon_1^{-1} - \alpha_1)^2 (\Upsilon_1^{-1}\Upsilon_2^{-1} -\alpha_2\Upsilon_1^{-1} - \alpha_1\Upsilon_2^{-1}),\\
a_8 & = 0,\\
a_7 & = -C_1(\Upsilon_1^{-1} - \alpha_1) (\alpha_1\alpha_2 \Upsilon_1^{-1} - 3\alpha_2 \Upsilon_1^{-2} + \Upsilon_2^{-1} - 3\alpha_1 \Upsilon_1^{-1} \Upsilon_2^{-1} + 3\Upsilon_1^{-2} \Upsilon_2^{-1}),\\
a_6 & = - C_2(-\alpha_2\Upsilon_1^{-1} + 4\alpha_1\alpha_2 \Upsilon_1^{-2} - 3\alpha_2 \Upsilon_1^{-3} - \alpha_1 \Upsilon_2^{-1} + 3\Upsilon_1^{-1} \Upsilon_2^{-1} - 6\alpha_1 \Upsilon_1^{-2} \Upsilon_2^{-1} + 3\Upsilon_1^{-3} \Upsilon_2^{-1}).
\end{split}
\end{displaymath}

First, assume that $|\Upsilon_1^{-1} - \alpha_1| \geq T^{-a}$ and $|\Upsilon_1^{-1}\Upsilon_2^{-1} - \alpha_2\Upsilon_1^{-1} - \alpha_1 \Upsilon_2^{-1}| \geq T^{-a}$, where $\frac{1}{4} > a > 6b$ is a constant to be chosen in a moment.
We conclude that
$$|t_1 - q_i|^9 \ll \frac{T^{\frac{17}{2}+\frac{21}{2}b+\varepsilon}}{|a_9|} \leq T^{\frac{17}{2}+3a+\frac{21}{2}b+\varepsilon}, $$
for some $i \in \{1, \ldots, 9\}.$
Since $t_1$ is now in a fixed interval, independent of $t_2$, we may apply \eqref{t2} to obtain
\begin{equation}\label{lem2-1}
\text{{\rm meas}}(\mathcal{A}_{j, k}^{\alpha_1, \alpha_2} ) \ll T^{\frac{17}{18}+\frac{a}{3}+\frac{7}{6}b+\varepsilon} B_3 E_1^{-1/2} \leq T^{\frac{13}{9}+\frac{a}{3}+\frac{8}{3}b+\varepsilon}.
\end{equation}

If instead $|\Upsilon_1^{-1} - \alpha_1| \leq T^{-a}$, so that necessarily $\alpha_1 = 1$ and $\Upsilon_1^{-1} = 1 + O(T^{-a})$, then the coefficients simplify
\begin{displaymath}
\begin{split}
  a_9 &= (\Upsilon_1^{-1} - 1)^2\left(-\alpha_2 + O(T^{-a} (1 + \Upsilon_2^{-1}))\right) \ll T^{-2a} (1 + T^{-a} \Upsilon_2^{-1}) \ll T^{-2a},\\
  a_7 & = -C_1(\Upsilon_1^{-1} - 1)\left(-2\alpha_2 + \Upsilon_2^{-1} + O(T^{-a} (1 + \Upsilon_2^{-1}))\right) \ll T^{2-a+b},\\
  a_6 &= -C_2\left(-\Upsilon_2^{-1} + O(T^{-a}(1 + T^{-a}\Upsilon_2^{-1}))\right) \asymp - T^3 \Upsilon_2^{-1} \geq T^{3-3b}
    \end{split}
\end{displaymath}
(the assumptions on $a$ imply $T^b\geq \Upsilon_2^{-1} \geq T^{-3b} \geq T^{-2a}$), so
\begin{equation*}
\begin{split}
\sum_{i=0}^6 a_it^i \ll T^{\varepsilon} \left(T^{\frac{17}{2}+\frac{21}{2}b}+ |U_1|^9 T^{-2a}  + |U_1|^7 T^{2-a+b}\right) 
	\ll T^{9-a+8b+\varepsilon},
\end{split}
\end{equation*}
and we apply the same reasoning as before to obtain
\begin{equation}\label{lem2-2}
\begin{split}
\text{{\rm meas}}(\mathcal{A}_{j, k}^{1, \alpha_2}) \ll & T^{1-\frac{a}{6}+\frac{11}{6}b+\varepsilon} B_3 E_1^{-1/2} \leq T^{\frac{3}{2}-\frac{a}{6}+\frac{10}{3}b+\varepsilon}.
\end{split}
\end{equation}

Finally, if  $|\Upsilon_1^{-1}\Upsilon_2^{-1} - \alpha_2\Upsilon_1^{-1} - \alpha_1 \Upsilon_2^{-1}| \leq T^{-a}$, so that $\Upsilon_1^{-1} \Upsilon_2^{-1} = \alpha_2\Upsilon_1^{-1} + \alpha_1\Upsilon_2^{-1} + O(T^{-a})$, then
\begin{displaymath}
\begin{split}
  a_9 &\ll (\Upsilon_1^{-1} - \alpha_1)^2 T^{-a} \ll T^{-a+2b},\\
  a_7 & = -C_1\left(\alpha_1\alpha_2 \Upsilon_1^{-2} + O((1 + \Upsilon_1^{-1} + \Upsilon_1^{-2})T^{-a})\right) \asymp - \alpha_1\alpha_2 C_1 \Upsilon_1^{-2} \gg T^{2-6b},
\end{split}
\end{displaymath} 
(recall that we are assuming $\Upsilon_1^{-2} \geq T^{-6b} \ggg T^{-a  }$), so
\begin{displaymath}
\begin{split}
\sum_{i=0}^7 a_it^i \ll  T^{\varepsilon} \left( T^{\frac{17}{2}+\frac{21}{2}b}+ |U_1|^9   T^{-a+2b}\right) \ll T^{\frac{17}{2}+\frac{21}{2}b+\varepsilon}+ T^{9-a+11b+\varepsilon}, 
\end{split}
\end{displaymath}
and again we apply the same reasoning to obtain
\begin{equation}\label{lem2-3}
\begin{split}
  \text{{\rm meas}}(\mathcal{A}_{j, k}^{\alpha_1, \alpha_2}) \ll T^{1-\frac{a}{7}+\frac{17}{7}b+\varepsilon} B_3 E_1^{-1/2} \leq T^{\frac{3}{2}-\frac{a}{7}+\frac{55}{14}b+\varepsilon}.
 \end{split}
\end{equation}

Since $b < 1/50$, the choice $a := (7+159b)/60$ satisfies $\frac{1}{4} > a > 6b$, and combining \eqref{lem2-1} -- \eqref{lem2-3}, we obtain in all cases
$$\text{meas}(\mathcal{A}_{j, k}^{\alpha_1, \alpha_2}) \ll T^{\frac{3}{2} - \frac{1}{60} + \frac{71}{20}b + \varepsilon}.$$
After substituting into \eqref{mainint} we obtain
\begin{equation}\label{lem2-all}
\mathcal{I} \ll  \frac{\text{meas}(\mathcal{A}_{j, k}) }{T^{2-\varepsilon}} \left(\frac{B_3}{B_1B_2}\right)^{1/2}\ll T^{-1-\frac{1}{60} + \frac{91}{20}b + \varepsilon}.
\end{equation}

\subsection{Another  special case}\label{sec7} Before proceeding to the general cases, it is convenient to deal with the following special case:
\begin{equation}\label{case7}
|U_1| \geq T^{1+b},  \quad T^{1- b/4} \leq |U_2| \leq T^{1 + b/4}.
\end{equation}
Here we have $B_1 \asymp B_3 \asymp |U_1|$. 
 By \eqref{2b} we conclude $\text{meas}(\mathcal{A}_{j, k}) \ll T^{\varepsilon} |U_1| \min(B_2, |U_2|)^{1/2}$, so that
$$\mathcal{I} \ll \frac{T^{\varepsilon} |U_1| \min(B_2, |U_2|)^{1/2}}{|U_1| (T+|U_2|) B_2^{1/2}} \leq  T^{-1+\varepsilon}.$$
Moreover, by \eqref{consist} we have 
$$\Upsilon_2 \asymp \frac{(T+ |U_2|)^2 B_2 }{|U_1| |U_2|^2} \ll \frac{(T+ |U_2|)^3   }{|U_1| |U_2|^2} \ll \frac{T^{1+b/2}}{|U_1|} \leq T^{-b/2}.$$
 Since   $E_1  \asymp |U_1|$, it follows from
$$T^{\varepsilon} |U_1|^{-1/2}   \gg g_1(t) = \sum_{l=1}^3 \log\left|1 - \frac{r_j}{t_1}\right| - \log\left|1 - \frac{t_2}{t_1}\right| - \log \Upsilon_1 = - \log \Upsilon_1 + O\left(T^{1 + b/4 }|U_1|^{-1}\right)$$
that $$|\Upsilon_1 - 1| \ll T^{1+\frac{b}{4}  } |U_1|^{-1} + T^{\varepsilon} |U_1|^{-1/2}\ll T^{-b/2}.$$
Again, we will choose $b=\Oneb$, so this suffices for the proof of \eqref{11cbound1}, and of course the same argument works with exchanged indices. \\

Having this case out of the way, we will show the bound
\begin{equation}\label{toshow}
  \mathcal{I} \ll T^{-1-\frac{b}{8} +\varepsilon}
\end{equation}
in all other cases; choosing $b=\Oneb$ here and in \eqref{lem2-all} then gives \eqref{11cbound2}. To this end we distinguish the following principal cases
\begin{displaymath}
\begin{split}
& (1) \quad |U_1| \asymp |U_2| \asymp T, \quad\quad \,\,\,  (2) \quad  T^{1+b} \geq |U_1|  \ggg  T \asymp |U_2|, \quad \quad   (3)  \quad  |U_1|  \lll  T \asymp |U_2|,  \\
& (4) \quad  |U_1| \geq |U_2| \ggg T,  \quad\quad (5) \quad   |U_1| \ggg T \ggg  |U_2|,   \quad\quad\quad\quad\quad  (6) \quad |U_1| \leq |U_2| \lll T \end{split}
\end{displaymath}
with the understanding that those situations covered by \eqref{case1a} and \eqref{case7} and its version with exchanged indices are excluded.  By symmetry, this covers all possibilities. 

\subsection{Case 1} The present assumption  $|U_1| \asymp |U_2| \asymp T$  implies $B_1, B_2, B_3 \ll T$. 
We distinguish the following subcases:
\begin{displaymath}
\begin{split}
& (1a) \quad B_3 \leq T^{1-\frac{b}{5}}, \qquad  (1b) \quad B_1 \leq  B_2 \lll T \asymp B_3,   \qquad (1c) \quad B_1 \leq T^{1-b}, \,  T^{1-\frac{b}{5}} \leq B_3 \asymp B_2. 
\end{split}
\end{displaymath}
They cover all cases where $B_1 \leq B_2$ and \eqref{case1a} does not hold (the   cases with $B_2 \leq B_1$ being completely analogous). Indeed, if $B_3 \leq T^{1 - b/5}$, we are in Case 1a. If $B_3 \asymp T$, we are in Case 1b if $B_2 \lll B_3$ and in Case 1c if $B_2 \asymp B_3$. If $T^{1-b/5} \leq B_3 \lll T$ and in addition $B_1 \geq T^{1-b}$ (which we can assume by \eqref{case1a}), then by the triangle inequality $B_2 \asymp B_3$, and we are again in Case 1c. 

Before we give a detailed analysis, we explain briefly why we obtain a non-trivial bound in all cases. In Case 1a, the trivial bounds suffices, because small $B_3$ is advantageous in \eqref{mainint}. In Case 1b we know that $t_1 \approx r_j$ and $t_2\approx - r_k$, which lets us control the size of $P_2(t)$ defined in \eqref{p1}, and hence the size of the second derivative $g''(t)$. Case 1c is a bit more difficult. In typical situations we have control over the first derivative $\partial g_2/\partial t_2$ since we know that $t_1 \asymp r_j$. There are certain degenerate configurations, however, but they restrict  $t_2$ to a small interval, so that then a simple bound suffices.  \\

\emph{Case 1a:} We may assume without loss of generality that $B_1 \leq B_2$. It follows from the relation between $B_1, B_2, B_3$ that the only possibilities are $B_3 \ll B_1 \asymp B_2$ and $B_1 \leq B_2 \asymp B_3$. In the first case, \eqref{1b} gives $\text{meas}(\mathcal{A}_{j, k}) \ll T^{\varepsilon} B_1B_3^{1/2}$, so that 
\begin{equation}\label{bound1a}
\mathcal{I} \ll T^{-2+\varepsilon} B_3 \leq T^{-1 - \frac{b}{5} + \varepsilon}.
\end{equation}
In the second case, $\text{meas}(\mathcal{A}_{j, k}) \ll T^{\varepsilon} B_1^{1/2}B_3$, and we obtain the same bound. 

\emph{Case 1b:} 
It follows from the triangle inequality that the current assumptions imply $j \not= k$. Moreover $t_1 = r_j + o(T)$ and $t_2 = -r_k + o(T)$, so that  $P_2(t)$, defined in \eqref{p1}, satisfies $ P_2(t) = P_2(r_j, -r_k) + o(T^6) \asymp T^6$ by \eqref{factP1}. We conclude that $|\det g''(t) | \asymp (B_1B_2)^{-1}$, and hence by \eqref{33} (whose assumption is automatically satisfied) 
$$\text{meas}(\mathcal{A}_{j, k}) \ll T^{\varepsilon} (B_1B_2)^{1/2}$$
so that after inserting into \eqref{mainint} we obtain 
\begin{equation}\label{bound1b}
\mathcal{I} \ll \frac{B_3^{1/2}}{T^{2-\varepsilon}}\ll  T^{-3/2 + \varepsilon}.
\end{equation}

\emph{Case 1c:} This case follows closely with Case 2 of Section \ref{B1small}. 
If $\{j, l, m\} = \{1, 2, 3\}$ (again $k$ may or may not be different from $j$), we conclude from \eqref{factP2} and the present assumptions that 
$$\frac{\partial g_2}{\partial t_2}(t) = \frac{r_j(t_2 - 2r_lr_m/r_j)}{t_2(t_2+r_l)(t_2 + r_m)} + O(T^{-1 - \frac{3}{5}b}).$$
Let $0 < a < \frac{3}{5}b$. We distinguish two cases. If 
$$\Bigl| t_2 - 2\frac{r_lr_m}{r_j}\Bigr| \geq T^{1-a},$$
then $|\frac{\partial g_2}{\partial t_2}(t)| \gg  T^{-1-a} \gg E_2^{-5/4 + \varepsilon}$ (here we use $a < 3b/5$), so that \eqref{1a} with $M = \Bbb{R}$ implies
$$\text{meas}(\mathcal{A}_{j, k}) \ll T^{1+a+\varepsilon} B_2^{-1/2} B_1 \leq T^{1 + a - \frac{2}{5}b + \varepsilon} B_1^{1/2}$$
and hence
$$\mathcal{I} \ll \frac{\text{meas}(\mathcal{A}_{j, k}) }{T^{2-\varepsilon}} \left(\frac{B_3}{B_1B_2}\right)^{1/2} \ll T^{-1 + a - \frac{2}{5}b + \varepsilon} .$$
In the opposite case, we have
$$\frac{\partial g_1}{\partial t_1}(t) = \sum_{l=1}^3 \frac{1}{t_1 - r_l} - \frac{1}{t_1 + t_2} - \frac{2}{t_1}  = \frac{1}{t_1 - r_j}   +O\left(\frac{1}{T} + \frac{1}{B_3}\right) \asymp \frac{1}{B_1}.$$ 
By \eqref{2a} (whose assumption is trivially satisfied) with 
$$M = \left(2 \frac{r_l r_m}{r_j} - T^{1-a}, 2 \frac{r_l r_m}{r_j} + T^{1-a}\right)$$
we obtain 
$\text{meas}(\mathcal{A}_{j, k}) \ll T^{1 - a + \varepsilon} B_1^{1/2},$ 
so that 
$$\mathcal{I} \ll T^{-1 -a+ \varepsilon}. $$
Choosing $a = b/5$, we   obtain the final bound
\begin{equation}\label{bound1c}
\mathcal{I} \ll T^{-1  - \frac{1}{5}b + \varepsilon}
\end{equation}
in the present Case 1c. \\

Combining all three bounds \eqref{bound1a} -- \eqref{bound1c} from Cases 1a -- 1c, the final bound in Case  1 is $\mathcal{I} \ll T^{-1-\frac{b}{5} + \varepsilon}$ in agreement with \eqref{toshow}. The remaining 5 cases are simple.

\subsection{Case 2} The present assumption $|U_2| \asymp T \lll |U_1| \ll T^{1+b}$ implies $B_1, B_3, E_1 \asymp |U_1|$, $E_2 \asymp B_2 \ll T$. By excluding \eqref{case1a}, we may assume $B_2 \leq T^{1-b}$. Then we conclude from \eqref{modif1}  that $\text{meas}(\mathcal{A}_{j, k}) \ll T^{-1/2 + \varepsilon} B_2 |U_1|$, so that
$$\mathcal{I} \ll T^{\varepsilon} \frac{\text{meas}(\mathcal{A}_{j, k})}{|U_1| T B_2^{1/2} } \ll  T^{-3/2 + \varepsilon} B_2^{1/2} \leq T^{-1 - \frac{b}{2}  + \varepsilon}.$$

\subsection{Case 3}   The present assumption $|U_1| \lll T \asymp |U_2|$ implies  $E_1 \asymp |U_1|$, $E_2 \asymp B_2 \ll T \asymp B_3 \asymp B_1$.
We have $P_1(t) = 2 C_2 t_2+o(T^4) \asymp T^4$ so that $|\frac{\partial g_1}{\partial t_1}(t)| \asymp |U_1|^{-1} \gg T^\epsilon E_1^{-5/4}$.  By \eqref{2a} with $M=\RR$ we have  
$\text{meas}(\mathcal{A}_{j, k}) \ll T^{\varepsilon} |U_1|^{1/2} B_2, $ 
and hence
$$\mathcal{I} \ll \frac{\text{meas}(\mathcal{A}_{j, k})}{T^{2-\varepsilon} B_2^{1/2}} \ll T^{-2 + \varepsilon} (|U_1| B_2)^{1/2}. $$
By excluding \eqref{case1a}, we have $\min (|U_1|, B_2) < T^{1-b}$, giving
$$\mathcal{I} \ll T^{-1-\frac{b}{2} + \varepsilon}. $$

\subsection{Case 4} The present assumption $|U_1| \gg |U_2| \ggg T$ implies $B_2 \asymp |U_2|$, $B_1 \asymp |U_1|$. By excluding \eqref{case1a}, we may assume either $|U_1| \geq T^{1+b}$ or $B_3 \leq T^{1-b}$.\\

\emph{Case 4a:} Let us first assume that $B_3 \leq T^{1-b}$.  Then in particular $|U_1| \asymp |U_2|$ and $E_1 \asymp B_3$, and \eqref{1b} implies
$\text{meas}(\mathcal{A}_{j, k}) \ll T^{\varepsilon} |U_1| B_3^{1/2},$ 
so
$$\mathcal{I} \ll T^{\varepsilon} \frac{B_3^{1/2}}{|U_2|^{3/2} |U_1|^{3/2}} \text{meas}(\mathcal{A}_{j, k}) \ll T^{\varepsilon} |U_1|^{-2} B_3 \ll T^{-1 - b+\varepsilon}.$$

\emph{Case 4b:} Next we assume $|U_1| \geq T^{1+b}$, and in view of the results from Section \ref{sec7} we may also assume that $|U_2| \geq T^{1 + b/4}$. 
  Then, using $B_3 \ll |U_1|$, \eqref{2b} implies
$$\text{meas}(\mathcal{A}_{j, k}) \ll T^{\varepsilon} \frac{|U_2| B_3}{\min(|U_2|, B_3)^{1/2}} \ll T^{\varepsilon} |U_1| |U_2|^{1/2},$$
so
$$\mathcal{I} \ll T^{\varepsilon} |U_2|^{-1} \ll T^{-1 - \frac{b}{4}+\varepsilon}.$$


\subsection{Case 5} The present assumption $|U_1| \ggg T \ggg |U_2|$ 
implies $B_1 \asymp B_3 \asymp |U_1|$ and  $B_2  \asymp T$, $E_2 \asymp |U_2|$.  
Notice that the case $\abs{U_2} \ge T^{1-b/4}$ is covered by \eqref{case1a} and \eqref{case7} as $\abs{U_1}$ is respectively smaller or larger than $T^{1+b}$.
Hence we may assume $\abs{U_2} < T^{1-b/4}$, and we simply use \eqref{2b} to conclude $\text{meas}(\mathcal{A}_{j, k}) \ll T^{\varepsilon} |U_1| |U_2|^{1/2}$, so that 
$$\mathcal{I} \ll T^{\varepsilon} \frac{\text{meas}(\mathcal{A}_{j, k})}{|U_1|T^{3/2}} \ll T^{-1-\frac{b}{8} + \varepsilon}.$$

\subsection{Case 6} Under the present assumption $|U_1| \leq |U_2| \lll T$, and excluding the case \eqref{case1a}, we may assume $\min(|U_1|, B_3) \leq T^{1-b}$ since $B_1 \asymp B_2 \asymp T$.  Clearly, $B_3 \ll U_2$. \\

\emph{Case 6a:} Let us first assume that $B_3 \lll |U_2|$. Then in particular $|U_1| \asymp |U_2|$, so that $B_3 \leq T^{1-b}$ and $E_1 \asymp E_2 \asymp B_3$.
Now \eqref{1b} implies $\text{meas}(\mathcal{A}_{j, k}) \ll T^{\varepsilon} |U_1| B_3^{1/2}$, so that 
$$\mathcal{I} \ll T^{-3+\varepsilon} B_3^{1/2} \text{meas}(\mathcal{A}_{j, k}) \ll T^{-1-b + \varepsilon}.$$

\emph{Case 6b:} Next we assume $B_3 \asymp |U_2|$. In this case, $E_1 \asymp |U_1|$, $E_2 \asymp |U_2|$, and we must have $|U_1| \leq T^{1-b}$.
Again using \eqref{1b}, $\text{meas}(\mathcal{A}_{j, k}) \ll T^{\varepsilon} |U_1|^{1/2} B_3$, so that 
$$\mathcal{I} \ll T^{-3+\varepsilon} B_3^{1/2} \text{meas}(\mathcal{A}_{j, k}) \ll T^{-1-\frac{b}{2} + \varepsilon}.$$


\begin{thebibliography}{EMOT}

\bibitem[Bl1]{Bl1} V. Blomer, \emph{Subconvexity for  twisted $L$-functions on ${\rm GL}(3)$},  Amer. J. Math. \textbf{134} (2012), 1385-1421

\bibitem[Bl2]{Bl} V. Blomer, \emph{Applications of the Kuznetsov formula on ${\rm GL}(3)$}, Invent. math. \textbf{194} (2013), 673-729


\bibitem[BKY]{BKY} V. Blomer, R. Khan, M. Young, \emph{Mass distribution of holomorphic cusp forms}, Duke Math. J. \textbf{162} (2013), 2609-2644

\bibitem[Br]{Br} F. Brumley, \emph{Second order average estimates on local data of cusp forms}, Arch. Math. \textbf{87} (2006), 19-32

\bibitem[Bum]{Bum} D. Bump, \emph{Automorphic forms on ${\rm GL}(3, \Bbb{R})$}, Lecture Notes in Mathematics 1083, Springer-Verlag, 1984

\bibitem[BFG]{BFG} D. Bump, S. Friedberg, D. Goldfeld, \emph{Poincar\'e series and Kloosterman sums for ${\rm SL}(3, \Bbb{Z})$}, Acta Arith. \textbf{50} (1988), 31-89

\bibitem[Bu1]{Bu1} J. Buttcane, \emph{Sums of ${\rm SL}(3, \Bbb{Z})$ Kloosterman sums}, \changed{Ramanujan J. \textbf{32} (2013), 371-419}

\bibitem[Bu2]{Bu} J. Buttcane, \emph{The spectral Kuznetsov formula on ${\rm SL}(3, \Bbb{Z})$}, {\tt arXiv:1411.7802}


\bibitem[EMOT]{EMOT} A. Erd\'elyi, W. Magnus, F. Oberhettinger, F. Tricomi, \emph{Higher transcendental functions II}, McGraw-Hill 1953

\bibitem[Go]{Go} D. Goldfeld, \emph{Automorphic forms and L-functions for the group $ {\rm GL}(n, \Bbb{R})$}, Cambridge studies in advanced mathematics \textbf{99} (2006)

\bibitem[GR]{GR} I. S. Gradshteyn, I. M. Ryzhik, \emph{Tables of integrals, series, and products}, 7th edition, Academic Press, New York, 2007.




\bibitem[Iw]{Iw} H. Iwaniec, \emph{The spectral growth of automorphic L-functions}, J. reine angew. Math. \textbf{428} (1992), 139-159

\bibitem[IK]{IK} H. Iwaniec, E. Kowalski, \emph{Analytic number theory}, Colloq. Publ. \textbf{53}, Amer. Math. Soc., Providence, RI, 2004.

\bibitem[Li1]{Li1} Xiannan Li, \emph{Upper bounds on $L$-functions at the edge of the critical strip}, IMRN 2010, 727-755

\bibitem[Li2]{Li} Xiaoqing Li, \emph{Bounds for ${\rm GL}(3) \times  {\rm GL}(2)$ $L$-functions and ${\rm GL}(3)$ $L$-functions},  Annals of Math. \textbf{173} (2011), 301-336

\bibitem[MV]{MV} P. Michel, A. Venkatesh, \emph{The subconvexity problem for ${\rm GL}_2$}, Publ. Math. IHES \textbf{111} (2010), 171-271

\bibitem[Mi]{Mi} J. Milnor, \emph{On the Betti numbers of real varieties},  Proc. Amer. Math. Soc. \textbf{15} (1964),  275-280

\bibitem[Mu1]{Mu1} R. Munshi, \emph{The circle method and bounds for $L$-functions - II: Subconvexity for twists of ${\rm GL}(3)$ $L$-functions},  to appear Amer.  J. Math. 

\bibitem[Mu2]{Mu2} R. Munshi, \emph{The circle method and bounds for $L$-functions - III: $t$-aspect subconvexity for ${\rm GL}(3)$ $L$-functions}, {\tt arXiv:1301.1007}

\bibitem[Mu3]{Mu3} R. Munshi,  \emph{The circle method and bounds for $L$-functions  - IV: Subconvexity for twists of ${\rm GL}(3)$ $L$-functions}, to appear in Annals of Math. 

\bibitem[Th]{Th} R. Thom, \emph{Sur l'homologie des vari\'et\'es alg\'ebriques r\'eelles}, in: Differential and Combinatorial Topology (A Symposium in Honor of Marston Morse) pp.\ 255-265 (1965), Princeton Univ. Press, Princeton, N.J.


\end{thebibliography}
\end{document}